\documentclass[reqno,10pt,centertags,draft]{amsart}
\usepackage{amsmath,amsthm,amscd,amssymb,latexsym,upref}
\date{\today}

\newcommand{\bbN}{{\mathbb{N}}}
\newcommand{\bbR}{{\mathbb{R}}}

\newcommand{\bbZ}{{\mathbb{Z}}}
\newcommand{\bbC}{{\mathbb{C}}}

\newcommand{\cA}{{\mathcal A}}
\newcommand{\cB}{{\mathcal B}}
\newcommand{\cC}{{\mathcal C}}

\newcommand{\cE}{{\mathcal E}}
\newcommand{\cF}{{\mathcal F}}

\newcommand{\cH}{{\mathcal H}}

\newcommand{\cJ}{{\mathcal J}}
\newcommand{\cK}{{\mathcal K}}
\newcommand{\cL}{{\mathcal L}}

\newcommand{\cP}{{\mathcal P}}
\newcommand{\cQ}{{\mathcal Q}}

\newcommand{\cU}{{\mathcal U}}
\newcommand{\cV}{{\mathcal V}}
\newcommand{\cW}{{\mathcal W}}
\newcommand{\cY}{{\mathcal Y}}

\newcommand{\dott}{\,\cdot\,}
\newcommand{\no}{\notag}
\newcommand{\lb}{\label}
\newcommand{\f}{\frac}

\newcommand{\ol}{\overline}

\newcommand{\wti}{\widetilde}

\newcommand{\Oh}{O}

\newcommand{\loc}{\text{\rm{loc}}}

\newcommand{\ran}{\text{\rm{ran}}}

\newcommand{\dom}{\text{\rm{dom}}}

\newcommand{\supp}{\text{\rm{supp}}}

\newcommand{\bi}{\bibitem}
\newcommand{\hatt}{\widehat}
\newcommand{\beq}{\begin{equation}}
\newcommand{\eeq}{\end{equation}}
\newcommand{\ba}{\begin{align}}
\newcommand{\ea}{\end{align}}

\newcommand{\tr}{\text{\rm{tr}}}

\newcommand{\Wr}{\text{\rm{Wr}}}

\renewcommand{\Re}{\text{\rm Re}}
\renewcommand{\Im}{\text{\rm Im}}

\renewcommand{\ge}{\geqslant}
\renewcommand{\le}{\leqslant}

\DeclareMathOperator{\sgn}{sgn}

\DeclareMathOperator{\diag}{diag}

\allowdisplaybreaks
\numberwithin{equation}{section}

\newtheorem{theorem}{Theorem}[section]

\newtheorem{lemma}[theorem]{Lemma}

\newtheorem{hypothesis}[theorem]{Hypothesis}
\newtheorem{definition}[theorem]{Definition}
\newtheorem{example}[theorem]{Example}
\theoremstyle{definition}
\newtheorem{remark}[theorem]{Remark}

\begin{document}

\title[Derivatives of (Modified) Fredholm Determinants]
{Derivatives of (Modified) Fredholm Determinants and Stability of
Standing and Traveling Waves}
\author[F.\ Gesztesy, Yu.\ Latushkin, and
K.\ Zumbrun]{Fritz Gesztesy, Yuri\ Latushkin, and Kevin Zumbrun \hspace{80pt}}
\address{Department of Mathematics,
University of Missouri, Columbia, MO 65211, USA}
\email{fritz@math.missouri.edu}
\urladdr{http://www.math.missouri.edu/personnel/faculty/gesztesyf.html}
\address{Department of Mathematics,
University of Missouri, Columbia, MO 65211, USA}
\email{yuri@math.missouri.edu}
\urladdr{http://www.math.missouri.edu/personnel/faculty/latushkiny.html}
\address{Mathematics Department, Indiana University,
Bloomington, IN 47405, USA}
\email{kzumbrun@indiana.edu}
\urladdr{http://www.math.indiana.edu/people/profile.phtml?id=kzumbrun}
\thanks{Partially supported by the US National Science
Foundation Grant Nos.\ DMS-0405526, DMS-0338743, DMS-0354339,
DMS-0070765, DMS-0300487, and by the CRDF grant UP1-2567-OD-03.}
\thanks{{\it J. Math. Pures Appliqu\'ees} {\bf 90}, 160--200 (2008).}
\date{\today}

\begin{abstract}
Continuing a line of investigation initiated in \cite{GLM07}
exploring the connections between Jost and Evans functions and
(modified) Fredholm determinants of Birman--Schwinger type integral
operators, we here examine the stability index, or sign of the
first nonvanishing derivative at frequency zero of the
characteristic determinant, an object that has found considerable
use in the study by Evans function techniques of stability of
standing and traveling wave solutions of partial differential
equations (PDE) in one dimension. This leads us to the derivation
of general perturbation expansions for analytically-varying
modified Fredholm determinants of abstract operators. Our main conclusion, similarly in
the analysis of the determinant itself, is that the derivative of
the characteristic Fredholm determinant may be efficiently
computed from first principles for integral operators with
semi-separable integral kernels, which
include in particular the general one-dimensional case, and for
sums thereof, which appears to offer
applications in the multi-dimensional case.

A second main result is to show that the multi-dimensional 
characteristic Fredholm determinant is the renormalized limit 
of a sequence of Evans functions
defined in \cite{LPSS00} on successive Galerkin subspaces,
giving a natural extension of the one-dimensional results
of \cite{GLM07} and
answering a question of \cite{N07} whether this sequence 
might possibly converge (in general, no, but with renormalization, yes).
Convergence is useful in practice for numerical error control
and acceleration.
\end{abstract}
\maketitle

{\bf R\'esum\'e:}

Nous poursuivons l'\'etude, initi\'ee dans \cite{GLM07}, des liens  entre  les fonctions de
Jost et d'Evans et les d\'eterminants
(modifi\'es) de Fredholm d'op\'erateurs int\'egraux de type Birman--Schwinger.
Nous examinons ici l'indice de stabilit\'e, c'est-\`a-dire le signe de   la premi\`ere
d\'eriv\'ee non nulle, \`a la fr\'equence z\'ero, du d\'eterminant
caract\'eristique. Cet indice a trouv\'e une utilisation   consid\'erable dans l'\'etude,
par des techniques de fonction d'Evans,  de la stabilit\'e de
solutions de type ondes progressives de syst\`emes d'\'equations aux   d\'eriv\'ees partielles
en une dimension d'espace. Cela nous am\`ene \`a \'ecrire des 
formules pour des d\'eveloppements g\'en\'eraux de type perturbatif pour
les d\'eterminants modifi\'es de Fredholm d'op\'erateurs analytiques   abstraits.
Notre conclusion principale
est que la d\'eriv\'ee du d\'eterminant caract\'eristique de Fredholm,
comme le d\'eterminant lui-m\^ eme, 
peut \^etre calcul\'ee   efficacement
pour des op\'erateurs int\'egraux  dont les
noyaux sont semi-s\'eparables et pour les sommes en. 
Le premier classe d'op\'erateurs inclut en   particulier le
cas g\'en\'eral en une dimension d'espace;
le dernier laisse envisager des applications
au cas multidimensionnel.

Le deuxi\`eme r\'esultat principal est la preuve que le
d\'eterminant caract\'eristique multidimensionnel de Fredholm est
la limite renormalis\'ee d'une suite de fonctions d'Evans,
d\'efinie dans \cite{LPSS00}, sur des sous-espaces successifs de
Galerkin; ce r\'esultat est une extension naturelle des r\'esultats   unidimensionnels
de \cite{GLM07}, et r\'epond \`a la question, pos\'ee dans \cite{N07},  de la
convergence de cette suite (la r\'eponse est qu'en g\'en\'eral  il n'y a pas
convergence, mais qu'on peut obtenir la convergence  apr\`es  renormalisation). 
Convergence
est utile dans la pratique
pour le contr\^ole d'erreur et l'acc\'el\'eration des calculs num\'eriques.

\section{Introduction}\label{s1}

A problem of general interest is to determine the spectrum
of a general variable-coefficient
linear differential operator $L=\sum_{|\alpha|=0}^N
a_\alpha(x) \partial_{x}^{\alpha}$,
 $a_\alpha(x)\in \bbR^{n\times n}$, $x\in \bbR^d$,
with prescribed behavior of the coefficients as $|x|\to \infty$.
This arises naturally, for example, in the study of traveling-
or standing-wave solutions of nonlinear PDEs in a wide variety
of applications, as described, for example, in the survey articles
\cite{Sa02}, \cite{Zu03}, and the references therein.
In the one-dimensional case, $d=1$, a very useful and general
tool for this purpose is the {\it Evans function}
\cite{AGJ90}, \cite{Ev72}--\cite{Ev75}, \cite{PW92},
defined as a Wronskian $\cE(z)$
of bases of the set of solutions $\Psi_\pm$ of the associated eigenvalue
ODE $(L-\lambda)\Psi=0$ decaying at $x=+\infty$ and $x=-\infty$,
respectively, whose zeros correspond in location and multiplicity
with the eigenvalues of $L$.

Among the many applications of the Evans function, perhaps the
simplest and most general is the computation of the {\it stability
index}
\begin{equation}
\Gamma= \sgn\big(d_z^k \cE(0)\big) \sgn (\cE(+\infty)),
\end{equation}
whose sign determines the parity of the number of unstable
eigenvalues $\lambda$, or eigenvalues with positive real part
$\Re (\lambda) >0$, where $d_z^k \cE(0)$ is the first nonvanishing
derivative of $\cE(z)$ at $z=0$;
see, for example, \cite{Ev75}, \cite{PW92}, \cite{Zu03}.
(A standard property of the Evans function
is that it may be constructed so as to respect complex conjugation;
in particular, it may be taken real-valued for $z\in \bbR$.)
A problem that has received considerable recent interest\footnote{
For example, this was a focus topic of the workshop
``Stability Criteria for Multi-Dimensional
Waves and Patterns'', at the American Institute of
Mathematics (AIM) in Palo Alto (California/USA), May 16-20, 2005.}
is to extend the Evans function, and in particular the stability index,
to the more general setting of multi-dimensions in a way that
is useful for practical computations.
Here, we refer mainly to numerical computation, as presumably the
only feasible way to treat large-scale problems associated
with multi-dimensions.

Various different constructions have been suggested toward this end;
see, for example, \cite{DN06}, \cite{DN08}, \cite{LPSS00}. 
However, only one of these, the Galerkin approximation method
of \cite{LPSS00} (described in Section \ref{s4}), seems in principle computable,
and the computations involved appear quite numerically intensive.
(So far, no such computations have satisfactorily been carried out,
though, see the proposed methods discussed in \cite{HZ06}, \cite{N07}.)
It is therefore highly desirable to explore other directions that
may be more computationally efficient.

Here, we follow a very natural direction first proposed in \cite{GLM07}.
Specifically, it is shown in \cite{GLM07} for a quite
general class of one-dimensional operators $L$ that
the Evans function, appropriately normalized,
agrees with a (modified) characteristic {\it Fredholm determinant}, thus
generalizing the classical relation known for the {\it Jost function} associated
with Schr\"odinger operators of mathematical physics; for further discussion
of this problem and its history, see \cite{GLM07} and the references therein.
A central point of the analysis is the observation
\cite{GM04} that {\it semi-separability} of the integral
kernel of the resolvent operator is the key property of
one-dimensional operators that makes possible the reduction of
the infinite-dimensional characteristic determinant to a 
finite-dimensional determinant expressed by the Jost or Evans function.

The identification of Evans functions and Fredholm determinants
yields a natural generalization of the Evans function to multi-dimensions,
since the definition of (modified) Fredholm determinants extends to 
higher-dimensional problems.
What is not immediately clear is whether this extension leads
to a practically useful, computable formulation of either the Evans function
or the stability index.
For, up to now, the main approach to computation of the characteristic
determinant was to express it as a Jost function or the
usual Wronskian expression for the (one-dimensional) Evans function.

In the present paper, we extend some of the investigations of \cite{GLM07}
in two ways.
First, we derive in Theorem \ref{t2.7} a general perturbation formula
for analytically-varying modified Fredholm determinants,
by which we may express the stability index
as a product of a finite-dimensional minor (Lyapunov--Schmidt
decomposition) and a finite-rank perturbation of the original
characteristic determinant: that is, directly in terms of
Fredholm determinants, without reference to a Jost or Evans
function formulation. This is carried out in Section \ref{s2},
with the main result given in Theorems \ref{t2.3} and \ref{t2.7}.

Second, we discuss in Theorem \ref{t3.6a} an important case when the finite-rank perturbation part, like the original
characteristic determinant, may be reduced to a finite-dimensional
determinant whenever the resolvent of the operator $L$ has a semi-separable
integral kernel, in particular, in the one-dimensional case.
This is a consequence of the simple observation that a
sum of operators with semi-separable integral kernels may be expressed as
an operator with a matrix-valued semi-separable integral kernel and evaluated
in the same way; indeed, the analysis of \cite{GM04} on which this reduction is based
is actually presented in the more general,
matrix-valued setting. 
We illustrate this procedure in Section \ref{s3} by explicit computations
for the example of a scalar Schr\"odinger operator that arises in
the study of stability of standing-wave patterns of
one-dimensional reaction--diffusion equations,
in the process illuminating various relations between Jost functions
and characteristic Fredholm determinants.

To explain our main results in Section \ref{s3}, we recall that a classical formula
by Jost and Pais equates the Jost function and the Fredholm determinant
of a Birman--Schwinger-type operator (cf.\ Section \ref{s3} for
its discussion and definitions).  Since the Jost function is the
Wronskian of the Jost solutions  $\Psi_{\pm}$, this result can
be viewed as a calculation of the Fredholm determinant via the
solutions $\Psi_{\pm}$ of the homogenous Schr\"odinger equation
that are asymptotic to the exponential plane waves.  We 
prove a
new formula in this spirit, see \eqref{3.71}, and 
compute the
{\em derivative} of the Jost function, that is, the derivative of the Fredholm determinant, via some solutions $\psi_{\pm}$ of a {\em nonhomogenous} Schr\"odinger equation
(cf.\ \eqref{3.98a}) that, in turn, are asymptotic to $\Psi_{\pm}$ (cf.\ \eqref{3.72}). 
We are not aware of any earlier references mentioning these solutions  $\psi_{\pm}$.

In addition, we obtain in Section \ref{s3} in passing an elementary proof of
an interesting formula derived by Simon \cite{Si00}
for the Jost solutions $\Psi_\pm(z,\dott)$  in terms of Fredholm determinants.

Of course, the approach of Section \ref{s3}
applies equally well to the general one-dimensional case,
yielding in principle a similarly compact formula for the stability
index obtained entirely through Fredholm determinant manipulations.
However, in this paper we do not pursue this any further, 
leaving this topic for future research.

More generally, this suggests an approach to the multi-dimensional case
by a limiting procedure based on Galerkin approximation, but
carried out within the Fredholm determinant framework.
More precisely, we propose in place of the standard approach of reducing to
a large one-dimensional system by Galerkin approximation then defining
a standard Evans function, to first relate the Evans function and a
Fredholm determinant, then evaluate the latter by Galerkin 
approximation/semi-separable reduction.
This offers the advantage that successive levels of approximation
are embedded in a hierarchy of convergent problems useful for
error control, without the need to prescribe appropriate normalizations
by hand.

Our main result in the multi-dimensional case,
generalizing the one-dimensional results of \cite{GLM07},
is that the sequence of approximate $2$-modified Jost
functions $\cF_{2,J}$
generated by Galerkin approximation of the $2$-modified
Fredholm determinant
at wave number $J$ {\it agrees, up to appropriate
normalization}, with the Evans functions $\cE_J$ constructed
in \cite{LPSS00} for the sequence of one-dimensional equations
obtained by Galerkin approximation at the same wave number;
see Theorem \ref{equivalence}.
This includes the information that the sequence $\cE_J$, 
introduced in \cite{LPSS00} as a tool to compute a topologically-defined
stability index, in fact determines a well-defined
$2$-modified Jost function; that is,
an appropriate renormalization $\cF_{2,J}= e^{\Theta_J} \cE_J$
of the sequence $\cE_J$ converges to a limiting $2$-modified Jost
 function $\cF_{2}$ (cf.\ Theorem \ref{convergence}),
a fact that is not apparent from the construction of \cite{LPSS00}.
The Jost and Evans functions of course carry considerably more information
than the stability index alone.

We obtain at the same time a slightly different algorithm for computing
the stability index, which avoids some logistical
difficulties of the existing Galerkin schemes.
We discuss these issues in Section \ref{s4}, illustrating our approach with respect
to the basic multi-dimensional examples of
flow in an infinite cylinder and solutions with radial limits. 

Finally, we note that the ODE systems arising in computation of
the Galerkin-based Evans functions $\cE_J$ become extremely stiff
as $J\to \infty$, featuring growth/decay modes of order $\pm J$.
Thus, numerical conditioning becomes a crucial consideration for the
large-scale systems that result in multi-dimensions
($J\sim 100$, as described in \cite{HZ06}).
It may well be that, for sufficiently large $J$,
direct computation of the Fredholm determinant by discretization
may be more efficient than either Evans function computations
or simple discretization of the linearized operator $L$;
see Subsection \ref{conditioning}. 
Efficient numerical realization of this approach would be
a very interesting direction for future investigation.

{\bf Plan of the paper.}
In Section \ref{s2}, we provide a general perturbation
formula for analytically varying (modified) Fredholm determinants.
In Section \ref{s3}, we use this result together with the
reduction method of \cite{GM04}
to compute the stability index in the case of a one-dimensional
self-adjoint Schr\"odinger operator.
Finally, in Section \ref{s4}, we describe
extensions to multi-dimensions.

\section{A general perturbation expansion for Fredholm Determinants}
\label{s2}

In this section we describe the analytic behavior of Fredholm determinants
${\det}_{\cH}(I-A(z))$ and modified Fredholm determinants
${\det}_{2,\cH}(I-A(z))$ in a
neighborhood of $z=0$ with $A(\dott)$ analytic in a neighborhood of $z=0$ in
trace norm, respectively, Hilbert--Schmidt norm. Special emphasis
will be put on
the case where $[I-A(0)]$ is not boundedly invertible in the Hilbert
space $\cH$.

In the first part of this section we suppose that all relevant
operators belong to
the trace class and consider the associated Fredholm determinants. In
the second
part we consider $2$-modified Fredholm determinants in the case where the
relevant operators are Hilbert--Schmidt operators.

\subsection{Trace class operators}
In the course of the proof of our first result we repeatedly will
have to use some
of the standard properties of determinants, such as,
\begin{align}
& {\det}_{\cH}((I_\cH-A)(I_\cH-B))={\det}_{\cH}(I_\cH-A) \,
{\det}_{\cH}(I_\cH-B),
\quad A, B \in\cB_1(\cH), \lb{2.1} \\
& {\det}_{\cH'}(I_{\cH'}-AB)={\det}_{\cH}(I_\cH-BA) \;\, \text{ for all
$A\in\cB(\cH,\cH')$, $B\in\cB(\cH',\cH)$} \lb{2.2} \\
& \hspace*{5.2cm} \text{ such that $BA\in \cB_1(\cH)$,
$AB\in \cB_1(\cH')$,} \no
\intertext{and}
&{\det}_{\cH}(I_\cH-A)={\det}_{\bbC^k}(I_k-D_k) \, \text{ for } \,
A=\begin{pmatrix}
0 & C \\ 0 & D_k \end{pmatrix}, \;\, \cH=\cK\dotplus \bbC^{k}, \lb{2.3}
\intertext{since}
&I_\cH -A=\begin{pmatrix} I_\cK & -C \\ 0 & I_k-D_k \end{pmatrix} =
\begin{pmatrix} I_\cK & 0 \\ 0 & I_k-D_k \end{pmatrix}
\begin{pmatrix} I_\cK & -C \\ 0 & I_k \end{pmatrix}. \lb{2.4}
\end{align}
Finally, assuming $A, B \in\cB_1(\cH)$, we also mention the following estimates:
\begin{align}
& |{\det}_{\cH}(I_{\cH}-A)| \leq \exp(\|A\|_{\cB_1(\cH)}),   \lb{2.4a} \\
& |{\det}_{\cH}(I_{\cH}-A) - {\det}_{\cH}(I_{\cH}-B)| \leq \|A-B\|_{\cB_1(\cH)}
\lb{2.4b} \\
& \hspace*{5.3cm}
\times \exp(\|A\|_{\cB_1(\cH)}+\|B\|_{\cB_1(\cH)}+1).  \no
\end{align}
Here $\cH$ and $\cH'$ are complex separable Hilbert spaces,
$\cB(\cH)$ denotes the set of bounded linear operators on $\cH$, $\cB_p(\cH)$,
$p\geq 1$, denote the usual trace ideals of $\cB(\cH)$, and $I_\cH$ denotes
the identity operator in $\cH$ (similarly, $I_k$ abbreviates the
identity operator in $\bbC^k$). The ideal of compact operators on
$\cH$ will be denoted by $\cB_{\infty}(\cH)$. Moreover,
${\det}_{p,\cH}(I_\cH-A)$,
$A\in\cB_p(\cH)$, denotes the ($p$-modified) Fredholm determinant of $I_\cH-A$
with ${\det}_{1,\cH}(I_\cH-A)={\det}_{\cH}(I_\cH-A)$,
$A\in\cB_1(\cH)$, the standard Fredholm determinant of a trace class
operator, and $\tr_{\cH}(A)$, $A\in\cB_1(\cH)$, the trace of a trace
class operator in $\cH$. Finally, $\dotplus$ in \eqref{2.3} denotes a
direct but not necessary orthogonal direct decomposition of $\cH$ into $\cK$
and the $k$-dimensional subspace $\bbC^k$. These results can be found, for
instance, in \cite{GGK96}, \cite[Sect.\ IV.1]{GK69}, \cite[Ch.\ 17]{RS78},
\cite{Si77},  \cite[Ch.\ 3]{Si05}. In the following, $\sigma(T)$
denotes the spectrum of a densely defined, closed linear operator $T$
in $\cH$, and $\sigma_{\rm d}(T)$ denotes the discrete spectrum of
$T$ (i.e., isolated eigenvalues of $T$  of finite algebraic multiplicity).

For the general theory of (modified) Fredholm determinants we refer,
for instance, to
\cite[Sect.\ XI.9]{DS88}, \cite{GGK96}, \cite{GGK97}, \cite[Ch.\ X.III]{GGK00},
\cite[Ch.\ IV]{GK69}, \cite{Si77}, and \cite[Sects.\ 3, 9]{Si05}.

\begin{hypothesis}  \lb{h2.1}
Suppose $A(\dott)\in \cB_1(\cH)$ is a family of trace class operators on
$\cH$ analytic on an open neighborhood $\Omega_0\subset\bbC$ of $z=0$ in trace
class norm $\|\dott\|_{\cB_1(\cH)}$.
\end{hypothesis}

Given Hypothesis \ref{h2.1} we write
\begin{equation}
A(z) \underset{z\to 0}{=} A_0+A_1 z+\Oh\big(z^2\big) \, \text{ for
$z\in\Omega_0$ sufficiently small, }
\, A_\ell \in \cB_1(\cH), \; \ell=0,1.  \lb{2.5}
\end{equation}

We start by noting the following well-known result.

\begin{lemma}  \lb{l2.2}
Assume Hypothesis \ref{h2.1} and suppose $(I_{\cH}-A_0)^{-1}\in\cB(\cH)$. Then,
\begin{equation}
{\det}_{\cH}(I_\cH-A(z)) \underset{z\to 0}{=} {\det}_{\cH}(I_{\cH}-A_0)
-{\det}_{\cH}(I_{\cH}-A_0) \, \tr_{\cH}\big((I_{\cH}-A_0)^{-1}A_1\big) z
+ \Oh\big(z^2\big).  \lb{2.6}
\end{equation}
\end{lemma}
\begin{proof}
This follows from
\begin{align}
\begin{split}
{\det}_{\cH}(I_{\cH}-A(z))& \underset{z\to 0}{=}
{\det}_{\cH}\big((I_{\cH}-A_0)\big[I_{\cH}-(I_{\cH}-A_0)^{-1}A_1 z
+\Oh\big(z^2\big)\big]\big)  \\
&  \underset{z\to 0}{=}
{\det}_{\cH}(I_{\cH}-A_0)\big[1-\tr_{\cH}\big((I_{\cH}-A_0)^{-1}A_1\big) z
+ \Oh\big(z^2\big)\big], \lb{2.87}
\end{split}
\end{align}
where we used the fact that
\begin{equation}
{\det}_{\cH}(I_{\cH}-B z)=\exp\bigg[-\sum_{k\in\bbN}
\f{\tr_{\cH}\big(B^k\big)}{k}z^k\bigg] \,
\text{ for $|z|$ sufficiently small}  \lb{2.8}
\end{equation}
with $B\in\cB_1(\cH)$.
\end{proof}

Next we turn to the case where $(I_{\cH}-A_0)$ is not boundedly
invertible. Before
we state the analog of Lemma \ref{l2.2} in this more general setting,
we need some preparations.

We temporarily assume
\begin{equation}
A_0 \in \cB(\cH) \, \text{ and $1\in\sigma_{\rm d}(A_0)$,}   \lb{2.13b}
\end{equation}
and abbreviate by $P_0$ the Riesz projection associated with $A_0$
and the discrete eigenvalue $1$ of $A_0$,
\begin{equation}
P_0 =\f{-1}{2\pi i} \oint_{\cC_0} d\zeta \, (A_0 - \zeta
I_{\cH})^{-1},    \lb{2.13c}
\end{equation}
where $\cC_0$ denotes a sufficiently small counterclockwise oriented
circle centered at $1$ such that no part of
$\sigma(A_0)\backslash\{1\}$ intersects $\cC_0$ and its open
interior. We denote by
\begin{equation}
n_0=\dim(\ran(P_0))    \lb{2.14}
\end{equation}
the algebraic multiplicity of the eigenvalue $1$ of $A_0$. In
addition, we introduce the
quasinilpotent operator $D_0$ associated with $A_0$ and its discrete
eigenvalue $1$ by
\begin{equation}
D_0=(A_0-I_{\cH})P_0   \lb{2.14a}
\end{equation}
such that
\begin{equation}
D_0=D_0P_0=P_0D_0=P_0D_0P_0.   \lb{2.14b}
\end{equation}

In the following we denote by
$\Delta(z_0;r_0)\subset\bbC$ the open disc centered at $z_0\in\bbC$
of radius $r_0>0$ and by $\cC(z_0;r_0)=\partial \Delta(z_0;r_0)$ the
counterclockwise oriented circle of radius $r_0>0$ centered at $z_0$.
Assuming that the analytic family of operators
$A(\dott)$ satisfies Hypothesis \ref{h2.1}, we prescribe an
$\varepsilon_0 >0$ and choose a sufficiently small open neighborhood
$\Omega_0$ of $z=0$ such that all eigenvalues $\lambda_j(z)$ of
$A(z)$ for $z\in\Omega_0$, which satisfy
$\lambda_j(0)=1$, $1\leq j \leq \nu_0$ for some $\nu_0\in\bbN$ with
$\nu_0\leq n_0$, stay in the disc $\Delta(1;\varepsilon_0/2)$.
Moreover we assume that $\Omega_0$ is chosen sufficiently small that
no other eigenvalue branches of $A(z)$, $z\in\Omega_0$, intersect the
larger disc $\Delta(1;\varepsilon_0)$. Introducing the Riesz
projection
$P(z)$ associated with $A(z)$, $z\in\Omega_0$ (cf., e.g, \cite[Sect.\
III.6]{Ka80}),
\begin{equation}
P(z)=\f{-1}{2\pi i} \oint_{\cC(1;\varepsilon_0)} d\zeta \, (A(z) -
\zeta I_{\cH})^{-1}, \quad
z\in\Omega_0,   \lb{2.15}
\end{equation}
then $P(\dott)$ is analytic in $\Omega_0$ and we expand
\begin{equation}
P(z) \underset{z\to 0}{=} P_0 + P_1 z + \Oh\big(z^2\big) \, \text{
for $|z|$ sufficiently small.}
\lb{2.16}
\end{equation}
Moreover, we introduce the projections
\begin{equation}
Q(z)=I_{\cH} - P(z), \quad z\in\Omega_0,  \quad Q_0=I_{\cH} - P_0,  \lb{2.17}
\end{equation}
and expand
\begin{equation}
\quad Q(z) \underset{z\to 0}{=} Q_0 + Q_1 z + \Oh\big(z^2\big) \,
\text{ for $|z|$ sufficiently small.}   \lb{2.18}
\end{equation}
Since $P(z)^2=P(z)$, \eqref{2.16} implies
\begin{equation}
P_0P_1+P_1P_0=P_1 \, \text{ and hence } \, P_0P_1P_0=0.    \lb{2.19}
\end{equation}
Following Wolf \cite{Wo52} we now introduce the transformation
\begin{equation}
T(z)=P_0P(z)+Q_0Q(z) = P_0P(z)+[I_{\cH}-P_0][I_{\cH}-P(z)], \quad z\in\Omega_0,
\lb{2.20}
\end{equation}
such that
\begin{equation}
P_0T(z)=T(z)P(z), \quad Q_0T(z)=T(z)Q(z), \quad z\in\Omega_0.   \lb{2.21}
\end{equation}
In addition, for $|z|$ sufficiently small,
\begin{align}
T(z)& \underset{z\to 0}{=} I_{\cH} + (P_0P_1-P_1P_0)z +
\Oh\big(z^2\big),   \lb{2.22} \\
T(z)^{-1}& \underset{z\to 0}{=} I_{\cH} - (P_0P_1-P_1P_0)z +
\Oh\big(z^2\big),   \lb{2.23}
\end{align}
and hence
\begin{equation}
P_0=T(z)P(z)T(z)^{-1}, \quad Q_0=T(z)Q(z)T(z)^{-1}   \lb{2.24}
\end{equation}
for $|z|$ sufficiently small. Below we will use \eqref{2.24} to
reduce determinants in the Hilbert space $P(z)\cH$ to that in the
fixed Hilbert space $P_0\cH$.

Next, we introduce the following notation: We denote by
$S\in\cB\big(P_0\cH,\bbC^{n_0}\big)$ the boundedly invertible linear
operator which puts the nilpotent operator $P_0D_0P_0$ into its
$n_0\times n_0$ Jordan canonical form $SP_0D_0P_0S^{-1}$, and
abbreviate by $\nu_0$ the number of entries $1$ in the canonical
Jordan representation of $SP_0D_0P_0S^{-1}$, where\footnote{In
particular, $\nu_0$ equals the sum of the dimensions of all
nontrivial (i.e., nondiagonal) Jordan blocks in the canonical Jordan
representation $SP_0D_0P_0S^{-1}$ of $P_0D_0P_0$. Thus,
$SP_0D_0P_0S^{-1}$ contains $\nu_0$ entries $1$ at certain places
right above the main diagonal and $0$'s everywhere else including on
the main diagonal. In particular, $n_0 -\nu_0$ represents the geometric multiplicity 
of the eigenvalue $1$ of $A_0$.} $0\leq \nu_0\leq n_0-1$. Moreover, we denote by
$\wti A_1$ the $(n_0-\nu_0)\times (n_0-\nu_0)$-matrix obtained from
$SP_0A_1P_0S^{-1}$ by striking from it the $\nu_0$ columns and rows
in which $SP_0D_0P_0S^{-1}$ contains an entry $1$. With this notation
in mind, we now formulate our first abstract result on expansions of
Fredholm determinants
${\det}_{\cH}(I_\cH-A(z))$ as $z\to 0$, with $[I_\cH-A(0)]$ not
boundedly invertible in
$\cH$:

\begin{theorem}  \lb{t2.3}
Assume Hypothesis \ref{h2.1} and let
$1\in\sigma_{\rm d}(A_0)$. Then, given the notation in the paragraph
preceding this theorem,
\begin{align}
&{\det}_{\cH}(I_\cH-A(z))  \underset{z\to 0}{=}
\big[{\det}_{Q_0\cH}(I_{Q_0\cH}-Q_0A_0Q_0)+\Oh(z)\big]
\no \\
& \hspace*{3cm}
\times (-1)^{n_0} {\det}_{P_0\cH}(P_0D_0P_0+P_0[A_1+\Oh(z)]P_0 z)  \no \\
&\quad \underset{z\to 0}{=} {\det}_{Q_0\cH}(I_{Q_0\cH}-Q_0A_0Q_0)
{\det}_{\bbC^{n_0-\nu_0}}\big(\wti A_1\big)
(-z)^{n_0-\nu_0} + \Oh\big(z^{n_0-\nu_0 +1}\big).        \lb{2.25}
\end{align}
Here $n_0 -\nu_0$ represents the geometric multiplicity 
of the eigenvalue $1$ of $A_0$. 
In the special case where $1$ is a semisimple eigenvalue of $A_0$
$($i.e., where
$D_0=0$ and $\nu_0=0$$)$ one obtains,
\begin{align}
{\det}_{\cH}(I_\cH-A(z)) & \underset{z\to 0}{=}
{\det}_{Q_0\cH}(I_{Q_0\cH}-Q_0A_0Q_0)
{\det}_{P_0\cH}(P_0A_1P_0) (-z)^{n_0} + \Oh\big(z^{n_0+1}\big)  \no \\
& \underset{z\to 0}{=}  {\det}_{\cH}(I_{\cH}-P_0-A_0)
{\det}_{P_0\cH}(P_0A_1P_0) z^{n_0} + \Oh\big(z^{n_0+1}\big).  \lb{2.26}
\end{align}
In particular, if $1$ is a simple eigenvalue of $A_0$ $($i.e., if
$n_0=1$, $D_0=0$, and
$\nu_0=0$$)$, one obtains
\begin{equation}
{\det}_{\cH}(I_\cH-A(z)) \underset{z\to 0}{=} {\det}_{\cH}(I_{\cH}-P_0-A_0)
\, {\det}_{P_0\cH}(P_0A_1P_0) z + \Oh\big(z^2\big).  \lb{2.28}
\end{equation}
\end{theorem}
\begin{proof} Since
\begin{equation}
\cH = P(z)\cH \, \dot+ \, Q(z)\cH, \quad P(z)A(z)=A(z)P(z), \quad
P(z)Q(z)=Q(z)P(z)=0,
\end{equation}
one computes using \eqref{2.21},
\begin{align}
& {\det}_{\cH}(I_{\cH}-A(z)) =
{\det}_{\cH}\big(I_{\cH}-T(z)A(z)T(z)^{-1}\big)  \no \\
& \quad
={\det}_{\cH}\big(I_{\cH}-T(z)[P(z)A(z)P(z)+Q(z)A(z)Q(z)]T(z)^{-1}\big)
\no \\
& \quad
={\det}_{\cH}\big(I_{\cH}-P_0T(z)A(z)T(z)^{-1}P_0-Q_0T(z)A(z)T(z)^{-1}Q_0\big)
\no \\
& \quad ={\det}_{P_0\cH}\big(I_{P_0\cH}-P_0T(z)A(z)T(z)^{-1}P_0\big) \no \\
& \qquad \times
{\det}_{Q_0\cH}\big(I_{Q_0\cH}-Q_0T(z)A(z)T(z)^{-1}Q_0\big). \lb{2.29}
\end{align}
Using \eqref{2.19}, \eqref{2.5}, \eqref{2.22}, and \eqref{2.23}, one computes
\begin{align}
P_0T(z)A(z)T(z)^{-1}P_0&  \underset{z\to 0}{=} P_0A_0P_0+P_0A_1P_0 z
+ P_0\Oh\big(z^2\big)P_0,  \lb{2.30} \\
Q_0T(z)A(z)T(z)^{-1}Q_0& \underset{z\to 0}{=} Q_0A_0Q_0+ Q_0\Oh(z)Q_0
\lb{2.31}
\end{align}
for $|z|$ sufficiently small. Relation \eqref{2.31} implies
\begin{equation}
{\det}_{Q_0\cH}(I_{Q_0\cH}-Q_0T(z)A(z)T(z)^{-1}Q_0)  \underset{z\to 0}{=}
{\det}_{Q_0\cH}(I_{Q_0\cH}-Q_0 A_0 Q_0) + \Oh(z) \neq 0   \lb{2.31a}
\end{equation}
for $|z|$ sufficiently small, and hence we next focus on the first
factor on the right-hand side of \eqref{2.29}. Applying
\eqref{2.14a}, \eqref{2.14b}, and \eqref{2.30} one obtains
\begin{align}
& {\det}_{P_0\cH}\big(I_{P_0\cH}-P_0T(z)A(z)T(z)^{-1}P_0\big) \no \\
& \quad \underset{z\to 0}{=}
{\det}_{P_0\cH}(I_{P_0\cH}-P_0A_0P_0-P_0[A_1+\Oh(z)]P_0 z)  \no \\
& \quad  \underset{z\to 0}{=}
(-1)^{n_0} {\det}_{P_0\cH}(P_0D_0P_0+P_0[A_1+\Oh(z)]P_0 z).  \lb{2.31b}
\end{align}
Next, let $S\in\cB\big(P_0\cH,\bbC^{n_0}\big)$ be the transformation
which puts $P_0D_0P_0$ into its Jordan canonical form $\widehat
D_0=SP_0D_0P_0S^{-1}$ and
denote
\begin{equation}
\widehat A_1(z) \underset{z\to 0}{=} SP_0[A_1+\Oh(z)]P_0 S^{-1}
\underset{z\to 0}{=} \widehat A_1 + \Oh(z).
\end{equation}
Then,
\begin{align}
& {\det}_{P_0\cH}\big(I_{P_0\cH}-P_0T(z)A(z)T(z)^{-1}P_0\big) \no \\
& \quad  \underset{z\to 0}{=}
(-1)^{n_0} {\det}_{P_0\cH}(P_0D_0P_0+P_0[A_1+\Oh(z)]P_0 z)  \no \\
& \quad  \underset{z\to 0}{=}
(-1)^{n_0} {\det}_{\bbC^{n_0}}\big(\widehat D_0 + \widehat A_1(z) z\big) \no \\
& \quad  \underset{z\to 0}{=}
(-1)^{n_0} {\det}_{\bbC^{n_0}}\big(\widehat D_0 + \widehat A_1 z +
\Oh\big(z^2\big)\big)
\no \\
& \quad \underset{z\to 0}{=}
  (-z)^{n_0-\nu_0} {\det}_{\bbC^{n_0-\nu_0}}\big(\wti A_1 + \Oh(z)\big) \no \\
& \quad \underset{z\to 0}{=} (-z)^{n_0-\nu_0}
{\det}_{\bbC^{n_0-\nu_0}}\big(\wti A_1\big)
+ \Oh\big(z^{n_0-\nu_0+1}\big),  \lb{2.31c}
\end{align}
by applying the Laplace determinant expansion formula (cf., e.g.,
\cite[Sect.\ 3.3]{VD99})
to ${\det}_{\bbC^{n_0}}\big(\widehat D_0 + \widehat A_1 z
+\Oh\big(z^2\big)\big)$ with respect to the $\nu_0$ columns in
$\widehat D_0$ which contain a $1$. Combining \eqref{2.31c} and
\eqref{2.31a} then proves \eqref{2.25}.

If $1$ is a semisimple eigenvalue of $A_0$ and hence $D_0=0$,
$\nu_0=0$, the first line on the right-hand side in \eqref{2.26} is
clear from \eqref{2.25}. To prove the second line in the right-hand
side of \eqref{2.26}, we recall that
\begin{equation}
P_0A_0P_0=P_0A_0=A_0P_0=P_0 \, \text{ and } \,
Q_0A_0Q_0=Q_0A_0=(I_{\cH}-P_0)A_0=A_0-P_0,   \lb{2.31d}
\end{equation}
and hence,
\begin{align}
{\det}_{\cH}(I_{\cH}-P_0-A_0) & = {\det}_{\cH}(I_{\cH}-P_0-P_0A_0P_0-Q_0A_0Q_0)
\no \\
& = {\det}_{\cH}(-P_0+Q_0-Q_0A_0Q_0)   \no \\
& = (-1)^{n_0} {\det}_{Q_0\cH}(I_{Q_0\cH}-Q_0A_0Q_0) \neq 0.  \lb{2.31e}
\end{align}
The special case $n_0=1$ in \eqref{2.26} then yields \eqref{2.28}.
\end{proof}

\begin{remark} \lb{r2.4}
Since $\nu_0$ can take on any particular value from $0$ to $n_0 -1$
and $A_1$ is generally independent of $A_0$ and hence $D_0$, the
power $n_0-\nu_0$ of $z$ in \eqref{2.25} can take on any value
between $1$ and $n_0$.
\end{remark}

\subsection{Hilbert--Schmidt operators}
Next, we treat the case of $2$-modified Fredholm determinants, where
all relevant
operators are only assumed to lie in the Hilbert--Schmidt class. In addition to
\eqref{2.1}--\eqref{2.3} we recall the following standard facts for
$2$-modified Fredholm determinants ${\det}_{2,\cH}(I_{\cH}-A)$,
$A\in\cB_2(\cH)$ (cf., e.g., \cite{GGK97}, \cite[Ch.\ XIII]{GGK00},
\cite[Sect.\
IV.2]{GK69}, \cite{Si77}, \cite[Ch.\ 3]{Si05}),
\begin{align}
& {\det}_{2,\cH} (I_{\cH}-A)= {\det}_{\cH}((I_{\cH}-A)\exp(A)), \quad
A\in\cB_2 (\cH), \lb{2.32} \\
& {\det}_{2,\cH}((I_{\cH}-A)(I_{\cH}-B))={\det}_{2,\cH}(I_{\cH}-A) \,
{\det}_{2,\cH}(I_{\cH}-B) \, e^{-\tr_{\cH}(AB)},   \lb{2.33} \\
& \hspace*{8.75cm} A, B \in\cB_2(\cH),    \no  \\
& {\det}_{2,\cH}(I_{\cH}-A)={\det}_{\cH}(I_{\cH}-A)
\, e^{\tr_{\cH}(A)}, \quad A\in\cB_1(\cH),   \lb{2.34} \\
& {\det}_{2,\cH'}(I_{\cH'}-AB)={\det}_{2,\cH}(I_\cH-BA) \;\, \text{ for all
$A\in\cB(\cH,\cH')$, $B\in\cB(\cH',\cH)$}  \no \\
& \hspace*{5.2cm} \text{ such that $BA\in \cB_2(\cH)$,
$AB\in \cB_2(\cH')$.}  \lb{2.34aa}
\end{align}
Moreover, in analogy to \eqref{2.8} one now has
\begin{equation}
{\det}_{2,\cH}(I_{\cH}-B z)=\exp\bigg[-\sum_{k=2}^\infty
\f{\tr_{\cH}\big(B^k\big)}{k}z^k\bigg] \, \text{ for $|z|$
sufficiently small, $B\in\cB_2(\cH)$.}  \lb{2.34a}
\end{equation}
Finally, assuming $A, B \in\cB_2(\cH)$, we mention some estimates to be useful in Section \ref{s4}: For some $C>0$,
\begin{align}
& |{\det}_{2,\cH}(I_{\cH}-A)| \leq \exp\big(C\|A\|_{\cB_2(\cH)}^2\big),   \lb{2.34b} \\
& |{\det}_{2,\cH}(I_{\cH}-A) - {\det}_{2,\cH}(I_{\cH}-B)| \leq \|A-B\|_{\cB_2(\cH)}
\lb{2.34c} \\
& \hspace*{5.75cm}
\times \exp\big(C[\|A\|_{\cB_2(\cH)}+\|B\|_{\cB_2(\cH)}+1]^2\big)    \no
\end{align}

\begin{hypothesis}  \lb{h2.5}
Suppose $A(\dott)\in \cB_2(\cH)$ is a family of Hilbert--Schmidt operators on
$\cH$ analytic on an open neighborhood $\Omega_0\subset\bbC$ of $z=0$ in the
Hilbert--Schmidt norm $\|\dott\|_{\cB_2(\cH)}$.
\end{hypothesis}

Given Hypothesis \ref{h2.5}, we write again
\begin{equation}
A(z) \underset{z\to 0}{=} A_0+A_1 z+\Oh\big(z^2\big) \, \text{ for
$z\in\Omega_0$ sufficiently small, }  \, A_\ell \in \cB_2(\cH), \;
\ell=0,1.  \lb{2.35}
\end{equation}

We start with the analog of Lemma \ref{2.2}.

\begin{lemma}  \lb{l2.6}
Assume Hypothesis \ref{h2.5} and suppose $(I_{\cH}-A_0)^{-1}\in\cB(\cH)$. Then,
\begin{align}
\begin{split}
{\det}_{2,\cH} (I_\cH-A(z)) & \underset{z\to 0}{=} {\det}_{2,\cH}
(I_{\cH}-A_0)    \\
& \quad -{\det}_{2,\cH} (I_{\cH}-A_0) \,
\tr_{\cH}\big((I_{\cH}-A_0)^{-1}A_0 A_1\big) z
+ \Oh\big(z^2\big).   \lb{2.36}
\end{split}
\end{align}
\end{lemma}
\begin{proof}
This follows most easily from rewriting \eqref{2.6} in terms of modified
Fredholm determinants ${\det}_{2,\cH}(I_{\cH}+\dott)$, using
\eqref{2.34}, and then approximating Hilbert--Schmidt operators by trace
class (or finite-rank) operators (cf., e.g., \cite[Theorem
III.7.1]{GK69}). Indeed, one computes using \eqref{2.34} repeatedly,
\begin{align}
& {\det}_{2,\cH} (I_{\cH} - A(z)) = {\det}_{\cH}(I_{\cH} -
A(z))e^{\tr_{\cH}(A(z))}  \no \\
& \quad = {\det}_{\cH}\big((I_{\cH}-A_0)\big[I_{\cH} - (I_{\cH}-A_0)^{-1}
(A(z)-A_0)\big]\big)e^{\tr_{\cH}(A(z))} \no \\
& \quad =  {\det}_{\cH}(I_{\cH}-A_0) {\det}_{\cH}\big(I_{\cH} -
(I_{\cH}-A_0)^{-1}
(A(z)-A_0)\big)e^{\tr_{\cH}(A(z))} \no \\
& \quad =  {\det}_{2,\cH}(I_{\cH}-A_0) e^{-\tr_{\cH}(A_0)}
{\det}_{2,\cH}\big(I_{\cH} - (I_{\cH}-A_0)^{-1}(A(z)-A_0)\big) \no \\
& \qquad \times e^{-\tr_{\cH}((I_{\cH}-A_0)^{-1}(A(z)-A_0))}
e^{\tr_{\cH}(A(z))}
\no \\
& \quad = {\det}_{2,\cH}(I_{\cH}-A_0)
{\det}_{2,\cH}\big(I_{\cH} - (I_{\cH}-A_0)^{-1}(A(z)-A_0)\big) \no \\
& \qquad \times e^{\tr_{\cH}([I_{\cH}-(I_{\cH}-A_0)^{-1}](A(z)-A_0))} \no \\
& \quad = {\det}_{2,\cH}(I_{\cH}-A_0)
{\det}_{2,\cH}\big(I_{\cH} - (I_{\cH}-A_0)^{-1}(A(z)-A_0)\big) \no \\
& \qquad \times e^{- \tr_{\cH}((I_{\cH}-A_0)^{-1} A_0(A(z)-A_0))} \no \\
& \quad \underset{z\to 0}{=} {\det}_{2,\cH}(I_{\cH}-A_0)
{\det}_{2,\cH}\big(I_{\cH} - (I_{\cH}-A_0)^{-1}A_1 z +
\Oh\big(z^2\big)\big)  \no \\
& \qquad \times e^{- \tr_{\cH}((I_{\cH}-A_0)^{-1} A_0A_1 z + \Oh(z^2))} \no \\
& \quad \underset{z\to 0}{=} {\det}_{2,\cH}(I_{\cH}-A_0)
{\det}_{2,\cH}\big(I_{\cH} - (I_{\cH}-A_0)^{-1}A_0 A_1 z +
\Oh\big(z^2\big)\big)  \no \\
& \qquad \times \big[1 - \tr_{\cH}\big((I_{\cH}-A_0)^{-1}A_0 A_1\big)
z + \Oh\big(z^2\big)\big]
\no \\
& \quad \underset{z\to 0}{=} {\det}_{2,\cH}(I_{\cH}-A_0)  \no \\
& \qquad - {\det}_{2,\cH}(I_{\cH}-A_0)
\tr_{\cH}\big((I_{\cH}-A_0)^{-1}A_0 A_1\big) z + \Oh\big(z^2\big).
\end{align}
Here we used (cf.\ \eqref{2.34a})
\begin{equation}
{\det}_{2,\cH}(I_{\cH}-B(z) z) \underset{z\to 0}{=} 1 + \Oh\big(z^2\big)
\end{equation}
for $B(\dott)$ analytic in $\cB_2(\cH)$-norm near $z=0$.
\end{proof}

In exactly the same manner one obtains the Hilbert--Schmidt operator version of
Theorem \ref{t2.3}. Again we rely on the notation introduced in the
paragraph preceding Theorem \ref{t2.3}.

\begin{theorem}  \lb{t2.7}
Assume Hypothesis \ref{h2.5} and let
$1\in\sigma_{\rm d}(A_0)$. Then,
\begin{align}
&{\det}_{2,\cH}(I_\cH-A(z))  \underset{z\to 0}{=}
\big[{\det}_{2,Q_0\cH}(I_{Q_0\cH}-Q_0A_0Q_0)+\Oh(z)\big] e^{n_0}   \no \\
& \hspace*{3cm}
\times (-1)^{n_0} {\det}_{P_0\cH}(P_0D_0P_0+P_0[A_1+\Oh(z)]P_0 z)
\lb{2.40}  \\
&\quad \underset{z\to 0}{=} {\det}_{2,Q_0\cH}(I_{Q_0\cH}-Q_0A_0Q_0) e^{n_0}
{\det}_{\bbC^{n_0-\nu_0}}\big(\wti A_1\big)
(-z)^{n_0-\nu_0} + \Oh\big(z^{n_0-\nu_0 +1}\big).    \no
\end{align}
Here $n_0 -\nu_0$ represents the geometric multiplicity 
of the eigenvalue $1$ of $A_0$. 
In the special case where $1$ is a semisimple eigenvalue of $A_0$
$($i.e., where
$D_0=0$ and $\nu_0=0$$)$ one obtains,
\begin{align}
{\det}_{2,\cH}(I_\cH-A(z)) & \underset{z\to 0}{=}
{\det}_{2,Q_0\cH}(I_{Q_0\cH}-Q_0A_0Q_0) e^{n_0}
{\det}_{P_0\cH}(P_0A_1P_0) (-z)^{n_0} \no \\
& \quad + \Oh\big(z^{n_0+1}\big)   \lb{2.41} \\
& \underset{z\to 0}{=}  {\det}_{2,\cH}(I_{\cH}-P_0-A_0) e^{n_0}
{\det}_{P_0\cH}(P_0A_1P_0) z^{n_0} + \Oh\big(z^{n_0+1}\big).   \no
\end{align}
In particular, if $1$ is a simple eigenvalue of $A_0$ $($i.e., if
$n_0=1$, $D_0=0$, and
$\nu_0=0$$)$, one obtains
\begin{equation}
{\det}_{2,\cH}(I_\cH-A(z)) \underset{z\to 0}{=}
{\det}_{2,\cH}(I_{\cH}-P_0-A_0)
\, {\det}_{P_0\cH}(P_0A_1P_0) z + \Oh\big(z^2\big).  \lb{2.42}
\end{equation}
\end{theorem}
\begin{proof}
Again this follows from rewriting \eqref{2.25}--\eqref{2.28} in terms
of modified
Fredholm determinants ${\det}_{2,\cH}(I_{\cH}+\dott)$, using
\eqref{2.34},
\begin{align}
\tr_{\cH}(A(z))& =
\tr_{\cH}\big((P_0+Q_0)T(z)A(z)T(z)^{-1}(P_0+Q_0)\big)  \no \\
& =  \tr_{\cH}\big(P_0T(z)A(z)T(z)^{-1}P_0\big) +
  \tr_{\cH}\big(Q_0T(z)A(z)T(z)^{-1}Q_0\big)  \no \\
& = \tr_{P_0\cH}\big(P_0T(z)A(z)T(z)^{-1}P_0\big) +
  \tr_{Q_0\cH}\big(Q_0T(z)A(z)T(z)^{-1}Q_0\big)  \no \\
& \hspace*{-1.5mm} \underset{z\to 0}{=} \tr_{P_0\cH}(P_0A_0P_0) + \Oh(z) +
  \tr_{Q_0\cH}\big(Q_0T(z)A(z)T(z)^{-1}Q_0\big)   \no \\
  & \hspace*{-1.5mm}  \underset{z\to 0}{=}  n_0 + \Oh(z) +
  \tr_{Q_0\cH}\big(Q_0T(z)A(z)T(z)^{-1}Q_0\big),   \lb{2.43}
\end{align}
and then approximating Hilbert--Schmidt operators by trace
class (or finite-rank) operators (cf., e.g., \cite[Theorem
III.7.1]{GK69}). Here the
fact that $\tr_{P_0\cH}(P_0A_0P_0)=n_0$ follows from the canonical
Jordan structure of
$SP_0A_0P_0S^{-1}$ .  Explicitly, one computes
\begin{align}
& {\det}_{2,\cH} (I_{\cH}-A(z)) = {\det}_{\cH} (I_{\cH}-A(z))
e^{\tr_{\cH}(A(z))}  \no \\
& \quad = {\det}_{P_0\cH}\big(I_{P_0 \cH}-P_0T(z)A(z)T(z)^{-1}P_0\big)  \no \\
& \qquad \times {\det}_{2,Q_0 \cH} \big(I_{Q_0
\cH}-Q_0T(z)A(z)T(z)^{-1}Q_0\big) \no \\
& \qquad \times e^{-\tr_{Q_0 \cH}(Q_0T(z)A(z)T(z)^{-1}Q_0)}
e^{\tr_{\cH}(A(z))}  \no \\
& \quad \underset{z\to 0}{=} {\det}_{P_0\cH}\big(I_{P_0 \cH}
- P_0T(z)A(z)T(z)^{-1}P_0\big)
[{\det}_{2,Q_0 \cH} (I_{Q_0 \cH}-Q_0A_0Q_0) + \Oh(z)]   \no \\
& \qquad \times e^{n_0 + \Oh(z)}  \no \\
& \quad \underset{z\to 0}{=} (-1)^{n_0} {\det}_{P_0\cH}(P_0D_0P_0
+ P_0[A_1 + \Oh(z)]P_0z))   \no \\
& \qquad \times [{\det}_{2,Q_0 \cH} (I_{Q_0 \cH}-Q_0A_0Q_0) + \Oh(z)] e^{n_0}.
\end{align}
\end{proof}

\subsection{Higher modified determinants}
It is now obvious how to proceed in connection with higher modified
Fredholm determinants ${\det}_p(\dott)$, $p\in\bbN$, $n\geq 3$, when
the family $A(\dott)$ in Hypothesis \ref{h2.5} is replaced by one
analytic near $z=0$ in $\cB_p(\cH)$-norm. We omit further details at
this point.

\section{One-Dimensional Reaction--Diffusion Equations
and Schr\"odinger Operators with Nontrivial Spatial Asymptotics} \label{s3}

As an elementary illustration of formula \eqref{2.28} and a view
toward applications to reaction-diffusion equations, we now
illustrate the abstract results of Section \ref{s2} in the context of
one-dimensional Schr\"odinger operators. (For background literature
on reaction-diffusion equations we refer, e.g., to \cite{Pa92},
\cite{Ro84}, \cite{Sm94}.)

To motivate our considerations of one-dimensional Schr\"odinger
operators in this section, we start with a brief discussion of a
simple model for one-dimensional scalar reaction-diffusion equations
of the type
\begin{equation}
w_t = w_{xx} + f(w), \; t>0,    \quad  w(0) = w_0.   \lb{3.0}
\end{equation}
Here, we assume, for simplicity,
\begin{equation}
f\in C^1(\bbR), \; f, f' \in L^\infty(\bbR),   \quad
w_0\in H^2(\bbR)\cap C^\infty(\bbR), \; w_0' \in H^2(\bbR),
\end{equation}
with $w=w(x,t)$, $(x,t)\in \bbR\times [0,T]$,
for some $T=T(w_0)>0$, a
mild solution of \eqref{3.0}  satisfying $w \in C_b([0,T];
L^2(\bbR;dx))$.

As usual, $C_b([0,T]; L^2(\bbR;dx))$ denotes the space
of bounded
continuous maps from $[0,T]$ with values in $L^2(\bbR;dx)$ and the
sup norm
\begin{equation}
\|v\|_{C_b([0,T]; L^2(\bbR;dx))}=\sup_{s\in[0,T]}\|v(s)\|_{L^2(\bbR;dx)},
\quad v \in C_b([0,T]; L^2(\bbR;dx)).
\end{equation}
We refer, for instance, to \cite{LLMP05} for more details in this
context (see also,
\cite[Sect.\ 3.2]{He81}). Also, $H^{m}(\bbR)$,
$m \in\bbN$, denotes the standard Sobolev spaces of regular
distributions which together with their derivatives up to the
$m$th-order  lie in $L^2(\bbR;dx)$.

Assuming that $U$ is a stationary (steady state) solution of
\eqref{3.0}, that is,
$U_t=0$, and hence
\begin{equation}
U'' + f(U) =0,  \lb{3.0a}
\end{equation}
we now linearize \eqref{3.0} around $U$ and obtain
the linearized problem,
\begin{equation}
v_t = L v,
\end{equation}
where $L$ in $L^2(\bbR;dx)$ is given by
\begin{equation}
L = \f{d^2}{d x^2} + f'(U),  \quad \dom(L)=H^{2}(\bbR),
\end{equation}
Since $U' \in H^{2}(\bbR)\cap C^\infty(\bbR)$, one infers upon
differentiating \eqref{3.0a} with respect to $x$ that
\begin{equation}\label{3.7a}
L(U') = (U')'' + f'(U) U' = [U'' + f(U)]' = 0,
\end{equation}
and hence, $0 \in \sigma_{\rm p} (L)$. Here $\sigma_{\rm p}(\dott)$
abbreviates the point spectrum (i.e., the set of  eigenvalues). Thus,
the Schr\"odinger operator
\begin{equation}
H = - L= -\f{d^2}{d x^2} - f'(U),  \quad \dom(H)=H^{2}(\bbR),    \lb{3.0b}
\end{equation}
in $L^2(\bbR;dx)$ with the potential
\begin{equation}
V(x) = - f'(U(x)), \quad x\in\bbR,    \lb{3.0c}
\end{equation}
has the special property of a zero eigenvalue, that is,
\begin{equation}
0 \in \sigma_{\rm p} (H),   \lb{3.0d}
\end{equation}
and we will be studying the situation where in addition $0$ is a
discrete (and hence simple) eigenvalue of $H$ in great detail in the
remainder of this section. In particular, under appropriate
assumptions on the the ``potential'' term $-f'(U)$ in $H$, the
Birman--Schwinger-type integral operator $K(z)$ associated with $H$
and the complex eigenvalue parameter $z\in\bbC$
will be a trace class operator with $K(0)$ having the eigenvalue $1$.
Consequently, the behavior of
$\det_{L^2(\bbR;dx)} (I_{L^2(\bbR;dx)} - K(z))$ for $z$ in a
sufficiently small open neighborhood of $z=0$, which determines
linear stability of \eqref{3.0} around its stationary solution $U$
(cf.\ \cite[Sect.\ 5.1, Ex. 6, Sect.\ 6.2]{He81},
\cite[Sects.\ 9.1.4, 9.1.5]{Lu95}), can be determined in accordance with
Theorem \ref{t2.3} or alternatively, directly from Jost
(respectively, Evans) function considerations, as discussed in Remark
\ref{r3.4a}.

We conclude these illustrations with the following elementary example:

\begin{example} \lb{e3.0}
Let $n\in\bbN$, $n\geq 2$, $c>0$, $\kappa>0$, $x\in\bbR$, and consider
\begin{align}
f_n(w) &= - (n-1)^2  \kappa^2 w + (n-1)n \kappa^2 c^{-2/(n-1)}
w^{(n+1)/(n-1)},  \\
U_n(x) &= c [\cosh(\kappa x)]^{-n+1}.
\end{align}
Then the potential $V_n(x) = -f'(U_n(x))$ in the corresponding
Schr\"odinger operator $H_n=- (d^2/d x^2) + V_n$ coincides with a
particular family of $n$-soliton
Korteweg--de Vries (KdV) potentials
\begin{align}
& V_n(x) = - f'(U_n(x)) = (n-1)^2 \kappa^2 - n(n+1) \kappa^2
[\cosh(\kappa x)]^{-2},  \\
& \lim_{x\to\pm\infty} V_n(x) = (n-1)^2 \kappa^2 >0,
\end{align}
and the zero-energy eigenfunction $U'_n$ satisfying $H_n (U'_n)=0$ is given by
\begin{equation}
U'_n(x) =  - (n-1) \kappa c \sinh(\kappa x) [\cosh(\kappa x)]^{-n}.
\end{equation}
\end{example}

It seems a curious coincidence that $V_n$ should coincide with a particular family of
$n$-soliton KdV potentials (there are many other such $n$-soliton KdV
potentials, cf., e.g., \cite[Example 1.31]{GH03}) in this reaction-diffusion
equation context.

\medskip

Motivated by these considerations, we now start to investigate one-dimensional
Schr\"odinger operators with a scalar potential displaying 
nonzero
asymptotics as
$|x|\to\infty$.

\begin{hypothesis}  \lb{h3.1}
Assume $V\colon\bbR\to\bbR$ is Lebesgue measurable and satisfies
\begin{equation}
\int_\bbR dx\, |V(x)-V_{\infty}|<\infty \, \text{ for some } \,
V_{\infty}>0.   \lb{3.1}
\end{equation}
\end{hypothesis}

Assuming Hypothesis \ref{h3.1} we introduce
\begin{align}
& W(x) =V(x)-V_{\infty}, \label{defW}\\
& W(x)=v(x) u(x), \quad u(x)=\sgn(W(x)) v(x), \;\; v(x)=|W(x)|^{1/2},
\end{align} 
for a.e.\ $x\in\bbR$, and define $H$ to be the (maximally defined)
self-adjoint realization in $L^2(\bbR;dx)$ of the differential
expression $\tau = - (d^2/d x^2) +V(x)$, $x\in\bbR$,
obtained by the method of quadratic forms, or equivalently, by using
the limit point theory for self-adjoint 2nd order  ordinary
differential operators,
\begin{align}
& Hf=\tau f, \\
& f\in\dom(H)=\big\{g\in L^2(\bbR;dx)\,\big|\, g,g'\in
AC_{\loc}(\bbR); \tau g\in L^2(\bbR;dx)\big\}.
\end{align}
(Here $\prime$ denotes differentiation with respect to $x\in\bbR$.)
This also implies
\begin{equation}
\dom\big(H^{1/2}\big)=H^1(\bbR).
\end{equation}
We also introduce the self-adjoint operator $H^{(0)}$ in $L^2(\bbR;dx)$ associated with the differential
expression $\tau^{(0)}= - (d^2/d x^2) +V_\infty$, $x\in\bbR$,
replacing $V(\dott)$ by its asymptotic value $V_\infty$ (in the sense of \eqref{3.1}),
\begin{equation}\lb{H0}
H^{(0)}f=\tau^{(0)} f, \quad f \in \dom\big(H^{(0)}\big)=H^2(\bbR).
\end{equation}
By well-known results, \eqref{3.1} implies
\begin{equation}
\sigma_{\rm ess}(H)=\sigma_{\rm ess}\big(H^{(0)}\big)
=[V_{\infty},\infty), \quad \sigma_{\rm d}(H) \subset (-\infty,V_{\infty}),
\end{equation}
where $\sigma_{\rm ess}(\dott)$ and $\sigma_{\rm d}(\dott)$ denote
the essential
and discrete spectrum, respectively.

Applying the Birman--Schwinger principle (cf., e.g., \cite[Sect.\ 3]{GLMZ05}),
\begin{equation}
H\Psi(\lambda_j)=\lambda_j \Psi(\lambda_j), \quad
\lambda_j<V_{\infty}, \; \lambda_j\in\sigma_{\rm d}(H), \;
\Psi(\lambda_j)\in\dom(H),
\lb{3.8}
\end{equation}
is equivalent to
\begin{equation}
K(\lambda_j) \Phi(\lambda_j)=\Phi(\lambda_j), \quad \lambda_j <V_{\infty}, \;
\Phi(\lambda_j)\in L^2(\bbR;dx),   \lb{3.9}
\end{equation}
with equal finite geometric multiplicity of either eigenvalue problem
\eqref{3.8} and \eqref{3.9}. In particular, in this special
one-dimensional context, the eigenvalue $\lambda_j$ of $H$ as well as
the eigenvalue $1$ of $K(\lambda_j)$ are necessarily simple. Here we
abbreviated
\begin{equation}
K(z)= -\ol{u\big(H^{(0)}-z I_{L^2(\bbR;dx)}\big)^{-1}v},
\quad z\in\bbC\backslash [V_{\infty},\infty),   \lb{3.10}
\end{equation}
with $\ol T$ denoting the operator closure of $T$. We recall that the integral
kernel $\big(H^{(0)}-z
I_{L^2(\bbR;dx)}\big)^{-1}(x,x')$ of the resolvent of $H^{(0)}$ is
explicitly given by
\begin{align}
\begin{split}
& \big(H^{(0)}-z I_{L^2(\bbR;dx)}\big)^{-1}(x,x')= (i/ 2) (z-V_\infty)^{-1/2}
e^{i (z-V_\infty)^{1/2}|x-x'|},   \lb{3.12}   \\
& \hspace*{2.35cm}
z\in \bbC\backslash [V_\infty,\infty), \; \Im\big((z-V_\infty)^{1/2}\big) >0,
\quad x, x'\in\bbR,
\end{split}
\end{align}
and hence
\begin{equation}
K(z) \in \cB_1\big(L^2(\bbR;dx)\big), \quad z \in\bbC\backslash [V_\infty,\infty), 
\end{equation}
since $K(z)$ is the product of the two Hilbert--Schmidt operators 
 (cf.\ \cite[Ch.\ 4]{Si05})
\begin{equation}
u \big(H^{(0)}-z I_{L^2(\bbR;dx)}\big)^{-1/2} \, \text{ and } \,  
\ol{\big(H^{(0)}-z I_{L^2(\bbR;dx)}\big)^{-1/2} v}
\end{equation}
(no operator closure necessary in the first factor).   
In addition, $\Psi(\lambda_j)$ and $\Phi(\lambda_j)$ are related by
\begin{equation}
\Phi(\lambda_j,x)=u(x)\Psi(\lambda_j,x), \quad x\in\bbR,
\end{equation}
and we note that $\Psi(\lambda_j,\dott)$ is also bounded,
\begin{equation}
\Psi(\lambda_j,\dott) \in L^\infty(\bbR)
\end{equation}
(in fact,  even exponentially decaying with respect to $x$ by standard iterations of the Volterra integral equations \eqref{2.51}). Here the discrete eigenvalues
$\sigma_{\rm d}(H)=\{\lambda_j\}_{j\in J}$ of $H$, with $J\subseteq
\bbN$ an appropriate (finite or infinite) index set, are ordered as
follows:
\begin{equation}
\lambda_1 < \lambda_2 < \cdots < V_{\infty}.    \lb{3.14}
\end{equation}

Moreover, one obtains
\begin{equation}
K(\ol z)^*= -\ol{v\big(H^{(0)}-z I_{L^2(\bbR;dx)}\big)^{-1}u},
\quad z\in\bbC\backslash [V_{\infty},\infty),
\end{equation}
and hence
\begin{equation}
K(\ol z)^*= SK(z)S^{-1}, \quad z\in\bbC\backslash [V_{\infty},\infty),
\end{equation}
where $S$ denotes the unitary operator of multiplication by $\sgn(W(\dott))$
in $L^2(\bbR;dx)$,
\begin{equation}
(S f)(x) = \sgn(W(x)) f(x) \, \text{ for a.e.\ $x\in\bbR$, } \, f \in
L^2(\bbR;dx).
\end{equation}
Thus,
\begin{equation}
K(\lambda_j)^* \wti\Phi(\lambda_j)=\wti\Phi(\lambda_j), \quad
\lambda_j <V_{\infty}, \;
\wti\Phi(\lambda_j)\in L^2(\bbR;dx),
\end{equation}
where
\begin{equation}
\wti\Phi(\lambda_j,x)=C(\lambda_j) (S \Phi(\lambda_j))(x)
=C(\lambda_j) v(x) \Psi(\lambda_j,x) \, \text{ for some } \,
C(\lambda_j) \in\bbC\backslash\{0\}.
\end{equation}
To fix the normalization constant $C(\lambda_j)$ we require
\begin{equation}
\big(\wti\Phi(\lambda_j),\Phi(\lambda_j)\big)_{L^2(\bbR;dx)} = 1.
\end{equation}
This yields
\begin{equation}
C(\lambda_j)=\bigg[\int_{\bbR} dx\, W(x) |\Psi(\lambda_j,x)|^2 \bigg]^{-1}
=(v\Psi(\lambda_j), u\Psi(\lambda_j))^{-1}_{L^2(\bbR; dx)}    \lb{3.22}
\end{equation}
and hence
\begin{align}
\wti\Phi(\lambda_j,x)&=\bigg[\int_{\bbR} dx\, W(x)
|\Psi(\lambda_j,x)|^2 \bigg]^{-1}
\sgn(W(x)) \Phi(\lambda_j,x) \\
&=\bigg[\int_{\bbR} dx\, W(x) |\Psi(\lambda_j,x)|^2 \bigg]^{-1}
v(x) \Psi(\lambda_j,x).
\end{align}
In particular, the corresponding one-dimensional Riesz projection
$P(\lambda_j)$ onto the eigenspace corresponding to the eigenvalue
$1$ of $K(\lambda_j)$ is then of the simple form
\begin{equation}
P(\lambda_j) = \big(\wti\Phi(\lambda_j),\dott \big)_{L^2(\bbR;dx)}
\Phi(\lambda_j), \quad
\lambda_j  < V_{\infty}.
\end{equation}

\subsection{The Jost--Pais derivative formula}
Next, consider the generalized Jost-type (distributional) solutions
\begin{align}
\Psi_\pm(z,x) &= e^{\pm i (z-V_{\infty})^{1/2}x}  \no \\
& \quad - \int_x^{\pm\infty} dx' \,
\f{\sin((z-V_{\infty})^{1/2}(x-x'))}{(z-V_{\infty})^{1/2}}
[V(x')-V_{\infty}] \Psi_\pm(z,x'),
\lb{2.51}  \\
& \hspace*{2.55cm} z\in\bbC\backslash[V_{\infty},\infty), \;
\Im\big((z-V_{\infty})^{1/2}\big)>0,
\; x\in\bbR,   \no
\end{align}
of $L\Psi(z)=z\Psi(z)$. Then,
\begin{equation}
\Psi_{\pm}(\lambda,x) \, \text{ are real-valued for
$\lambda<V_{\infty}$, $x\in\bbR$,}
\lb{2.52}
\end{equation}
and
\begin{equation}
\Psi_{\pm}(\lambda,x)>0 \, \text{ for $\lambda<V_{\infty}$ and $\pm
x$ sufficiently large.}
\lb{2.52a}
\end{equation}

The Jost function $\cF$ associated with $H$ is then given by
\begin{align}
\cF(z) & = \f{\Wr(\Psi_-(z),\Psi_+(z))}{2i (z-V_{\infty})^{1/2}} \lb{2.53} \\
& = 1 + \f{i}{2(z-V_{\infty})^{1/2}} \int_{\bbR} dx \, e^{\mp i
(z-V_{\infty})^{1/2}x} [V(x)-V_{\infty}]
\Psi_{\pm}(z,x),   \lb{2.54}  \\
& \hspace*{7.1cm} z\in\bbC\backslash [V_{\infty},\infty),  \no
\end{align}
where $\Wr(f,g)(x)=f(x)g'(x)-f'(x)g(x)$, $x\in\bbR$, $f, g \in
C^1(\bbR)$, denotes the
Wronskian of $f$ and $g$. The limits to the real axis
\begin{equation}
\lim_{\varepsilon\downarrow 0}\cF(\lambda\pm i \varepsilon) =
\cF(\lambda \pm i 0)
\end{equation}
exist and are continuous for all
$\lambda\in\bbR\backslash\{V_{\infty}\}$. In addition, one verifies
\begin{equation}
\cF(z)\underset{|z|\to\infty}{\longrightarrow} 1.    \lb{2.54aa}
\end{equation}

Moreover, one can prove the following result originally due to Jost and Pais
\cite{JP51} (for $V_\infty=0$) in the context of half-line
Schr\"odinger operators. The actual case at hand of Schr\"odinger
operators on the whole real line (again for
$V_\infty=0$) was discussed by Newton \cite{Ne80} and we refer to
\cite{GM04} for more background and details,
\begin{equation}
\cF(z) = {\det}_{L^2(\bbR;dx)}(I_{L^2(\bbR;dx)}-K(z)), \quad
z\in\bbC\backslash[V_{\infty},\infty).
\lb{2.54a}
\end{equation}

Since we are interested especially in the $z$-derivative of $\cF(z)$
at a discrete  eigenvalue of $H$, we now prove the following result,
Lemma \ref{l3.2}. For the remainder of this paper we abbreviate
differentiation with respect to the complex-valued spectral
parameter $z\in\bbC$ by $\bullet$ (to distinguish it from
differentiation with respect to the space variable $x\in\bbR$).

\begin{lemma}  \lb{l3.2}
Assume Hypothesis \ref{h3.1} and $z\in\bbC\backslash
[V_\infty,\infty)$. Moreover, let
$\lambda_j<V_{\infty}$, $\lambda_j \in \sigma_{\rm d}(H)$. Then,
$\cF(\lambda_j)= 0$ and
\begin{equation}
\cF^{\bullet}(\lambda_j)=
\f{-1}{2(V_{\infty}-\lambda_j)^{1/2}}\int_{\bbR} dx \,
\Psi_-(\lambda_j,x) \Psi_+(\lambda_j,x).   \lb{2.54b}
\end{equation}
\end{lemma}
\begin{proof}
Consider
\begin{equation}
\Psi''_{\pm}(z,x)=[V(x)-z]\Psi_{\pm}(z,x), \quad
\Psi^{\bullet\prime\prime}_{\pm}(z,x)=[V(x)-z]\Psi^{\bullet}_{\pm}(z,x)-\Psi_{\pm}(z,x)
\end{equation}
(we recall that $\bullet$ abbreviates $d/dz$) to derive the identities
\begin{align}
\f{d}{dx}\Wr(\Psi_-(z,x),\Psi^{\bullet}_+(z,x)) &=-\Psi_-(z,x)\Psi_+(z,x),  \\
\f{d}{dx}\Wr(\Psi^{\bullet}_-(z,x),\Psi_+(z,x)) &= \Psi_-(z,x)\Psi_+(z,x).
\end{align}
Then one obtains for all $R>0$,
\begin{align}
& \Wr(\Psi_-(z), \Psi^{\bullet}_+(z))(x) + \Wr(\Psi^{\bullet}_-(z),
\Psi_+(z))(x)  \no \\
& \qquad-\Wr(\Psi_-(z), \Psi^{\bullet}_+(z))(R) -
\Wr(\Psi^{\bullet}_-(z), \Psi_+(z))(-R) \no \\
& \quad = -\int_x^R dx' \, \f{d}{dx'}\Wr(\Psi_-(z,x'),\Psi^{\bullet}_+(z,x'))
+ \int_{-R}^x dx' \, \f{d}{dx'}\Wr(\Psi^{\bullet}_-(z,x'),\Psi_+(z,x'))  \no \\
& \quad = \int_{-R}^R dx' \, \Psi_-(z,x') \Psi_+(z,x').
\end{align}
Next, we note that
\begin{equation}
\f{d}{dz} \Wr(\Psi_-(z), \Psi_+(z)) = \Wr(\Psi_-(z), \Psi^{\bullet}_+(z))
+ \Wr(\Psi^{\bullet}_-(z), \Psi_+(z)),
\end{equation}
and choosing $z=\lambda_j<V_{\infty}$, $\lambda_j \in \sigma_{\rm
d}(H)$, one concludes  $\Psi_{\pm}(\lambda_j,\dott)\in L^2(\bbR; dx)$
and hence,
\begin{equation}
\f{d}{dz} \Wr( \Psi_-(z), \Psi_+(z))\big|_{z=\lambda_j} = \int_\bbR dx
\, \Psi_-(\lambda_j,x) \Psi_+(\lambda_j,x).
\end{equation}
Together with \eqref{2.53} and
\begin{equation}
2i (\lambda_j -V_{\infty})^{1/2} \cF(\lambda_j)
=\Wr( \Psi_-(\lambda_j), \Psi_+(\lambda_j)) = 0,
\end{equation}
this establishes \eqref{2.54b}.
\end{proof}

Next, we specialize to the case $\lambda_j=0$ and hence assume
\begin{equation}
0 \in \sigma_{\rm d}(H).
\end{equation}
In this context we then denote
\begin{equation}
\Psi_0=\Psi(0), \; \Phi_0=\Phi(0), \; \wti\Phi_0=\wti\Phi(0), \;
C_0=C(0), \; K_0=K(0), \; P_0=P(0), \, \text{ etc.,}
\end{equation}
and without loss of generality (cf.\ \eqref{2.52}) we assume that
$\Psi_0(x)$ is real-valued for all $x\in\bbR$.

We summarize the results for $\cF^{\bullet}(0)$:

\begin{lemma}  \lb{l3.3}
Assume Hypothesis \ref{h3.1} and suppose $0 \in \sigma_{\rm d}(H)$. Then,
$\cF(0)=0$ and
\begin{align}
\cF^{\bullet}(0)&=\f{-1}{2 V_{\infty}^{1/2}} \int_{\bbR} dx \,
\Psi_-(0,x) \Psi_+(0,x)
\lb{2.80}   \\
&=\begin{cases} <0 & \text{if $\Psi_{\pm}(0,\dott)$ has an even number of zeros
on $\bbR$,} \\
>0 & \text{if $\Psi_{\pm}(0,\dott)$ has an odd number of zeros on $\bbR$,}
\end{cases}  \lb{2.81}  \\
&=\begin{cases} <0 & \text{if $0$ is an odd eigenvalue of $H$,} \\
>0 & \text{if $0$ is an even eigenvalue of $H$.}
\end{cases}   \lb{2.82}
\end{align}
Here the eigenvalues $\{\lambda_j\}_{j\in J}$, $J\subseteq \bbN$ an appropriate
$($finite or infinite$)$ index set, are ordered in magnitude according to
$\lambda_1 < \lambda_2 < \cdots < V_{\infty}$ $($cf.\ \eqref{3.14}$)$.
\end{lemma}

\begin{proof}
Equation \eqref{2.81} is immediate from \eqref{2.52a} and
$\Psi_+(0,x)=c \Psi_-(0,x)$ for some $c\in\bbR\backslash\{0\}$ (cf.\
\eqref{2.79}). Relation \eqref{2.82} is a direct consequence of the
fact that $H$ is bounded from below,
$H\geq \lambda_1 I_{L^2(\bbR; dx)}$, the discrete eigenvalues of $H$ are in a
one-to-one correspondence with the zeros of $\cF$ on $[\lambda_1,
V_{\infty})$ (the zeros necessarily being all simple), and the fact
that
$\cF(z)\underset{|z|\to\infty}{\longrightarrow} 1$,
$z\in\bbC\backslash[V_{\infty},\infty)$ (cf.\ \eqref{2.54aa}).
\end{proof}

\begin{remark}  \lb{r3.4a}
Since $\cF(\lambda \pm i
0)\underset{\lambda\to\pm\infty}{\longrightarrow} 1$ by
\eqref{2.54aa}, $\cF(0)=0$ and $\cF^{\bullet}(0) > 0$ imply that $H$
has at least one negative eigenvalue and hence $L=-H$ (cf.\
\eqref{3.0b}) has at least one positive eigenvalue. This implies linear 
instability of the stationary solution $U$ in the context of the
reaction-diffusion equation \eqref{3.0} identifying $V(x)$
and $-f'(U(x))$, $x\in\bbR$ (cf.\ \eqref{3.0c}).
For, it is easily seen by consideration of the standing-wave equation
\eqref{3.0a}, a scalar nonlinear oscillator,
that the derivative $U'$ of a pulse-type solution has precisely one
zero, whereas $\Psi_\pm$ are nonzero multiples of the zero eigenfunction $U'$.
For discussion of spectral stability and some of its applications, we
refer, for instance, to \cite{PW92}, \cite{PZ04}, \cite{Sa02},
\cite{Zu03}, \cite{ZH98}, and the references cited therein.
An equivalent formula for $\cF^{\bullet}$ yielding the same
conclusions may be derived in straightforward
fashion by Evans function techniques, following the standard approach
introduced in \cite{Ev72}, \cite{Ev72a}, \cite{Ev72b}, \cite{Ev75}.
We recall (cf.\ \cite{GLM07}) that Jost and Evans functions, suitably normalized,
agree.
\end{remark}

\subsection{Fredholm determinant version}
We now turn to the connection with the abstract approach to the
asymptotic behavior of Fredholm determinants presented in Section
\ref{s2}.

Applying \eqref{2.28} to \eqref{2.54a} one then obtains
\begin{equation}
\cF^{\bullet}(0)
= {\det}_{L^2(\bbR;dx)}(I_{L^2(\bbR;dx)}-K_0-P_0) \, {\det}_{P_0
L^2(\bbR;dx)}(P_0K_1P_0),
\lb{2.55}
\end{equation}
where
\begin{equation}
K(z) \underset{z\to 0}{=} K_0 + K_1 z + \Oh\big(z^2\big)
\end{equation}
with
\begin{align}
K_0& = -\ol{u{H^{(0)}}^{-1}v}, \\
K_1& = K^{\bullet}(z)\big|_{z=0}= - \ol{u{H^{(0)}}^{-2}v}.
\end{align}

\subsubsection{Evaluation of the second factor}
We start by determining the second factor ${\det}_{P_0 L^2(\bbR;dx)}(P_0K_1P_0)$
on the right-hand side of \eqref{2.55}:

\begin{theorem}  \lb{t3.4}
Assume Hypothesis \ref{h3.1} and suppose $0 \in \sigma_{\rm d}(H)$. Then,
\begin{align}
{\det}_{P_0 L^2(\bbR;dx)}(P_0K_1P_0)
& = - \|\Psi_0\|^2_{L^2(\bbR; dx)} \big/ (v\Psi_0,
u\Psi_0)_{L^2(\bbR; dx)}  \lb{2.60} \\
& = \big[\|\Psi_0'\|_{L^2(\bbR;dx)}^2 + V_{\infty}
\|\Psi_0\|_{L^2(\bbR;dx)}^2\big]^{-1}
\|\Psi_0\|_{L^2(\bbR;dx)}^2 > 0.   \lb{2.61}
\end{align}
\end{theorem}
\begin{proof}
First, we choose a compactly supported sequence $W_n \in C^\infty_0(\bbR)$,
$n\in\bbN$, such that
\begin{equation}
W_n=V_n-V_\infty=u_n v_n, \quad u_n= \sgn(W_n)v_n, \;\; v_n=|W_n|^{1/2}, \;
n\in\bbN,
\end{equation}
and
\begin{equation}
\lim_{n\to\infty}\|v_n - v\|_{L^2(\bbR; dx)} =0.   \lb{3.70a}
\end{equation}
Given $W_n$, we introduce the self-adjoint operator sum
$H_n=H^{(0)}+W_n$ in $L^2(\bbR;dx)$ defined on the domain $H^2(\bbR)$
associated with the differential expression $L_n=-(d^2/d x^2) + V_n(x)$,
$x\in\bbR$, $n\in\bbN$. Then $H_n$ converges to $H$ in norm resolvent
sense as $n\to\infty$, that is,
\begin{equation}
\lim_{n\to\infty} \big\|  \big(H_n - z I_{L^2(\bbR;dx)}\big)^{-1} -
  \big(H - z I_{L^2(\bbR;dx)}\big)^{-1} \big\| = 0, \quad
z\in\bbC\backslash\bbR.
\end{equation}
This follows from the resolvent identities,
\begin{align}
& \big(H_n - z I_{L^2(\bbR;dx)}\big)^{-1} = \big(H^{(0)} - z
I_{L^2(\bbR;dx)}\big)^{-1}  \no \\
& \quad -
\big(H^{(0)} - z I_{L^2(\bbR;dx)}\big)^{-1}  v_n
\big[I_{L^2(\bbR;dx)} - K_n(z)\big]^{-1} u_n \big(H^{(0)} - z
I_{L^2(\bbR;dx)}\big)^{-1},  \no \\
& \hspace*{8.5cm}  z\in \bbC\backslash\sigma(H_n),  \\
& \big(H - z I_{L^2(\bbR;dx)}\big)^{-1} = \big(H^{(0)} - z
I_{L^2(\bbR;dx)}\big)^{-1}  \no \\
& \quad -
\big(H^{(0)} - z I_{L^2(\bbR;dx)}\big)^{-1}  v \big[I_{L^2(\bbR;dx)}
- K(z)\big]^{-1} u \big(H^{(0)} - z I_{L^2(\bbR;dx)}\big)^{-1},  \no
\\
& \hspace*{8.5cm}  z\in \bbC\backslash\sigma(H),
\end{align}
where
\begin{equation}
K_n (z) = u_n \big(H^{(0)} - z I_{L^2(\bbR;dx)}\big)^{-1} v_n, \quad
z\in \bbC\backslash\sigma\big(H^{(0)}\big), \; n\in\bbN,
\end{equation}
and $K(z)$ is given by \eqref{3.10}, and the fact that \eqref{3.70a} implies
\begin{align}
&\lim_{n\to\infty}\big\|(v_n -v) \big(H^{(0)} - z I_{L^2(\bbR;dx)}\big)^{-1/2}
  \big\|_{\cB_2(L^2(\bbR;dx))} =0, \quad z\in
\bbC\backslash\sigma\big(H^{(0)}\big)
\end{align}
(cf.\ the detailed discussion in \cite{GM04}). Thus, the spectrum of
$H_n$ converges to that of $H$ as $n\to\infty$. In particular, for
$n\in\bbN$ sufficiently large, $H_n$ has a simple  eigenvalue
$\lambda_n$ in a small neighborhood of $z=0$ satisfying
\begin{equation}
\lambda_n \underset{n\to\infty}{\longrightarrow} 0.
\end{equation}
We denote by $\Psi_n$ the corresponding eigenfunction of $H_n$,
associated with the eigenvalue $\lambda_n$ of $H_n$, $H_n
\Psi_n=\lambda_n \Psi_n$, $n\in\bbN$. $\Psi_n$ is then a constant
multiple of the solutions satisfying \eqref{2.51} with $z=\lambda_n$
and $V$ replaced by $V_n$. We may choose the constant multiple in
$\Psi_n$ such that
\begin{equation}
\lim_{n\to \infty}\|\Psi_n - \Psi_0\|_{L^2(\bbR;dx)}=0.
\end{equation}
In addition, we recall that \eqref{2.51} also implies that $\Psi_n$
and $\Psi_0$ are exponentially bounded in $x\in\bbR$ with bounds
uniform with  respect to $n\in\bbN$.

In addition, we abbreviate
\begin{align}
& \Phi_n = u_n \Psi_n, \quad \wti \Phi_n = C_n S_n \Phi_n, \quad
\big(\wti \Phi_n, \Phi_n\big)_{L^2(\bbR;dx)} =1,   \\
& P_n = \big(\wti \Phi_n, \dott \big) \Phi_n,   \\
& S_n f = \sgn(W_n) f, \quad f \in L^2(\bbR; dx),  \\
& C_n = \bigg[\int_{\bbR} dx \, W_n(x) |\Psi_n(x)|^2 \bigg]^{-1},  \\
& K_n(z) \underset{z\to\lambda_n}{=} K_{0,n} + K_{1,n}(z-\lambda_n) +
\Oh\big((z-\lambda_n)^2\big), \\
& K_{0,n} = - u_n \big(H^{(0)}-\lambda_n I_{L^2(\bbR;dx)}\big)^{-1} v_n,
\quad K_{1,n} = - u_n \big(H^{(0)}-\lambda_n
I_{L^2(\bbR;dx)}\big)^{-2} v_n,  \no \\
& \hspace*{10cm} n\in\bbN,
\end{align}
and recall that
\begin{align}
& - u_n \big(H^{(0)}-\lambda_n I_{L^2(\bbR;dx)}\big)^{-1} v_n \Phi_n
= \Phi_n, \quad
- v_n \big(H^{(0)}-\lambda_n I_{L^2(\bbR;dx)}\big)^{-1} u_n \wti
\Phi_n = \wti \Phi_n,    \\
& H_n \Psi_n = \big(H^{(0)} + W_n\big) \Psi_n = \lambda_n \Psi_n, \quad
- \big(H^{(0)}-\lambda_n I_{L^2(\bbR;dx)}\big)^{-1} W_n \Psi_n =
\Psi_n, \; n\in\bbN.
\end{align}
We note that since $u_n, v_n, W_n$, $n\in\bbN$, are all bounded operators on
$L^2(\bbR; dx)$, no operator closure symbols in
$u_n \big(H^{(0)}-\lambda_n I_{L^2(\bbR;dx)}\big)^{-k} v_n$, $k=1,2$,
are needed in the following computation leading up to \eqref{3.80a}.

Next, one computes
\begin{align}
& {\det}_{P_n L^2(\bbR;dx)}(P_n K_{1,n} P_n)
= - {\det}_{P_n L^2(\bbR;dx)}\big(P_n u_n
\big(H^{(0)}-\lambda_n I_{L^2(\bbR;dx)}\big)^{-2} v_n P_n\big) \no \\
& \quad = - \big(\wti \Phi_n, u_n \big(H^{(0)}-\lambda_n
I_{L^2(\bbR;dx)}\big)^{-2}
v_n \Phi_n\big)_{L^2(\bbR;dx)}  \no \\
& \quad = - \big(\big(H^{(0)}-\lambda_n I_{L^2(\bbR;dx)}\big)^{-1}
u_n \wti\Phi_n,
\big(H^{(0)}-\lambda_n I_{L^2(\bbR;dx)}\big)^{-1} v_n
\Phi_n\big)_{L^2(\bbR;dx)}  \no \\
& \quad = - C_n \big(\big(H^{(0)}-\lambda_n
I_{L^2(\bbR;dx)}\big)^{-1} W_n \Psi_n,
\big(H^{(0)}-\lambda_n I_{L^2(\bbR;dx)}\big)^{-1} W_n
\Psi_n\big)_{L^2(\bbR;dx)}  \no \\
& \quad = - C_n \big\|\Psi_n\big\|^2_{L^2(\bbR;dx)}, \quad n\in\bbN.
\lb{3.80a}
\end{align}
Since
\begin{equation}
\lim_{n\to\infty} C_n = C_0 =\bigg[\int_{\bbR} dx\, W(x)
|\Psi_0(x)|^2 \bigg]^{-1},
\end{equation}
one finally obtains,
\begin{align}
{\det}_{P_0 L^2(\bbR;dx)}(P_0K_1P_0) &
= - {\det}_{P_0 L^2(\bbR;dx)}\Big(P_0 \ol{u{H^{(0)}}^{-2}v} P_0\Big) \no \\
&= - \Big(\wti\Phi_0, \ol{u{H^{(0)}}^{-2}v} \Phi_0\Big)_{L^2(\bbR;dx)}  \no \\
&= - \lim_{n\to\infty} \big(\wti \Phi_n, u_n
\big(H^{(0)}-\lambda_n I_{L^2(\bbR;dx)}\big)^{-2} v_n
\Phi_n\big)_{L^2(\bbR;dx)}  \no \\
&= - \lim_{n\to\infty} C_n \big\|\Psi_n\big\|^2_{L^2(\bbR;dx)} \no \\
& = - C_0 \|\Psi_0\|_{L^2(\bbR;dx)}^2    \lb{2.58}  \\
& = - \bigg[\int_{\bbR} dx\, W(x) |\Psi_0(x)|^2 \bigg]^{-1}
\|\Psi_0\|_{L^2(\bbR;dx)}^2   \lb{2.59}  \\
& = - \|\Psi_0\|^2_{L^2(\bbR; dx)} \big/ (v\Psi_0,
u\Psi_0)_{L^2(\bbR; dx)}  \lb{2.60a} \\
& = \big[\|\Psi_0'\|_{L^2(\bbR;dx)}^2 + V_{\infty}
\|\Psi_0\|_{L^2(\bbR;dx)}^2\big]^{-1}
\|\Psi_0\|_{L^2(\bbR;dx)}^2  >  0.   \lb{2.61a}
\end{align}
Here we applied the quadratic form equality
\begin{align}
\begin{split}
0 & < \|\Psi_0'\|_{L^2(\bbR;dx)}^2 + V_{\infty} \|\Psi_0\|_{L^2(\bbR;dx)}^2
=  - (v\Psi_0, u \Psi_0)_{L^2(\bbR;dx)}  \\
& = - \int_{\bbR} dx\, W(x) |\Psi_0(x)|^2
\end{split}
\end{align}
to \eqref{2.59}, to arrive at \eqref{2.61a}.
\end{proof}

\subsubsection{The first factor: A posteriori computation}
Before we proceed to a direct approach to compute the first factor on
the right-hand side of \eqref{2.55},
\begin{equation}
{\det}_{L^2(\bbR;dx)}(I_{L^2(\bbR;dx)}-K_0-P_0),
\end{equation}
we will next determine
${\det}_{L^2(\bbR;dx)}(I_{L^2(\bbR;dx)}-K_0-P_0)$ by using the final
answer \eqref{2.54b} for $\cF^{\bullet}(0)$.

\begin{theorem}  \lb{t3.6}
Assume Hypothesis \ref{h3.1} and suppose $0 \in \sigma_{\rm d}(H)$. Then,
\begin{align}
&{\det}_{L^2(\bbR;dx)}(I_{L^2(\bbR;dx)}-P_0-K_0)  \no \\
& \quad = \f{1}{2 V_{\infty}^{1/2}} \int_{\bbR} dx \, [V(x) - V_{\infty}]
\Psi_+(0,x) \Psi_-(0,x)  \lb{3.54}  \\
& \quad = - \big[2 V_{\infty}^{1/2}\big]^{-1}
\|\Psi_{\pm}(0)\|_{L^2(\bbR;dx)}^{-2}
\big[\|\Psi'_{\pm}(0)\|_{L^2(\bbR;dx)}^2 + V_{\infty}
\|\Psi_{\pm}(0)\|_{L^2(\bbR;dx)}^2\big]
\no  \\
& \qquad \times (\Psi_-(0),\Psi_+(0))_{L^2(\bbR;dx)}.   \lb{3.55}
\end{align}
$($Here the equations for the $+$ and $-$ sign should be read separately.$)$
\end{theorem}
\begin{proof}
Combining \eqref{2.52}, \eqref{2.54b} (setting $\lambda_j=0$), \eqref{2.55},
\eqref{2.60}, and \eqref{2.58},
and taking into account that for some constants
$C_{\pm}\in\bbR\backslash\{0\}$,
\begin{equation}
\Psi_0(x) = C_{\pm} \Psi_{\pm}(0,x), \quad x\in\bbR,   \lb{2.79}
\end{equation}
then yields for the first factor in \eqref{2.55},
\begin{align}
&{\det}_{L^2(\bbR;dx)}(I_{L^2(\bbR;dx)}-P_0-K_0)
= \f{\cF^{\bullet}(0)}{{\det}_{P_0 L^2(\bbR;dx)}(P_0K_1P_0)}  \no \\
& \quad = C_0^{-1} \|\Psi_0\|^{-2} \f{1}{2 V_{\infty}^{1/2}}
\int_{\bbR} dx' \, \Psi_-(0,x') \Psi_+(0,x')  \no \\
& \quad = \big[2 V_{\infty}^{1/2}\big]^{-1} \|\Psi_{0}\|^{-2}_{L^2(\bbR;dx)}
(v\Psi_{0},u\Psi_{0})_{L^2(\bbR;dx)}
(\Psi_-(0),\Psi_+(0))_{L^2(\bbR;dx)}  \no \\
&\quad = \big[2 V_{\infty}^{1/2}\big]^{-1}
\|\Psi_{\pm}(0)\|^{-2}_{L^2(\bbR;dx)}
(v\Psi_{\pm}(0),u\Psi_{\pm}(0))_{L^2(\bbR;dx)}
(\Psi_-(0),\Psi_+(0))_{L^2(\bbR;dx)}
\no \\
& \quad = \f{1}{2 V_{\infty}^{1/2}} \int_{\bbR} dx \, [V(z) - V_{\infty}]
\Psi_+(0,x) \Psi_-(0,x)   \no \\
& \quad = - \big[2 V_{\infty}^{1/2}\big]^{-1}
\|\Psi_{\pm}(0)\|_{L^2(\bbR;dx)}^{-2}
\big[\|\Psi'_{\pm}(0)\|_{L^2(\bbR;dx)}^2 + V_{\infty}
\|\Psi_{\pm}(0)\|_{L^2(\bbR;dx)}^2\big]
\no \\
& \qquad \times (\Psi_-(0),\Psi_+(0))_{L^2(\bbR;dx)}.
\end{align}
\end{proof}

\subsubsection{The first factor: Direct computation}
Next, we proceed to a direct approach to compute the first factor on
the right-hand side of \eqref{2.55},
${\det}_{L^2(\bbR;dx)}(I_{L^2(\bbR;dx)}-K_0-P_0)$. This will now be
an {\it ab initio} calculation entirely independent of the result
\eqref{2.80}.

\begin{theorem}  \lb{t3.6a}
Assume Hypothesis \ref{h3.1} and suppose $0 \in \sigma_{\rm d}(H)$. Then,
\begin{equation}
{\det}_{L^2(\bbR;dx)}(I_{L^2(\bbR;dx)}-P_0-K_0)
= \f{-1}{2 V_{\infty}^{1/2}} \int_{\bbR} dx \, [V(x) - V_{\infty}]
e^{\pm V_{\infty}^{1/2}x} \psi_{\pm}(x),   \lb{3.71}
\end{equation}
where $\psi_{\pm}$ are defined by
\begin{align}
& \psi_{\pm}(x)= - \Psi_{\pm}(0,x)
- \f{1}{V_\infty^{1/2}} \int_x^{\pm \infty} dx' \,
\sinh\big(V_\infty^{1/2}(x-x')\big) [V(x')-V_\infty] \psi_{\pm}(x'),   \no \\
& \hspace*{10cm}  x\in\bbR.   \lb{3.72}
\end{align}
$($Here the equations for the $+$ and $-$ sign should be read separately.$)$
\end{theorem}
\begin{proof}
Our strategy is to apply formulas (3.9) and (3.12) in \cite{GM04} to
the Fredholm determinant
${\det}_{L^2(\bbR;dx)}(I_{L^2(\bbR;dx)}-P_0-K_0)$ by appealing to the
semi-separable nature of $P_0+K_0$ upon an appropriate reformulation
involving block operator matrices. To this end we introduce
\begin{align}
f_1(x) &= \begin{pmatrix} -u(x) e^{-V_\infty^{1/2}x} & u(x) \Psi_0(x)
\end{pmatrix}, \\
g_1(x) &= \begin{pmatrix} [2(V_\infty)^{1/2}]^{-1}v(x) e^{V_\infty^{1/2}x} &
C_0 v(x) \Psi_0(x) \end{pmatrix}^{\top},  \\
f_2(x) &= \begin{pmatrix} -u(x) e^{V_\infty^{1/2}x} & u(x) \Psi_0(x)
\end{pmatrix}, \\
g_2(x) &= \begin{pmatrix} [2(V_\infty)^{1/2}]^{-1}v(x) e^{-V_\infty^{1/2}x} &
C_0 v(x) \Psi_0(x) \end{pmatrix}^{\top}
\end{align}
and note that $P_0+K_0$ is an integral operator with semi-separable
integral kernel
\begin{equation}
(P_0 + K_0)(x,x')=\begin{cases} f_1(x) g_1(x'), & x'<x, \\
f_2(x) g_2(x'), & x'>x. \end{cases}
\end{equation}
In addition, we introduce the integral kernel
\begin{equation}
H(x,x')=f_1(x)g_1(x')-f_2(x)g_2(x')=u(x)
\f{\sinh\big(V_\infty^{1/2}(x-x')\big)}{V_\infty^{1/2}} v(x'),
\end{equation}
and, temporarily assuming that
\begin{equation}
\supp \, (V - V_\infty) \, \text{ is compact,}    \lb{3.79}
\end{equation}
the pair of Volterra integral equations
\begin{align}
\hat f_1 (x) &= f_1(x) - \int_x^{\infty} dx' \, H(x,x') \hat f_1(x'),
\lb{3.80} \\
\hat f_2 (x) &= f_2(x) + \int_{-\infty}^x dx' \, H(x,x') \hat f_2(x')
\lb{3.81}
\end{align}
for a.e.\ $x\in\bbR$. Applying Theorem 3.2 of \cite{GM04}
(especially, (3.9) and (3.12) in \cite{GM04}) one then infers
\begin{align}
{\det}_{L^2(\bbR;dx)}(I_{L^2(\bbR;dx)}-P_0-K_0)
& = { \det}_{\bbC^2}\bigg(I_2-\int_{\bbR} dx \, g_1(x) \hat
f_1(x)\bigg)    \lb{3.82}  \\
& = { \det}_{\bbC^2}\bigg(I_2-\int_{\bbR} dx \, g_2(x) \hat
f_2(x)\bigg),   \lb{3.83}
\end{align}
with $I_2$ the identity matrix in $\bbC^2$. Introducing
$\hatt\Psi_\pm$ as the solutions of the pair of Volterra integral
equations
\begin{align}
\begin{split}
\hatt \Psi_{\pm}(x) &= \begin{pmatrix} e^{\mp V_\infty^{1/2}x} &
-\Psi_0(x) \end{pmatrix}
  \\
& \quad - \f{1}{V_\infty^{1/2}} \int_x^{\pm \infty} dx' \,
\sinh\big(V_\infty^{1/2}(x-x')\big) [V(x')-V_\infty]
\hatt \Psi_{\pm}(x'), \quad x\in\bbR,   \lb{3.84}
\end{split}
\end{align}
a comparison with \eqref{3.80} and \eqref{3.81} yields
\begin{equation}
\hat f_1(x)=-u(x)\hatt \Psi_+(x), \; \hat f_2(x)=-u(x)\hatt \Psi_-(x)
\, \text{ for a.e.\
$x\in\bbR$}.
\end{equation}
Writing
\begin{equation}
\hatt \Psi_\pm (x) = \begin{pmatrix} \Psi_\pm(0,x) & \hat \psi_\pm(x)
\end{pmatrix},
\quad x\in\bbR,
\end{equation}
\eqref{3.84} yields ($x\in\bbR$)
\begin{equation}
\hat \psi_{\pm}(x)= - \Psi_0 (x)
- \f{1}{V_\infty^{1/2}} \int_x^{\pm \infty} dx' \,
\sinh\big(V_\infty^{1/2}(x-x')\big) [V(x')-V_\infty]
\hat \psi_{\pm}(x'),
\end{equation}
and
\begin{equation}
\Psi_{\pm}(0,x)= e^{\mp V_\infty^{1/2}x}
- \f{1}{V_\infty^{1/2}} \int_x^{\pm \infty} dx' \,
\sinh\big(V_\infty^{1/2}(x-x')\big) [V(x')-V_\infty]
\Psi_{\pm}(0,x'),    \lb{3.88}
\end{equation}
in accordance with \eqref{2.51} for $z=0$. Because of \eqref{2.79},
$\Psi_0(x)=C_\pm \Psi_\pm(0,x)$, $x\in\bbR$, one infers that
\begin{equation}
\psi_\pm(x) = C_\pm^{-1} \hat \psi_\pm (x), \quad x\in\bbR,   \lb{3.89}
\end{equation}
with $\psi_\pm$ satisfying \eqref{3.72}.

Next, one computes
\begin{align}
g_1(x)\hat f_1(x)&=\left(\begin{smallmatrix}
- \f{\exp\big(V_\infty^{1/2}x\big)}{2V_\infty^{1/2}} [V(x)-V_\infty]
\Psi_+(0,x)
& - \f{\exp\big(V_\infty^{1/2}x\big)}{2V_\infty^{1/2}}
[V(x)-V_\infty] \hat \psi_+(x)  \\[1mm]
-C_0 [V(x)-V_\infty] \Psi_0(x)\Psi_+(0,x)  &
-C_0 [V(x)-V_\infty] \Psi_0(x)\hat \psi_+(x)
\end{smallmatrix}\right),  \lb{3.90} \\
g_2(x)\hat f_2(x)&=\left(\begin{smallmatrix}
- \f{\exp\big(-V_\infty^{1/2}x\big)}{2V_\infty^{1/2}} [V(x)-V_\infty]
\Psi_-(0,x)
& - \f{\exp\big(-V_\infty^{1/2}x\big)}{2V_\infty^{1/2}}
[V(x)-V_\infty] \hat \psi_-(x)  \\[1mm]
-C_0 [V(x)-V_\infty] \Psi_0(x)\Psi_-(0,x)  &
-C_0 [V(x)-V_\infty] \Psi_0(x)\hat \psi_-(x)
\end{smallmatrix}\right).   \lb{3.91}
\end{align}
Using the fact that $0\in\sigma_{\rm d}(H)$, and hence
${\det}_{L^2(\bbR;dx)}(I_{L^2(\bbR;dx)}-K_0)=0$, one obtains from
taking $z=0$ in \eqref{2.54} and \eqref{2.54a},
\begin{equation}
{\det}_{L^2(\bbR;dx)}(I_{L^2(\bbR;dx)}-K_0)
=1 + \f{1}{2V_\infty^{1/2}} \int_{\bbR} dx \, e^{\pm V_\infty^{1/2}x}
[V(x)-V_\infty] \Psi_\pm (0,x) =0.   \lb{3.92}
\end{equation}
Thus, \eqref{3.82} and \eqref{3.83} together with
\eqref{3.90}--\eqref{3.92} yield
\begin{align}
{\det}_{L^2(\bbR;dx)}(I_{L^2(\bbR;dx)}-P_0-K_0) & =
- C_0 \int_{\bbR} dx \, [V(x)-V_\infty] \Psi_0(x) \Psi_\pm (0,x)  \\
& \quad \times
\f{1}{2V_\infty^{1/2}} \int_{\bbR} dx \, e^{\pm V_\infty^{1/2}x}
[V(x)-V_\infty]
\hat \psi_\pm(x)  \no
\end{align}
(where the equations for the $+$ and $-$ sign should be read
separately). Applying
\eqref{3.22} for $\lambda_j=0$ and \eqref{2.79} one obtains
\begin{align}
{\det}_{L^2(\bbR;dx)}(I_{L^2(\bbR;dx)}-P_0-K_0) & =
\f{-1}{2V_\infty^{1/2}C_\pm} \int_{\bbR} dx \, e^{\pm
V_\infty^{1/2}x} [V(x)-V_\infty]
\hat \psi_\pm(x)  \no \\
& = \f{-1}{2V_\infty^{1/2}} \int_{\bbR} dx \, e^{\pm V_\infty^{1/2}x}
[V(x)-V_\infty]
\psi_\pm(x),
\end{align}
using \eqref{3.89} in the last line.

To remove the temporary compact support assumption \eqref{3.79} we
first note that
by a standard iteration argument, the Volterra equations
\begin{align}
e^{\pm V_\infty^{1/2}x} \psi_{\pm}(x) &= - e^{\pm V_\infty^{1/2}x}
\Psi_{\pm}(0,x)
\mp \f{1}{2V_\infty^{1/2}} \int_x^{\pm \infty} dx' \,
\big[e^{\pm 2V_\infty^{1/2}(x-x')} -1\big]  \no \\
& \hspace*{2.9cm} \times [V(x')-V_\infty]
e^{\pm V_\infty^{1/2}x'} \psi_{\pm}(x'), \quad x\in\bbR,   \lb{3.95}
\end{align}
have unique and bounded solutions on $\bbR$, which together with
their first derivatives are locally absolutely continuous on $\bbR$,
as long as the condition \eqref{3.1} is satisfied. This follows since
there exists a constant $C>0$ such that
$\big|e^{\pm V_\infty^{1/2}x} \Psi_{\pm}(0,x)\big|\leq C$,
$x\in\bbR$. Thus, the right-hand side of \eqref{3.71} remains
well-defined under condition \eqref{3.1} on $V$.

Next, similarly to the proof of Theorem \ref{t3.4}, we choose compactly
supported sequences $u_n, v_n \in L^2(\bbR; dx)$,
$n\in\bbN$, such that $W_n=V_n-V_\infty=u_n v_n$ converges to
$W=V-V_\infty=uv$ in $L^1(\bbR;dx)$ as $n\to\infty$ and introduce the
maximally defined operator $H_n$ in $L^2(\bbR;dx)$ associated with
the differential expression $L_n=-(d^2/d x^2) + V_n(x)$,
$x\in\bbR$. Since $H_n$ converges to $H$ in norm resolvent sense
(this follows in exactly the same manner as discussed in the proof of
Theorem \ref{t3.4}), the spectrum of $H_n$ converges to that of $H$ as
$n\to\infty$. In particular, for $n\in\bbN$ sufficiently large,
$H_n$ has a simple  eigenvalue $\lambda_n$ in a small neighborhood of
$z=0$ such that
$\lambda_n\underset{n\to\infty}{\longrightarrow} 0$ . Multiplying
$V_n$ by a suitable coupling constant $g_n\in\bbR$, where $g_n
\underset{n\to \infty}{\longrightarrow} 1$, then guarantees that the
maximally defined operator $H_n(g_n)$ in $L^2(\bbR;dx)$ associated
with the differential expression $L_n(g_n)=- (d^2/d x^2) + g_n V_n(x)$,
$x\in\bbR$, has a simple eigenvalue at $z=0$, in particular,
$0\in\sigma_{\rm d}(H_n(g_n))$. (Multiplying $V$ by $g_n$ changes the
essential spectrum of $H_n(g_n)$ into $[g_nV_\infty, \infty)$, but
since $\lambda_n\to 0$ and $g_n\to 1$ as $n\to\infty$, this shift in
the essential spectrum is irrelevant in this proof as long as
$n\in\bbN$ is sufficiently large.)

Finally, the approximation arguments described in the proof of
Theorem 4.3 of \cite{GM04} permit one to pass to the limit $n\to
\infty$ establishing \eqref{3.71} without the extra hypothesis
\eqref{3.79}.
\end{proof}

It remains to show that the two results \eqref{3.54} and \eqref{3.71}
for the Fredholm determinant
${\det}_{L^2(\bbR;dx)}(I_{L^2(\bbR;dx)}-P_0-K_0)$ coincide. This will
be undertaken next.

\begin{theorem}  \lb{l3.9}
Assume Hypothesis \ref{h3.1} and suppose $0 \in \sigma_{\rm d}(H)$. Then the
expressions \eqref{3.54} and \eqref{3.71} for
${\det}_{L^2(\bbR;dx)}(I_{L^2(\bbR;dx)}-P_0-K_0)$ coincide.
\end{theorem}
\begin{proof}
To keep the arguments as short as possible, we first prove that
\eqref{3.54} and
\eqref{3.71} coincide under the simplifying compact support
assumption \eqref{3.79} on
$V-V_\infty$. Again, the general case where $V$ satisfies Hypothesis
\ref{h3.1} then follows from an approximation argument.

More precisely, we suppose that
\begin{equation}
\supp \, (V-V_\infty) \subset [-R,R] \, \text{ for some $R>0$.}
\end{equation}
Proving that \eqref{3.54} and \eqref{3.71} coincide is then
equivalent to showing that
\begin{equation}
- \int_{-R}^{R} dx \, [V(x)-V_\infty] e^{\pm V_\infty^{1/2}x} \psi_\pm(x)
= \int_{-R}^{R} dx \, [V(x)-V_\infty] \Psi_+(0,x) \Psi_-(0,x).   \lb{3.97}
\end{equation}

We start with the right-hand side of \eqref{3.97}: First we note that
\begin{equation}
\f{d}{dx} \Wr(\psi_\pm,\Psi_\mp(0))(x) = [V(x)-V_\infty]\Psi_+(0,x)\Psi_-(0,x) \,
\text{ for a.e.\ $x\in\bbR$,}  \lb{3.98}
\end{equation}
where we used that
\begin{equation}
- \psi_{\pm}''(x) + V(x) \psi_{\pm}(x) = [V(x)-V_\infty] \Psi_{\pm}(0,x) \,
\text{ for a.e.\ $x\in\bbR$,}   \lb{3.98a}
\end{equation}
which in turn follows by twice differentiating \eqref{3.72}, and
\begin{equation}
- \Psi_{\pm}''(0,x) + V(x) \Psi_{\pm}(0,x) = 0 \,
\text{ for a.e.\ $x\in\bbR$.}   \lb{3.98b}
\end{equation}
Thus, one concludes that
\begin{align}
& \int_{-R}^{R} dx \, [V(x)-V_\infty] \Psi_+(0,x) \Psi_-(0,x) \no \\
& \quad =  \int_{-R}^{R} dx \, \f{d}{dx} \Wr(\psi_+,\Psi_-(0))(x) \no \\
& \quad = \Wr(\psi_+,\Psi_-(0))(R) - \Wr(\psi_+,\Psi_-(0))(-R) \no \\
& \quad = - \Wr(\Psi_+(0),\Psi_-(0))(R) - \Wr(\psi_+,\Psi_-(0))(-R) \no \\
& \quad = - \Wr(\psi_+,\Psi_-(0))(-R).   \lb{3.99}
\end{align}
Here we employed that
\begin{equation}
\psi_+(x) = - \Psi_+(0,x) \, \text{ for $x \geq R$}    \lb{3.100}
\end{equation}
(cf.\ \eqref{3.72}) and
\begin{equation}
2 V_\infty^{1/2} \cF(0)= \Wr(\Psi_+(0),\Psi_-(0)) = 0   \lb{3.101}
\end{equation}
since by hypothesis, $0\in\sigma_{\rm d}(H)$. Similarly, one computes
\begin{align}
& \int_{-R}^{R} dx \, [V(x)-V_\infty] \Psi_+(0,x) \Psi_-(0,x) \no \\
& \quad =  \int_{-R}^{R} dx \, \f{d}{dx} \Wr(\psi_-,\Psi_+(0))(x) \no \\
& \quad = \Wr(\psi_-,\Psi_+(0))(R) - \Wr(\psi_-,\Psi_+(0))(-R) \no \\
& \quad = \Wr(\psi_-,\Psi_+(0))(R) + \Wr(\Psi_-(0),\Psi_+(0))(-R) \no \\
& \quad = \Wr(\psi_-,\Psi_+(0))(R),   \lb{3.101a}
\end{align}
where we used
\begin{equation}
\psi_-(x) = - \Psi_-(0,x) \, \text{ for $x \leq -R$}    \lb{3.101b}
\end{equation}
(cf.\ \eqref{3.72}) and again \eqref{3.101}. In particular, one concludes that
\begin{equation}
\Wr(\psi_-,\Psi_+(0))(R) = - \Wr(\psi_+,\Psi_-(0))(-R).   \lb{3.101c}
\end{equation}

To compute the left-hand side of \eqref{3.97} we first note that
\begin{align}
& \f{d}{dx} \Wr\Big(\psi_\pm(x), e^{\pm V_\infty^{1/2}x}\Big) +
\f{d}{dx} \Wr\Big(\Psi_\pm(0,x), e^{\pm V_\infty^{1/2}x}\Big) \no \\
& \quad = - [V(x)-V_\infty] e^{\pm V_\infty^{1/2}x} \psi_{\pm}(x)  \,
\text{ for a.e.\ $x\in\bbR$,}   \lb{3.102}
\end{align}
where we employed again \eqref{3.98a} and \eqref{3.98b}. Thus, one infers that
\begin{align}
& - \int_{-R}^{R} dx \, [V(x)-V_\infty] e^{V_\infty^{1/2}x} \psi_+(x)   \no \\
& \quad =\int_{-R}^{R} dx \, \bigg[\f{d}{dx} \Wr\Big(\psi_+(x),
e^{V_\infty^{1/2}x}\Big) + \f{d}{dx} \Wr\Big(\Psi_+(0,x),
e^{V_\infty^{1/2}x}\Big)\bigg]  \no \\
& \quad = \Wr\Big(\psi_+(x),e^{V_\infty^{1/2}x}\Big)\Big|_{x=R}
- \Wr\Big(\psi_+(x),e^{V_\infty^{1/2}x}\Big)\Big|_{x=-R}  \no \\
& \qquad + \Wr\Big(\Psi_+(0,x),e^{V_\infty^{1/2}x}\Big)\Big|_{x=R}
- \Wr\Big(\Psi_+(0,x),e^{V_\infty^{1/2}x}\Big)\Big|_{x=-R}   \no \\
& \quad = - \Wr\Big(\Psi_+(0,x),e^{V_\infty^{1/2}x}\Big)\Big|_{x=R}
+ \Wr\Big(\Psi_+(0,x),e^{V_\infty^{1/2}x}\Big)\Big|_{x=R}   \no \\
& \qquad - \Wr\Big(\psi_+(x),e^{V_\infty^{1/2}x}\Big)\Big|_{x=-R}
- \Wr\Big(\Psi_+(0,x),e^{V_\infty^{1/2}x}\Big)\Big|_{x=-R}   \no \\
& \quad = - \Wr(\psi_+,\Psi_-(0))(-R)
- \Wr(\Psi_+(0),\Psi_-(0))(-R)   \no \\
& \quad = - \Wr(\psi_+,\Psi_-(0))(-R).
\end{align}
Here we used again \eqref{3.100} and \eqref{3.101} as well as
\begin{equation}
e^{V_\infty^{1/2}x} = \Psi_-(0,x) \, \text{ for $x\leq -R$}    \lb{3.105}
\end{equation}
(cf.\ \eqref{3.88}). Similarly one computes
\begin{equation}
- \int_{-R}^{R} dx \, [V(x)-V_\infty] e^{-V_\infty^{1/2}x} \psi_-(x)
  = \Wr(\psi_-,\Psi_+(0))(R).
\end{equation}
Taking into account \eqref{3.101c}, this completes the proof of \eqref{3.97}.
\end{proof}

\subsection{A formula of Simon}\lb{simonform}
Finally, we turn to an interesting formula for the Jost solutions
$\Psi_\pm(z,\dott)$ in terms of Fredholm determinants derived by
Simon \cite{Si00}.

To set the stage, we abbreviate $\bbR_{\pm}=(0,\pm\infty)$ and
introduce the one-dimensional Dirichlet and Neumann Laplacians
perturbed by the constant potential
$V_\infty$, $H_{\pm}^{(0),D}$ and $H_{\pm}^{(0),N}$  in $L^2(\bbR_{\pm};dx)$ by
\begin{align}
&H_{\pm}^{(0),D} = - \f{d^2}{d x^2} + V_\infty,  \quad
\dom\big(H_{\pm}^{(0),D}\big)=\{g\in H^{2}(\bbR_{\pm}) \,|\, g(0_{\pm})=0\}, \\
&H_{\pm}^{(0),N} = - \f{d^2}{d x^2} + V_\infty,  \quad
\dom\big(H_{\pm}^{(0),N}\big)=\{g\in H^{2}(\bbR_{\pm}) \,|\, g'(0_{\pm})=0\}.
\end{align}

Next, we recall that
\begin{equation}
\Psi_{\pm}(z,0)={\det}_{L^2(\bbR_{\pm};dx)}\Big(I+u(H_{\pm}^{(0),D} -
z)^{-1}v\Big), \quad
\Im\big((z-V_\infty)^{1/2}\big)
>  0,   \lb{3.58}
\end{equation}
a celebrated formula by Jost and Pais \cite{JP51} (in the case
$V_\infty =0$). For more details and background on \eqref{3.58} we
refer to \cite{GM04}
and the references cited therein. Moreover, it is known (cf.\
\cite{GLMZ05}, \cite{GM04}) that
\begin{align}
\begin{split}
\Psi_{\pm}^{\prime}(z,0) &= \pm i (z-V_\infty)^{1/2}
{\det}_{L^2(\bbR_{\pm};dx)}\Big(I+u(H_{\pm}^{(0),N}-z)^{-1}v\Big),   \\
&  \hspace*{5.05cm}  \Im\big((z-V_\infty)^{1/2}\big) > 0.  \lb{3.59}
\end{split}
\end{align}

We conclude this section by presenting a quick proof of the
representation of the Jost solutions $\Psi_{\pm}(z,x)$ and their
$x$-derivatives, $\Psi_{\pm}^{\prime}(z,x)$, in terms of symmetrized
perturbation determinants, starting from the Jost and Pais formula
\eqref{3.58} and its analog \eqref{3.59} for
$\Psi_{\pm}^{\prime}(z,0)$:

\begin{lemma} [\cite{Si00}]  \lb{l3.7}
Suppose $V$ satisfies \eqref{3.1} and let
$\Im\big((z-V_\infty)^{1/2}\big) > 0$,
$x\in\bbR$. Then,
\begin{align}
\Psi_{\pm}(z,x)&=e^{\pm i (z-V_\infty)^{1/2} x}  \no \\
& \quad \times
{\det}_{L^2(\bbR_{\pm};dx)}\Big(I_{L^2(\bbR_{\pm};dx))}
+u(\dott +x)(H_{\pm}^{(0),D} - z)^{-1}v(\dott +x)\Big),    \lb{3.60} \\
\Psi_{\pm}^{\prime}(z,x)&= \pm i (z-V_\infty)^{1/2} e^{\pm
i(z-V_\infty)^{1/2} x}  \no \\
& \quad \times
{\det}_{L^2(\bbR_{\pm};dx)}\Big(I_{L^2(\bbR_{\pm};dx))}
+u(\dott+x)(H_{\pm}^{(0),N}-z)^{-1}v(\dott+x)\Big).
\lb{3.61}
\end{align}
\end{lemma}
\begin{proof}
Denoting $V_y(x)=V(x+y)$, $x, y\in\bbR$, and by
$\Psi_{y,\pm}(z,\dott)$ the Jost solutions associated with $V_y$, an
elementary change of variables in the Volterra integral equation
\eqref{2.51} for $\Psi_{y,\pm}$ yields
\begin{align}
\begin{split}
\Psi_{y,\pm}(z,x) &= e^{\mp i (z-V_\infty)^{1/2}y} \Psi_{\pm}(z,x+y), \\
\Psi_{y,\pm}^{\prime}(z,x) &= e^{\mp i (z-V_\infty)^{1/2}y}
\Psi_{\pm}^{\prime}(z,x+y).   \lb{3.62}
\end{split}
\end{align}
Taking $x=0$ in \eqref{3.62} implies
\begin{align}
\Psi_{y,\pm}(z,0) &= e^{\mp i (z-V_\infty)^{1/2}y} \Psi_{\pm}(z,y),
\lb{3.63}  \\
\Psi_{y,\pm}^{\prime}(z,0) &= e^{\mp i (z-V_\infty)^{1/2}y}
\Psi_{\pm}^{\prime}(z,y).
\lb{3.64}
\end{align}
Using the Jost--Pais-type formulas
\begin{align}
\Psi_{y,\pm}(z,0)&={\det}_{L^2(\bbR_{\pm};dx)}\Big(I+u(\dott +y)
(H_{\pm}^{(0),D} - z)^{-1}v(\dott +y)\Big),
\lb{3.65} \\
\Psi_{y,\pm}^{\prime}(z,0)&= \pm i (z-V_\infty)^{1/2}
{\det}_{L^2(\bbR_{\pm};dx)}\Big(I+u(\dott+y)(H_{\pm}^{(0),N}-z)^{-1}v(\dott+y)\Big),
\lb{3.66}
\end{align}
an insertion of \eqref{3.65} into the left-hand side of \eqref{3.63}
proves \eqref{3.60}. Similarly, an insertion of \eqref{3.66} into the
left-hand side of \eqref{3.64} yields
\eqref{3.61}.
\end{proof}

\section{The multi-dimensional case}\label{s4}

In the previous section, we have illustrated
within the simple setting of one-dimensional scalar
reaction--diffusion equations how the stability index may
be equally well calculated from a Jost/Evans function
point of view, or else, using semi-separability
of the integral kernels of Birman--Schwinger-type operators,
directly from first principles using Fredholm determinants.
We conclude by describing, again within the reaction--diffusion
setting, an algorithm for multi-dimensional computations via
Fredholm determinants,
based on semi-separability of the
integral kernels combined with Galerkin approximations.

\subsection{Flow in an infinite cylinder}\label{cyl}
Consider a scalar reaction-diffusion equation
\begin{equation}\label{multeq}
w_t =\Delta w + f(w),
\end{equation}
on an infinite cylinder $x=(x_1,x_2,\dots,x_d)
\in \bbR\times \Omega$, where $\Delta=\Delta_x$ is the
Laplacian in the $x$-variables, $\Omega\subset \bbR^{d-1}$
is a bounded domain, $w$ and $f$ are real-valued functions,
\begin{equation}
f\in C^{3}(\bbR).
\end{equation}
In what follows
we will assume that $\Omega=[0,2\pi]^{d-1}$ and consider
only the physical cases $d=2,3$. Unless explicitly stated otherwise, we will 
always assume that periodic boundary conditions are used on the 
boundary $\partial\Omega$ of $\Omega$ (viewing $\Omega$ as a 
$(d-1)$-dimensional torus in the following) if $d=2,3$. For $x\in\bbR\times\Omega$
we will always write $x=(x_1,y)$, where $x_1\in\bbR$ and
$y=(x_2,\dots,x_d)\in\Omega$, and similarly,
$x'=(x_1',\dots,x_d')=(x_1',y')$, $y'=(x'_2,\dots,x_d')$. We will abbreviate
$dx=dx_1dx_2\dots dx_d$ and $dy=dx_2\dots dx_d$, and
frequently use the fact that the space $L^2(\bbR\times\Omega;dx)=
L^2\big(\bbR;dx_1;L^2(\Omega;dy)\big)=
L^2\big(\Omega;dy;L^2(\bbR;dx_1)\big)$
is isometrically isomorphic to the space
$\ell^2\big(\bbZ^{d-1};L^2(\bbR;dx_1)\big)$ via the discrete Fourier
transform in the $y$-variables: 
\begin{equation}
w(x) = \sum_{j\in\bbZ^{d-1}}
\hat{w}_j(x_1)e^{ij\dott y}, \quad x=(x_1,y)\in\bbR\times\Omega,   \lb{4.3}
\end{equation}
where
\begin{equation}
\hat w_j(x_1) = (2\pi)^{1-d} \int_{\Omega} dy \, w(x_1,y) e ^{-i j \dott y}, \quad 
x_1 \in \bbR.   \lb{4.4}
\end{equation}

\subsubsection{Galerkin-based Evans function}\label{galerkin}
We first review the Galerkin approach described in \cite{LPSS00},
in which a standard Evans function is defined for a one-dimensional
truncation of the linearized operator about a standing-wave solution 
in a series of remarks.

\begin{remark}\label{decoupled}
Under our standing assumption of periodic boundary conditions on $\partial\Omega$,
there exist planar steady-state
solutions $U=U(x_1)$,
where $U$ is the solution of the
corresponding
one-dimensional problem \eqref{3.0a} described in Section \ref{s3}.
Linearizing about $U$, cf.\ \eqref{3.0}--\eqref{3.0c},
denoting $V(x_1) = -f'(U(x_1))$, and taking the discrete Fourier transform in
directions $x_2, \dots, x_d$, one obtains a decoupled family of
one-dimensional eigenvalue problems
\begin{equation}
0=(L_j -\lambda)\psi= \bigg(\f{d^2}{d x_1^2} - |j|^2 -\lambda - V(x_1)\bigg)\psi,
\end{equation}
indexed by Fourier frequencies $j=(j_2, \dots, j_d)\in \bbZ^{d-1}$,
each of which possess a well-defined Evans function and stability
index.  At $j=0$ and $\lambda=0$, there is an eigenfunction
$U'(x_1)$ associated with translation invariance in
the $x_1$-direction of the underlying equations;
for other $j$, there is typically no eigenfunction at $\lambda=0$.
Asymptotic analysis as in \cite{AGJ90}, \cite{PW92} yields a trivial,
positive stability index for $|j|$ sufficiently large, so that
computations may be truncated at a finite value of $|j|$.
\end{remark}

\begin{remark}\label{real}
In the above example, the operators $L_j$ are real-valued
(i.e., map real-valued functions into real-valued ones),
hence a stability index makes sense.
For more general, non-selfadjoint operators, one may
expand in sines and cosines to obtain a family of real-valued
eigenvalue equations for which a stability index may again be defined.
This principle extends further to generalized Fourier expansions in
the case of general $\Omega$, requiring only real-valuedness (in the
above sense) of the original (multi-dimensional) operator $L$.
\end{remark}

\begin{remark}\label{coupled}
More generally,
consider a standing-wave solution $U=U(x)$ that is not
planar, but only converges as $x_1\to \pm \infty$ to a constant state
$U_\infty$. We assume
\begin{equation}
\big(x_1 \mapsto \|U(x_1,\dott)\|_{H^{3/2}(\Omega)}\big)
\in \big(L^1 \cap L^\infty\big)(\bbR; dx_1).
\end{equation}
This writing means that the map
$x_1\mapsto U(x_1,\dott)$ from $\bbR$ into
the fractional Sobolev space $H^{3/2}(\Omega)$ is both
an $L^1$- and $L^\infty$-function with respect to the variable $x_1$.
Linearizing about $U$, denoting by $\psi=\psi(x)$
the corresponding eigenfunction, and taking the Fourier transform in
directions $x_2, \dots, x_d$, one obtains a coupled family of
one-dimensional eigenvalue problems
\begin{equation}\label{kevalue}
0= \bigg(\f{d^2}{d x_1^2} - |j|^2-\lambda\bigg)\hat\psi_j -
\big(\hatt V(x_1, \dott)* \hat\psi_{(\dott)}\big)(j), \quad j\in\bbZ^{d-1}. 
\end{equation}
Here $*$ denotes convolution in $j$,
\begin{equation}
\psi(x)=\sum_{j\in\bbZ^{d-1}}\hat\psi_j(x_1)e^{ij\dott y}, \quad x=(x_1,y)
\in\bbR\times\Omega, 
\end{equation}
and $\hatt V(x_1,j)$ denotes the value of the Fourier transform
of $V(x_1,\,\dott)$ in the variable $y=(x_2, \dots, x_d)$.
Following the approach of \cite{LPSS00}, one may proceed by
{\it Galerkin approximation}, truncating
the system at some sufficiently
high-order mode $|j|\le J$, to obtain again a very large,
but finite, real-valued eigenvalue ODE in $x_1$, for which one may
define in the usual way an Evans function and a stability index.
\end{remark}

\begin{remark}\label{BCs}
While we focus primarily on periodic boundary conditions on $\partial\Omega$ throughout this section, one can treat other boundary conditions such as Dirichlet, Neumann, or more generally, Robin-type boundary conditions in an analogous fashion. 
The key fact used in \eqref{4.3} and \eqref{4.4} is the eigenfunction expansion associated with the discrete eigenvalue problem of the self-adjoint Laplacian in 
$L^2(\Omega;dy))$ with periodic boundary conditions on $\partial\Omega$. The latter can be replaced by analogous discrete eigenvalue problems of the Laplacian with other self-adjoint boundary conditions on $\partial\Omega$.
\end{remark}

\subsubsection{Fredholm determinant version}\label{fred}
We now describe an alternative method based on the Fredholm
determinant, in which the Jost and Evans functions are prescribed canonically
as characteristic determinants, but computed by Galerkin
approximation: that is, we approximate the determinant rather
than the system of equations.

Specifically, consider again the general situation of Remark \ref{coupled}
of a solution $U$ of \eqref{multeq} decaying as $x_1\to \pm \infty$
to some constant state $U_\infty$. As in \eqref{3.0c},
we denote $V(x)=-f'(U(x))$, $x\in\Omega\times\bbR$, and
assume that
\begin{equation}\label{posass}
V_\infty =-f'(U_\infty)>0.
\end{equation}
We define the Birman--Schwinger operator $K(z)$ similarly to
\eqref{3.10}, \eqref{H0} as
\begin{equation}\label{K}
K(z)= -u(H^{(0)}_{\Omega, \rm p} -zI_{L^2(\bbR\times\Omega; dx)})^{-1}v, 
\quad z\in\bbC\backslash[V_\infty,\infty),
\end{equation}
with $H^{(0)}_{\Omega, \rm p}$ the self-adjoint realization of the differential expression 
$-\Delta_x +V_\infty$ in $L^2(\bbR\times\Omega; dx)$ with periodic boundary conditions on $\partial\Omega$, and
\begin{equation}
u(x) = \sgn(V(x)-V_\infty) v(x), \quad v(x) = |V(x)-V_\infty|^{1/2}
\end{equation}
for a.e.\ $x \in \bbR\times\Omega$.

For later reference, we define also the asymmetric
rearrangement of $K(z)$ by the formula
\begin{align}\label{tK} 
\cK(z) &= - \big(H^{(0)}_{\Omega, \rm p}-zI_{L^2(\bbR\times\Omega; dx)}\big)^{-1}uv  
\no  \\
& = -\big(H^{(0)}_{\Omega, \rm p} -zI_{L^2(\bbR\times\Omega; dx)}\big)^{-1}(V-V_\infty),  
\quad z\in\bbC\backslash[V_\infty,\infty).
\end{align}

Then, under the assumption
\begin{equation}\label{L1}
\big(x_1 \mapsto \|V(x_1,\dott) - V_\infty\|_{L^\infty(\Omega; dy)}\big)
\in \big(L^1\cap L^\infty\big)(\bbR;dx_1),
\end{equation}
we have the following result generalizing
the one-dimensional case \cite[Lemma 2.9]{GLM07}.
Fix $z\in\bbC\backslash[V_\infty,\infty)$. Passing to
adjoint operators, if needed, with no loss of generality
we will assume below that $\Im (z) \ge0$ and
fix the branch of the square root
such that $\Re\big((V_\infty + |j|^2-z)^{1/2}\big)>0$
for each $j\in\bbZ^{d-1}$, a choice consistent with the choice of $Q$ in \eqref{EQ}.

\begin{lemma}  \lb{KHS}
Assume \eqref{L1} and $d\le 3$. Then $K(z), \, \cK(z)\in
\cB_2\big(L^2(\bbR\times\Omega;dx)\big)$  for each $z\in\bbC\backslash
[V_\infty,\infty)$.
Moreover, the condition on the dimensions is sharp.
\end{lemma}
\begin{proof}
The operator $H^{(0)}_{\Omega, \rm p}$, since constant-coefficient,
decouples under the Fourier transform in the variables $x_2, \dots, x_d$. Consequently, 
the integral kernel of the resolvent of $H^{(0)}_{\Omega, \rm p}$, denoted by 
$\big(H^{(0)}_{\Omega, \rm p} -zI_{L^2(\bbR\times\Omega; dx)}\big)^{-1}(x,x')$, may be found explicitly by a Fourier expansion
and, using \eqref{3.12}, can be expressed as a countable sum
of scalar integral kernels:
\begin{align}\label{semisep}
& \big(H^{(0)}_{\Omega, \rm p} -zI_{L^2(\bbR\times\Omega; dx)}\big)^{-1}(x,x')   \no \\
& \quad =\frac{i}{2}
\sum_{j\in\bbZ^{d-1}}
(z - V_\infty - |j|^2)^{-1/2}
e^{i(z - V_\infty - |j|^2)^{1/2}|x_1-x_1'|} e^{ij\dott (y-y')},  \\
& \hspace*{4.8cm}  x=(x_1,y), x'=(x_1',y') \in \bbR\times \Omega,   \no 
\end{align}
where
$y=(x_2, \dots, x_d)\in\Omega$, $y'=(x'_2, \dots, x'_d)\in\Omega$,
and $j\in \bbZ^{d-1}$ denote the Fourier wave numbers
in these directions.
Using Parseval's identity, we obtain for any fixed
$x'\in\bbR\times\Omega$ that
\begin{align}
& \big\|\big(H^{(0)}_{\Omega, \rm p} -zI_{L^2(\bbR\times\Omega; dx)}\big)^{-1}(\dott,x')
\big\|_{L^2(\bbR\times\Omega;dx)}^2  \no \\
& \quad = \f{1}{4}
\sum_{j\in \bbZ^{d-1}}\int_\bbR dx_1 \, 
\Big|(z - V_\infty - |j|^2)^{-1/2}
e^{i(z - V_\infty - |j|^2)^{1/2}|x_1-x_1'|}\Big|^2    \nonumber\\
& \quad =\f{1}{4}\sum_{j\in \bbZ^{d-1}}\big|z - V_\infty - |j|^2\big|^{-1}
\int_\bbR dx_1 \, e^{-2\Im((z - V_\infty - |j|^2)^{1/2})|x_1-x_1'|}   \nonumber\\
& \quad =\f{1}{4}\sum_{j\in \bbZ^{d-1}}\big|z - V_\infty - |j|^2\big|^{-3/2}
\Big(\sin\tfrac12\big(\arg(z - V_\infty - |j|^2)\big)\Big)^{-1},  \label{xprimeest}
\end{align}
where
$\Im\big((z - V_\infty - |j|^2)^{1/2}\big)>0$
for $z\in\bbC\backslash[V_\infty,\infty)$ due to $\Im (z) \ge0$. Since
$\arg\big(z - V_\infty - |j|^2\big)\to\pi$ as $|j|\to\infty$, there is
a constant $c=c(z)$ such that
\begin{equation}\label{4.13old}
\big\|\big(H^{(0)}_{\Omega, \rm p} - zI_{L^2(\bbR\times\Omega; dx)}\big)^{-1}(\dott,x')
\big\|_{L^2(\bbR\times\Omega;dx)}^2\le c
\sum_{j\in \bbZ^{d-1}} (V_\infty + |j|^2)^{-3/2},
\end{equation}
hence is finite if and only if $d\le 3$. We recall the formula for
the Hilbert--Schmidt norm of the Hilbert--Schmidt operator $K$
with integral kernel $K(x,x')$ (see, e.g., \cite[Thm.\ 2.11]{Si05}, \cite[Sect.\ 1.6.5]{Ya92}):
\begin{equation}\label{yafform}
\|K\|_{\cB_2(L^2(\bbR\times\Omega;dx))}=
\|K(\dott,\dott)\|_{L^2((\bbR\times\Omega)\times
(\bbR\times\Omega); dx dx')}.
\end{equation}
Using \eqref{yafform} and \eqref{4.13old} to estimate the integral kernels
of \eqref{K} and \eqref{tK}, one infers
\begin{align}
& \|K(z,\dott,\dott)\|_{L^2((\bbR\times\Omega)\times
(\bbR\times\Omega); dx dx')}^2 \no \\
& \quad \le \Big\|u(\dott) 
\big\|\big(H^{(0)}_{\Omega, \rm p} -zI_{L^2(\bbR\times\Omega; dx)}\big)^{-1}
(\dott,\dott)\big\|_{L^2(\bbR\times\Omega;dx')}
\Big\|_{L^2(\bbR\times\Omega;dx)}^2
\|v\|_{L^\infty(\bbR\times\Omega;dx')}^2\nonumber\\
& \quad \le
c\sum_{j\in \bbZ^{d-1}}
(V_\infty+|j|^2)^{-3/2}\,\dott\,
\|u\|_{L^2(\bbR\times\Omega;dx)}^2
\|v\|_{L^\infty(\bbR\times\Omega;dx')}^2,\label{Kest}\\
& \|\cK(z,\dott,\dott)\|_{L^2((\bbR\times\Omega)\times
(\bbR\times\Omega); dx dx')}^2   \no \\
&\quad \le \Big\| 
\big\|\big(H^{(0)}_{\Omega, \rm p} -zI_{L^2(\bbR\times\Omega; dx)}\big)^{-1}
(\dott,\dott)\big\|_{L^2(\bbR\times\Omega;dx)}
u(\dott)v(\dott)\Big\|_{L^2(\bbR\times\Omega;dx')}^2\nonumber\\
& \quad \le
c\sum_{j\in \bbZ^{d-1}} (V_\infty+|j|^2)^{-3/2}
\,\dott\,
\|uv\|_{L^2(\bbR\times\Omega;dx')}^2,\label{tKest}
\end{align}
and finds that $K(z)$ and $\cK(z)$ are Hilbert--Schmidt
operators for $d\le 3$, as claimed.
In the decoupled case, where $u=u(x_1)$, $v=v(x_1)$, these estimates are sharp,
showing that in general $K(z)$, $\cK(z)$ are Hilbert--Schmidt only for $d\le 3$.
\end{proof}

\begin{definition}\label{multievans}
Assume \eqref{L1}. 
Generalizing the one-dimensional case \eqref{2.54a}, we introduce in
dimensions $d=2,3$, a {\em $2$-modified} Jost function defined by 
\begin{equation}\label{multidef}
\cF_2(z)= {\det}_{2,L^2(\bbR\times\Omega; d x)}
(I_{L^2(\bbR\times\Omega; d x)}- \cK(z)), \quad z\in\bbC\backslash[V_\infty,\infty). 
\end{equation}
\end{definition}

By the determinant property \eqref{2.34aa} one then obtains
\begin{equation}
\cF_2(z)= {\det}_{2,L^2(\bbR\times\Omega; d x)}
(I_{L^2(\bbR\times\Omega; d x)}- K(z)), \quad z\in\bbC\backslash[V_\infty,\infty), 
\end{equation}
which could equivalently have been used to define $\cF_2(z)$.

\begin{remark}
Calculations similar to Lemma \ref{KHS} show that $K(z)$ belongs to
successively weaker trace ideal classes as $d$ increases, for
which a higher-modified Jost function may be defined as a higher-modified
Fredholm determinant.
We restrict our attention here to the main physical cases $d=2,3$.
\end{remark}

\begin{remark}\label{factor}
Comparing with formula \eqref{2.54a}  for the {\em non-modified}
Jost function given in
the one-dimensional $d=1$ case (when $K(z)$
and $\cK(z)$ are trace-class operators), \eqref{2.34} implies the relation
\begin{align}\label{4.21old}
\begin{split}
& {\det}_{2,L^2(\bbR; dx_1)} (I_{L^2(\bbR; dx_1)}- K(z))  \\
& \quad ={\det}_{L^2(\bbR; dx_1)} (I_{L^2(\bbR; dx_1)}- K(z))
e^{\tr_{L^2(\bbR; dx_1)} (K(z))}.
\end{split}
\end{align}
Thus, \eqref{multidef}
differs from \eqref{2.54a} by a nonvanishing analytic factor
$e^{\tr_{L^2(\bbR; dx_1)} (K(z))}$,
and hence for practical purposes the
use of ${\det}_{2,L^2(\bbR; dx_1)}(\dott)$
and ${\det}_{L^2(\bbR; dx_1)}(\dott)$ in
\eqref{multidef} for $d=1$ are equivalent.
For $d>1$ we make a related
comment in Remark \ref{remtrKJ}.
\end{remark}

\subsubsection{Galerkin approximations}\label{approximations} 
Next, we approximate $\cF_2$ by a Galerkin approximation, working
for convenience with the asymmetric version \eqref{multidef}. 
We will augment \eqref{L1} with the more restrictive, but
still typically satisfied, condition
\begin{equation}\label{strongL1}
\big(x_1 \mapsto \|V(x_1,\dott)-V_\infty\|_{H^{3/2}(\Omega)}\big)\in L^2 (\bbR;dx_1).
\end{equation}
As in \eqref{defW}, we introduce 
\begin{equation}
W(x)=V(x)-V_\infty, \quad x\in\bbR\times\Omega, 
\end{equation}
and expand $W$ into a Fourier series in variables $y=(x_2,\dots,x_d)$
so that 
\begin{equation}
W(x)=\sum_{m\in\bbZ^{d-1}}\hatt{W}_m(x_1)e^{im\dott y}, 
\quad x=(x_1,y)\in\bbR\times\Omega. 
\end{equation}
Substituting \eqref{semisep}
into \eqref{tK}, we obtain an expansion of
the integral kernel of the operator $\cK(z)$:
\begin{align}
\cK(z,x,x') &=- 
\sum_{j\in\bbZ^{d-1}}
e^{ij\dott y}\,\,\frac{e^{-(V_\infty +
 |j|^2-z)^{1/2}|x_1-x_1'|}}{2(V_\infty + |j|^2-z)^{1/2}}\,\,
(V(x')-V_\infty)e^{-ij\dott y'}\label{1stexp}\\
&=- \sum_{j,m\in\bbZ^{d-1}}
e^{ij\dott y}\,\,\frac{e^{-(V_\infty +
 |j|^2-z)^{1/2}|x_1-x_1'|}}{2(V_\infty + |j|^2-z)^{1/2}}
\, \hatt W_{j-m}(x'_1)e^{-im\dott y'}.\label{2ndexp}
\end{align}
Introducing 
\begin{align}
\begin{split}
f^j_1(x)&=-2^{-1}e^{ij\dott y}
(V_\infty + |j|^2-z)^{-1/2}
e^{(V_\infty + |j|^2-z)^{1/2}x_1},\\
f^j_2(x)&=-2^{-1}e^{ij\dott y}
(V_\infty + |j|^2-z)^{-1/2}
e^{-(V_\infty + |j|^2-z)^{1/2}x_1},\\
g^j_1(x')&=e^{-(V_\infty + |j|^2-z)^{1/2}x_1'}
(V(x'_1)-V_\infty)e^{-ij\dott y'},\\
g^j_2(x')&=e^{(V_\infty + |j|^2-z)^{1/2}x_1'}
(V(x'_1)-V_\infty)e^{-ij\dott y'},   \label{fg}
\end{split}
\end{align}
we obtain from \eqref{1stexp} an expansion of $\cK(z,x,x')$ as a countable sum
\begin{equation}\label{Ksemisep}
\cK(z,x,x') =\begin{cases}
\sum_{j\in\bbZ^{d-1}}
f^j_1(x)  g^j_1(x'), & x_1'>x_1, \\
\sum_{j\in\bbZ^{d-1}}
f^j_2(x)  g^j_2(x'), & x_1>x_1', \\
\end{cases}
\end{equation}
of scalar integral kernels that are semi-separable in $x_1$.

Truncating \eqref{2ndexp} at some finite wave number $J$ or,
equivalently, Fourier expanding $f^j_k$, $g^j_k$ in
\eqref{Ksemisep} in variables $y$ and $y'$
and truncating the resulting series at some finite wave number $J$,
we obtain a sequence of Galerkin approximations
\begin{align}
\cK_J(z,x,x')&= - \sum_{|m|, |j|\le J}e^{ij\dott y}
\,\,\frac{e^{-(V_\infty +
 |j|^2-z)^{1/2}|x_1-x_1'|}}{2(V_\infty + |j|^2-z)^{1/2}}\,\,
\hatt W_{j-m}(x'_1)e^{-im\dott y'}\label{1FTKJ}\\
& =\begin{cases}
\sum\limits_{|m|, |j|\le J} e^{ij\dott y} \big(\hatt {f^j_1}\big)_j (x_1)
 \big(\hatt {g^j_1}\big)_{-m} (x_1')
e^{-im\dott y'}, & x_1'>x_1, \\
\sum\limits_{|m|, |j|\le J} e^{ij\dott y} \big(\hatt {f^j_2}\big)_j (x_1)
 \big(\hatt {g^j_2}\big)_{-m} (x_1')
e^{-im\dott y'}, & x_1>x_1', \\
\end{cases}\label{FTKJ}
\end{align}
where $\big(\hatt {f^j_k}\big)_m$ denote the Fourier coefficients of 
\begin{equation}
f^j_k(x_1,y)=\sum_{m\in\bbZ^{d-1}}
\big(\hatt {f^j_k}\big)_m e^{im\dott y}
=\big(\hatt {f^j_k}\big)_j e^{ij\dott y}, 
\end{equation}
and
\begin{equation}
\big(\hatt {g_k^j}\big)_m (x'_1)=\hatt W_{j+m} (x'_1)
e^{(-1)^{k+1}(V_\infty + |j|^2-z)^{1/2}x_1'}, \quad k=1,2,
\end{equation} 
denote the Fourier coefficients of the function
\begin{equation}
g_{k}^j(x')=
e^{(-1)^{k+1}(V_\infty + |j|^2-z)^{1/2}x_1'}W(x')e^{-ij\dott y'}, \quad k=1,2.
\end{equation}
We denote by $\cK_J(z)$ the integral operator
on $L^2(\bbR\times\Omega;dx)$ with the integral kernel
\eqref{1FTKJ}, \eqref{FTKJ}.

\begin{theorem}\label{convergence} 
Let $z\in\bbC\backslash[V_{\infty},\infty)$. Then 
under assumptions \eqref{L1} and  \eqref{strongL1}, $\cK_J(z)\in
\cB_2\big(L^2(\bbR\times\Omega; dx)\big)$
for dimensions $d=2,3$.  Moreover,
for $d=2,3$, $\cK_J(z)$ converges in the
Hilbert--Schmidt norm to $\cK(z)$ at rate $J^{(d-4)/4}$
as $J\to\infty$ and hence the sequence $\cF_{2,J}(z)$ defined by 
\begin{equation}\label{FJ}
\cF_{2,J}(z)= {\det}_{2,L^2(\bbR\times\Omega; dx)}
(I_{L^2(\bbR\times\Omega; dx)}- \cK_J(z)), 
\end{equation}
converges to $\cF_2(z)= {\det}_{2,L^2(\bbR\times\Omega; dx)}
(I_{L^2(\bbR\times\Omega; dx)}- \cK(z))$ as
$J\to\infty$ at rate $J^{(d-4)/4}$.
\end{theorem}
\begin{proof}
Using \eqref{2ndexp} and \eqref{1FTKJ}, one obtains 
\begin{align}\label{difference}
& \cK(z,x,x') - \cK_J(z,x,x')=-\bigg(\sum_{|j|=0}^J\sum_{|m|=J+1}^\infty
+\sum_{|j|=J+1}^\infty\sum_{|m|=0}^\infty\bigg)\\
& \quad \times \Big(2^{-1}e^{ij\dott y}(V_\infty + |j|^2-z)^{-1/2}
e^{-(V_\infty + |j|^2-z)^{1/2}|x_1-x_1'|}
\hatt W_{j-m}(x'_1)e^{-im\dott y'}\Big).     \nonumber
\end{align}
By the triangle inequality and Parseval's identity one
therefore infers
\begin{align} 
& \frac12 \|\cK(z) - \cK_J(z)\|^2_{L^2((\bbR\times\Omega)
\times(\bbR\times\Omega);dx\times dx')}  \no \\
& \quad \le \bigg\|\sum_{|j|=0}^J\sum_{|m|=J+1}^\infty(\dott)
\bigg\|^2_{L^2((\bbR\times\Omega)
\times(\bbR\times\Omega);dx dx')}  \no \\
& \qquad +\bigg\|\sum_{|j|=J+1}^\infty
\sum_{|m|=0}^\infty(\dott)
\bigg\|^2_{L^2((\bbR\times\Omega)
\times(\bbR\times\Omega);dx dx')}\nonumber\\
& \quad =\sum_{|j|=0}^J\sum_{|m|=J+1}^\infty
\bigg\|\frac{e^{-(V_\infty +
 |j|^2-z)^{1/2}|x_1-x_1'|}}{2(V_\infty + |j|^2-z)^{1/2}}
\hatt W_{j-m}(x'_1) \bigg\|^2_{L^2(\bbR\times\bbR;
dx_1 dx'_1)}\label{firstsum}\\
&\qquad +\sum_{|j|=J+1}^\infty\sum_{|m|=0}^\infty
\bigg\|\frac{e^{-(V_\infty +
 |j|^2-z)^{1/2}|x_1-x_1'|}}{2(V_\infty + |j|^2-z)^{1/2}}
\hatt W_{j-m}(x'_1) \bigg\|^2_{L^2(\bbR\times\bbR;
dx_1 dx'_1)}\,.\label{secondsum}
\end{align}
We will now estimate \eqref{firstsum} and \eqref{secondsum}
separately. Using arguments similar to \eqref{xprimeest}--\eqref{4.13old},
one observes that the sum in \eqref{firstsum} can be estimated as follows:
\begin{align}
\text{\eqref{firstsum}}&=\sum_{|j|=0}^J
\bigg\|\frac{e^{-(V_\infty +
|j|^2-z)^{1/2}|\dott|}}{2(V_\infty + |j|^2-z)^{1/2}}
\bigg\|^2_{L^2(\bbR;dx_1)}
\sum_{|m|=J+1}^\infty
\big\|\hatt W_{j-m}(\dott)\big\|^2_{L^2(\bbR;dx'_1)}\nonumber\\
&\le  c
\sum_{j\in \bbZ^{d-1}} (V_\infty + |j|^2)^{-3/2}\,
\sum_{|m|=J+1}^\infty
\big\|\hatt W_{m}(\dott)\big\|^2_{L^2(\bbR;dx'_1)}\nonumber\\
&=c'\|W-W_J\|^2_{L^2(\bbR\times\Omega;dx')}.\label{fs2}
 \end{align}
In the last equality we used that the series
$\sum_{j\in\bbZ^{d-1}}|j|^{-3}$ converges due to $1 \le d \le 3$,
and Parseval's identity for $W-W_J$, where
\begin{equation}
W_J(x)=\sum_{|m|\le J}\hatt W_m (x_1)e^{im\dott y}, \quad 
x=(x_1,y)\in\bbR\times\Omega, 
\end{equation}
is the truncation of $W$. By the Sobolev embedding
$W^{3/2,2}(\Omega)\hookrightarrow W^{3/4,4}(\Omega)$
(cf., e.g., \cite[Theorem 1.6.1]{He81}, \cite[p.\ 328, Eq.\ (8)]{Tr95}, 
\cite[Sects.\ I.4--I.6]{Wl87}) and
a standard Cauchy--Schwartz argument, one infers for each $x_1\in\bbR$, 
\begin{align}\label{trankW}
& \|W(x_1,\dott) - W_J(x_1,\dott)\|^2_{L^2(\Omega;dy)}=
\sum_{|j|>J} \big|\hatt W_j(x_1)\big|^2=
\sum_{|j|>J} |j|^{3/2} \big|\hatt W_j(x_1)\big|^2|j|^{-3/2} \no \\
& \quad \le \bigg(\sum_{|j|>J} |j|^{3}
\big|\hatt W_j(x_1)\big|^4\bigg)^{1/2}
\bigg(\sum_{|j|>J} |j|^{-3}\bigg)^{1/2}  \no \\
& \quad \le\|W(x_1,\dott)\|^2_{W^{3/4, 4}(\Omega)}
\bigg(\sum_{|j|>J} |j|^{-3}\bigg)^{1/2} \no \\
& \quad \le c\|W(x_1,\dott)\|^2_{H^{3/2}(\Omega)} J^{(d-4)/2},
\end{align}
Here we used standard notation $W^{s,2}(\dott)=H^s(\dott)$ for Sobolev spaces. 
Thus, by \eqref{fs2},
\begin{equation}
\text{ \eqref{firstsum} }
\le C_1 \|W\|^2_{L^2\big(\bbR;dx_1; H^{3/2}(\Omega))} J^{(d-4)/2}.  \label{trunc}
\end{equation}
Likewise, similarly to \eqref{xprimeest}--\eqref{4.13old},
one observes that the sum in \eqref{secondsum} can be estimated as follows:
\begin{align}
\text{\eqref{secondsum}}&=\sum_{|j|=J+1}^\infty
\bigg\|\frac{e^{-(V_\infty +
|j|^2-z)^{1/2}|\dott|}}{2(V_\infty + |j|^2-z)^{1/2}}
\bigg\|^2_{L^2(\bbR;dx_1)}
\sum_{|m|=0}^\infty
\big\|\hatt W_{j-m}(\dott)\big\|^2_{L^2(\bbR;dx'_1)}\nonumber\\
&\le  c
\sum_{|j|=J+1}^\infty (V_\infty + |j|^2)^{-3/2}\,
\sum_{|m|=0}^\infty
\big\|\hatt W_{m}(\dott)\big\|^2_{L^2(\bbR;dx'_1)}\nonumber\\
&\le c'\sum_{|j|=J+1}^\infty
|j|^{-3}\, \|W\|^2_{L^2(\bbR\times\Omega;dx')}
\le C_2J^{d-4}\|W\|^2_{L^2(\bbR\times\Omega;dx')}.\label{fs3}
 \end{align}
Combining \eqref{firstsum}, \eqref{secondsum}, \eqref{trunc},
\eqref{fs3}, and using \cite[Section 1.6.5]{Ya92}
as in \eqref{yafform}, one arrives at the estimate
\begin{equation}
\|\cK-\cK_J\|^2_{L^2((\bbR\times\Omega)
\times(\bbR\times\Omega);dx dx')}
\le 2 \|W\|^2_{L^2(\bbR\times\Omega;dx')}
\big(C_1J^{(d-4)/2}+C_2J^{d-4}\big),
\end{equation}
yielding the claimed result.
\end{proof}

Next, we take a closer look at properties of the integral operator $\cK_J(z)$ in 
$L^2(\bbR\times\Omega;dx)$ with integral kernel given by \eqref{1FTKJ}, assuming at first that 
\begin{equation}
\hatt W_{j-m} \in L^2(\bbR;dx_1), \quad m, j \in\bbZ, \; |m|, |j|\le J.   \lb{4.45a}
\end{equation} 
Using the fact that $L^2(\bbR\times\Omega;dx)$ decomposes into 
\begin{equation}
L^2(\bbR\times\Omega;dx)=L^2(\bbR;dx_1) \otimes L^2(\Omega;dy),  \lb{4.46}
\end{equation}
we will exploit the natural tensor product structure of the individual terms 
$\cK_{m,j}(z)$ in 
\begin{equation}
\cK_J(z) = \sum_{|m|, |j|\le J} \cK_{m,j}(z),  \quad z\in\bbC\backslash [V_\infty,\infty),   \label{4.47}  
\end{equation}
where $\cK_{m,j}(z)$, $m, j \in \bbZ^{d-1}$, $|m|, |j|\le J$, are integral operators in 
$L^2(\bbR\times\Omega;dx)$ with integral kernels given by 
\begin{align}
\begin{split}
& \cK_{m,j}(z,x,x') = - e^{ij\dott y} \,\,\frac{e^{-(V_\infty +
 |j|^2-z)^{1/2}|x_1-x_1'|}}{2(V_\infty + |j|^2-z)^{1/2}}\,\,
\hatt W_{j-m}(x'_1)e^{-im\dott y'}, \\
& \hspace*{3.7cm}  z\in\bbC\backslash [V_\infty,\infty), 
\; m, j \in \bbZ^{d-1}, \; |m|, |j|\le J.  \lb{4.48}
\end{split}
\end{align}
With respect to the tensor product structure \eqref{4.46}, the operators $\cK_{m,j}(z)$
decompose as
\begin{equation}
\cK_{m,j}(z) = \cA_{m,j}(z) \otimes \cB_{m,j},  \quad z\in\bbC\backslash [V_\infty,\infty), 
\lb{4.49}
\end{equation}
where the operator 
\begin{equation}
\cA_{m,j}(z)=\ol{\big(-(d^2/dx_1^2) + (V_\infty + |j|^2 -z)I_{L^2(\bbR;dx_1)}
\big)^{-1} \hatt W_{j-m}},  \quad z\in\bbC\backslash [V_\infty,\infty),   \lb{4.50}
\end{equation}
in $L^2(\bbR;dx_1)$ has the integral kernel 
\begin{equation}
\cA_{m,j}(z,x_1,x_1') = - \frac{e^{-(V_\infty +
 |j|^2-z)^{1/2}|x_1-x_1'|}}{2(V_\infty + |j|^2-z)^{1/2}}\,\,
\hatt W_{j-m}(x'_1),   \lb{4.51}
\end{equation}
and $\cB_{m,j}$ in $L^2(\Omega;dy)$ has the integral kernel 
\begin{equation}
\cB_{m,j}(y,y') = e^{ij\dott y} \, e^{-im\dott y'}.   \lb{4.52}
\end{equation}
In particular, each $\cB_{m,j}$ is a rank-one and hence trace class 
operator on $L^2(\Omega;dy)$,
\begin{equation}
\cB_{m,j}\in \cB_1\big(L^2(\Omega;dy)\big).    \lb{4.53}
\end{equation}
Next, $\cA_{m,j}(z)$ is a Hilbert--Schmidt operator on $L^2(\bbR;dx_1)$,  
\begin{equation}
\cA_{m,j}(z)\in \cB_2\big(L^2(\bbR;dx_1)\big), \quad z\in\bbC\backslash [V_\infty,\infty), 
\lb{4.54}
\end{equation}
if and only if (cf., \cite[Thm.\ 2.11]{Si05}, \cite[Sect.\ 1.6.5]{Ya92} and \eqref{yafform})
\begin{equation}
\hatt W_{j-m} \in L^2(\bbR;dx_1).   \lb{4.55}
\end{equation}
In addition, applying \cite[Theorem\ 4.5, Lemma\ 4.7]{Si05}, $\cA_{m,j}(z)$ is a trace class operator on $L^2(\bbR;dx_1)$,
\begin{equation}
\cA_{m,j}(z)\in \cB_1\big(L^2(\bbR;dx_1)\big), \quad z\in\bbC\backslash [V_\infty,\infty), 
\lb{4.56}
\end{equation}
if and only if  
\begin{equation}
\hatt W_{j-m} \in \ell^1\big(L^2(\bbR;dx_1)\big).   \lb{4.57}
\end{equation}
Here the Birman--Solomyak space $\ell^1\big(L^2(\bbR;dx_1)\big)$ is defined by
\begin{equation}
\ell^1\big(L^2(\bbR;dx_1)\big) = \bigg\{f\in L^2_{\loc}(\bbR)\,\bigg|\, 
\sum_{n\in\bbZ} \bigg(\int_{Q_n} dx_1 \, |f(x_1)|^2\bigg)^{1/2}<\infty\bigg\},  \lb{4.58}
\end{equation} 
with $Q_n$ the unit cube in $\bbR$ centered at $n\in\bbZ$. We recall that 
(cf.\ \cite[Ch.\ 4]{Si05})
\begin{align}
\begin{split}
& L^1(\bbR;(1+|x_1|)^\delta dx_1) \subsetneqq \ell^1\big(L^2(\bbR;dx_1)\big) 
\subsetneqq L^1(\bbR; dx_1) \cap L^2(\bbR; dx_1)   \lb{4.58a} \\
& \text{for all $\delta > 1/2$.}
\end{split}
\end{align}

We note in passing, that the symmetrized version $A_{m,j}(z)$ of $\cA_{m,j}(z)$, 
given by
\begin{align}
\begin{split}
& A_{m,j}(z) = \ol{\hat u_{j-m}
\big(-(d^2/dx_1^2) + (V_\infty + |j|^2 -z)I_{L^2(\bbR;dx_1)}\big)^{-1}\hat v_{j-m}},  \\
&\hspace*{7.65cm} z\in\bbC\backslash [V_\infty,\infty),   \lb{4.58b}
\end{split} 
\end{align}
where
\begin{equation}
\hat u_{j-m}(x_1) = \sgn\big(\hatt W_{j-m}(x_1)\big) \hat v_{j-m}(x_1), \quad 
\hat v_{j-m}(x_1) =  \big|\hatt W_{j-m}(x_1)\big|^{1/2}   \lb{4.58c}
\end{equation}
for a.e.\ $x_1\in\bbR$, is a trace class operator under the weaker assumption 
\begin{equation}
\hatt W_{j-m} \in L^1(\bbR;dx_1).   \lb{4.58d}
\end{equation}

Given these preparations, we can now summarize Hilbert--Schmidt and trace class properties of $\cK_J(z)$ as follows:

\begin{lemma} \lb{l4.10}
Assume $z\in\bbC\backslash [V_\infty,\infty)$. Then, \\
$(i)$ $\cK_J(z)$ is a Hilbert--Schmidt operator on $L^2(\bbR\times\Omega;dx)$ if 
$\hatt W_{j-m} \in L^2(\bbR;dx_1)$ for all $m, j \in \bbZ$, $|m|, |j|\le J$. \\
$(ii)$ $\cK_J(z)$ is a trace class operator on $L^2(\bbR\times\Omega;dx)$ if   
$\hatt W_{j-m} \in \ell^1\big(L^2(\bbR;dx_1)\big)$ for all $m, j \in \bbZ$, $|m|, |j|\le J$. 
\end{lemma}
\begin{proof}
Since the sum in \eqref{4.47} is finite, it suffices to prove the Hilbert--Schmidt and trace class properties of $\cK_{m,j}(z)$ for fixed $m, j$. Since by \eqref{4.49},
\begin{equation}
|\cK_{m,j}(z)| = |\cA_{m,j}(z)| \otimes |\cB_{m,j}|,  \lb{4.59}
\end{equation}
where as usual, $|T|=(T^*T)^{1/2}$, the singular values of $\cK_{m,j}(z)$ (i.e., the eigenvalues of $|\cK_{m,j}(z)|$) are square summable, respectively, summable, if and only if the singular values of $\cA_{m,j}(z)$ are square summable, respectively, summable, since $\cB_{m,j}$ is a rank-one operator and hence has precisely one nonzero singular value. This follows from the well-known fact that the spectrum of a tensor product $A_1\otimes A_2$ in $\cH_1 \otimes \cH_2$ of bounded operators $A_j$ in the complex separable Hilbert spaces $\cH_j$, $j=1,2$, is given by the product of the individual spectra
(cf., e.g., \cite[Theorem\ XIII.34]{RS78}), that is,
\begin{align}
\begin{split}
\sigma(A_1\otimes A_2)&=\sigma(A_1) \cdot \sigma(A_2) \\
&=\{\lambda\in\bbC \,|\, \lambda = \lambda_1\lambda_2, \, 
\lambda_j \in \sigma(A_j), \, j=1,2\}.   \lb{4.60}
\end{split}
\end{align}
Thus, one can apply \eqref{4.54}, \eqref{4.55}, respectively, \eqref{4.56}, \eqref{4.57}.
\end{proof}

Of course, the condition on $W$ in the Hilbert--Schmidt context in 
Lemma \ref{l4.10} is much weaker than condition \eqref{strongL1} since only finitely many Fourier coefficients $\hatt W_k$ of $W$ are in involved in the former, while the stronger condition \eqref{strongL1} is used to prove the convergence of $\cK_J$ in 
Theorem \ref{convergence}.

\begin{remark}\label{exp}
More generally, any useful approximation of $H^{(0)}_{\Omega, \rm p}$ may
be employed, not necessarily an eigenfunction expansion or
one attached to a Fourier basis.
In particular, in the case that $V_\infty$ is not constant
in $x_2, \dots, x_d$,
one may proceed alternatively by Galerkin approximation
to approximate $H^{(0)}_{\Omega, \rm p}$ as the limit of operators with 
semi-separable integral kernels 
corresponding to the (no longer decoupled)
truncated operator $H^{(0)}_{\Omega, \rm p, K}$, in a spirit similar to \cite{LPSS00}.
Likewise, it is not essential to assume that $V$ has common
limits at $x_1=+\infty$ and $x_1=-\infty$; one may consider
also ``front-type'' solutions with
$\lim_{x_1\to \pm \infty} U(x_1,\dott)=U_\pm$,
though this introduces some additional technicalities in the analysis
connected with invertibility of $H^{(0)}_{\Omega, \rm p}$.
\end{remark}

\subsubsection{Connection with Galerkin-based Evans functions}\label{connection}
At this point, adopting the point of view of \cite{LPSS00}, we
consider $\cK_J$ as an operator with a matrix-valued integral kernel,
and acting on the subspace $\cL_J$ of $L^2\big(\Omega;dy; L^2(\bbR;dx_1)\big)$
spanned by Fourier modes with wave-number of modulus less than or
equal to $J$, that is, on
\begin{equation}
\cL_J= \bigg\{w(x)=\sum_{|m|\le J}\hat{w}_m(x_1)e^{im\dott y} \,\bigg|\, 
\hat{w}_m\in L^2(\bbR;dx_1)\bigg\}. 
\end{equation}
We let $N_J$ denote
the number of these  modes.
Using Lemma \ref{l4.10}, one verifies that
$\cK_J(z)\in\cB_2\big(L^2(\bbR\times\Omega;dx)\big)$ and thus
\begin{equation}
{\det}_{2, L^2(\bbR\times\Omega;dx)} (I_{L^2(\bbR\times\Omega;dx)}-\cK_J(z))=
{\det}_{2, \cL_J} (I_{\cL_J}-\cK_J(z)) 
\end{equation}
is well-defined.
Equivalently, since the Fourier modes form an orthonormal basis, and
hence the determinant is invariant under the Fourier transform,
we compute (on the subspace $L^2\big(\bbR;dx_1;\bbC^{N_J}\big)$
of $L^2\big(\bbR;dx_1;\ell^2(\bbZ^{d-1})\big)$ isomorphic
to $\cL_J$ via the Fourier transform)
instead of \eqref{FJ} the Fredholm determinant
\begin{equation}
{\det}_{2, L^2(\bbR;dx_1;\bbC^{N_J})}
\big(I_{L^2(\bbR;dx_1;\bbC^{N_J})}-\hatt \cK_J(z)\big),
\end{equation}
where $\hatt \cK_J(z)$ on $L^2\big(\bbR;dx_1;\bbC^{N_J}\big)$ is defined 
in terms of its integral kernel  
\begin{align}\label{4.36old}
\hatt \cK_J(z,x_1,x'_1)&=\bigg(
\frac{-e^{-(V_\infty +
 |j|^2-z)^{1/2}|x_1-x_1'|}}{2(V_\infty + |j|^2-z)^{1/2}}
 \hatt W_{j-m}(x'_1)\bigg)_{|j|,|m|\le J}\\
&=\begin{cases}
F_1(x_1) G_1^\top (x_1'), & x_1'>x_1, \\
F_2(x_1) G_2^\top (x_1'), & x_1>x_1', \\
\end{cases}\label{vsemisep}
\end{align}
as an operator with a single, matrix-valued semi-separable integral kernel,
where $F_k$ and $G_k$ denote the $N_J\times N_J$ matrices
\begin{equation}\label{FG}
(F_k)_{m,j}=\big(\hatt {f_k^j}\big)_m, \quad
(G_k)_{m,j}=\big(\hatt {g_k^j}\big)_{-m}, \quad k=1,2,
\end{equation}
and $G_k^\top =((G_k)_{j,m})$ denotes the transpose of the matrix 
$G_k=((G_k)_{m,j})$, $k=1,2$.  

\smallskip
We briefly pause for a moment and recall the principal underlying idea here:
The operator $\cK_J$ acts on the space
$L^2\big(\bbR;dx_1; L^2(\Omega;dy)\big)$
and leaves invariant its subspace $\cL_J$ which,
in fact, is isomorphic to
$L^2\big(\bbR;dx_1; L^2_{N_J}(\Omega;dy)\big)$,
where $L^2_{N_J}(\Omega;dy)$ is the
subspace of $L^2(\Omega;dy)$
spanned by the $N_J$ exponentials
$\{e^{i y\dott j}\}_{|j|\le J}$.
Via the Fourier transform,
$L^2_{N_J}(\Omega;dy)$ is isometrically
isomorphic to $\bbC^{N_J}$ viewed as a subspace
in $\ell^2(\bbZ^{d-1})$. Indeed, if $j\in\bbZ^{d-1}$ and
$|j|\le J$, then the Fourier transform
maps the element $e^{i j\dott y}$ of the basis of
$L^2(\Omega;dy)$ into the sequence 
$\delta_j=\{\delta_{m,j}\}_{m\in\bbZ^{d-1}}\in\ell^2(\bbZ^{d-1})$. Indexing a basis in $\bbC^{N_J}$ by means of the indices $j\in\bbZ^{d-1}$,
$|j|\le J$, we fix an isomorphism between $\bbC^{N_j}$ and
the subspace of $\ell^2(\bbZ^{d-1})$
spanned by $\delta_j$, $|j|\le J$, and thus between
$\bbC^{N_j}$ and $L^2_{N_J}(\Omega;dy)$. Clearly,
via the Fourier transform, the operator
$\cK_J$ on $L^2\big(\bbR;dx_1; L^2_{N_J}(\Omega;dy)\big)$
then becomes the operator $\hatt \cK_J$ on
$L^2\big(\bbR;dx_1;\bbC^{N_J}\big)$.
In particular,
if $\cK_J$ (and hence $\hatt \cK_J$) is
also of trace class, then the traces of $\cK_J$
and $\hatt \cK_J$ are equal as will be used below
in the proof of the second equality in \eqref{trformKJone}.
It is the operator $\hatt \cK_J$ which finally is an operator
with a semi-separable integral kernel.

\begin{remark}\label{vector}
The vector-valued case $w \in \bbR^n$ may be treated similarly,
with $V_\infty$ now a positive-definite $n\times n$ matrix, and
$F_k, \, G_k \in \bbR^{nN_J \times nN_J}$.
\end{remark}

Our objective is to relate the truncated $2$-modified Jost function 
\eqref{FJ} and the Evans function for the eigenvalue problem
\eqref{kevalue}. We rewrite \eqref{kevalue} in matrix form as
\begin{equation}\label{matsys}
\bigg(- \f{d^2}{d x_1^2} +V_\infty I_{\bbC^{N_J}}
 +\diag\big\{|j|^2\big\}_{|j|\le J} - zI_{\bbC^{N_J}}
 + \mathcal{W}_J\bigg)\Psi=0,
\end{equation}
where $j$ is the Fourier wave number, $|j|\le J$,
$\Psi=\Psi(x_1)$ is an $\bbC^{N_J}$-valued function
on $\bbR$, $\diag\big\{|j|^2\big\}_{|j|\le J}$ is a diagonal matrix of dimensions 
$N_J \times N_J$, and $\mathcal{W}_J+V_\infty I_{\bbC^{N_J}}$ 
is the matrix representation of
some chosen truncation of the convolution operator
$\hatt V * \dott$.\ This choice is to be followed consistently
in both Galerkin-based and Fredholm computations. Specifically,
if $\Psi=(\hat\psi_m)_{|m|\le J}$ for the
eigenfunction $\psi$ in \eqref{kevalue},
then 
\begin{equation}
\mathcal{W}_J(x_1)=\big(\hatt W_{j-m}(x_1)\big)_{|j|,|m|\le J}, \quad x_1 \in\bbR, 
\end{equation}
where $\hatt W_\ell(x_1)$ are the Fourier coefficients of 
$W(x_1,\dott)=V(x_1,\dott)-V_\infty$. Introducing 
the self-adjoint operator $\cH_J$ in $L^2\big(\bbR;dx_1;\bbC^{N_J}\big)$, 
\begin{equation}
\cH_J = - \f{d^2}{d x_1^2} +V_\infty I_{L^2(\bbR; dx_1;\bbC^{N_J})} +
\diag\big\{|j|^2\big\}_{|j|\le J}+ \mathcal{W}_J, 
\quad \dom(\cH_J)=H^2\big(\bbR; \bbC^{N_J}\big),  
\lb{4.52a}
\end{equation}
we note that the asymptotic operator for the
operator $\cH_J$ in \eqref{4.52a},
viewed as a one-dimensional matrix-valued second-order
differential operator in $L^2\big(\bbR;dx_1; \bbC^{N_J}\big)$, is given by 
\begin{equation}
\cH_J^{(0)} = - \f{d^2}{d x_1^2} + V_\infty I_{L^2(\bbR; dx_1;\bbC^{N_J})} +
\diag\big\{|j|^2\big\}_{|j|\le J}, 
\quad \dom\big(\cH_J^{(0)}\big)=H^2\big(\bbR; \bbC^{N_J}\big),  
\lb{4.53a}
\end{equation}
and thus the operator
\begin{equation}\label{defcK}
\hatt \cK_J(z)=-\big(\cH^{(0)}_J -zI_{L^2(\bbR; dx_1;\bbC^{N_J})}\big)^{-1}\cW_J, 
\quad z\in \bbC\backslash\sigma\big(\cH^{(0)}_J\big), 
\end{equation}
is the Birman--Schwinger-type operator (cf.\ \eqref{tK}) 
for the pair of the truncated operators $\cH_J$ and $\cH^{(0)}_J$.

\begin{lemma} \lb{l4.13}
Let $z\in\bbC$ such that $z\notin \sigma\big(\cH_J^{(0)}\big)$. In addition, 
assume that $\hatt W_{j-m} \in \ell^1\big(L^2(\bbR;dx_1)\big) \cap C(\bbR)$
for all $m, j \in \bbZ$, $|m|, |j|\le J$. Then the operator
$\hatt \cK_J(z)$ is of trace class on $L^2(\bbR;
dx_1;\bbC^{N_J})$, the operator $\cK_J(z)$ is of trace class on
$L^2(\bbR\times\Omega; dx)$, their traces are equal and given by
the following expression $\Theta_J(z)$:
\begin{align}
\Theta_J(z)& =\tr_{L^2(\bbR\times\Omega; dx)}(\cK_J(z))
=\tr_{L^2(\bbR; dx_1;\bbC^{N_j})}\big(\hatt
\cK_J(z)\big)\label{trformKJone}\\
&=-\frac12\bigg(\int_{\bbR\times\Omega} dx\, W(x) \bigg)
\sum_{|j|\le J}(V_\infty + |j|^2-z)^{-1/2}.      \label{trformKJ}
\end{align}
\end{lemma}
\begin{proof}  
By Lemma \ref{l4.10}\,$(ii)$, $\cK_J(z)$ and each $\cK_{m,j}(z)$ in \eqref{4.47} is a  trace class operator on $L^2(\bbR\times\Omega; dx)$. In addition, each of the integral kernels of $\cK_J(z)$ and $\cK_{m,j}(z)$ is continuous on the diagonal. Thus, \cite[Corollary\ 3.2]{Br91} applies and hence 
\begin{align}
& \tr_{L^2(\bbR\times\Omega; dx)} (\cK_J(z))
= \sum_{|m|, |j| \le J} \tr_{L^2(\bbR\times\Omega; dx)} (\cK_{m,j}(z))  \no \\
& \quad = - \f{1}{2}\sum_{|m|, |j| \le J} \int_{\bbR} dx_1 \, \hatt W_{j-m}(x_1) 
\f{1}{(V_\infty + |j|^2 -z)^{1/2}} \int_{\Omega} dy \, e^{i(j-m)y}  \no \\
& \quad = - \f{1}{2}\sum_{|m|, |j| \le J} \int_{\bbR} dx_1 \, \hatt W_{j-m}(x_1) 
\f{1}{(V_\infty + |j|^2 -z)^{1/2}} \, (2 \pi)^{d-1} \delta_{m,j}  \no \\
& \quad = - \f{1}{2} \bigg((2 \pi)^{d-1} \int_{\bbR} dx_1 \, \hatt W_0 (x_1)\bigg) 
\sum_{|j| \le J} \f{1}{(V_\infty + |j|^2 -z)^{1/2}}  \no \\
& \quad = - \f{1}{2} \bigg(\int_{\bbR \times \Omega} dx \, W(x)\bigg) 
\sum_{|j| \le J} \f{1}{(V_\infty + |j|^2 -z)^{1/2}},    \lb{4.77}
\end{align}
proving \eqref{trformKJ}. 

Finally, denote by $\cP_J$ the orthogonal projection in
$L^2(\bbR\times\Omega;dx)$ onto $\cL_J$ and by
$\cQ_J = I_{L^2(\bbR\times\Omega; dx)} - \cP_J$ the
complementary projection.
Since $\cL_J$ is a reducing subspace for
$\cK_J$,
\begin{equation}
\cK_J \cP_J = \cP_J \cK_J,    \lb{4.78}
\end{equation}
one can write $\cK_J$ in $L^2(\bbR\times\Omega; dx)$ in terms of
the
$2\times 2$ block decomposition
\begin{equation}
\cK_J(z) = \begin{pmatrix} \cK_J\big|_{\ran(\cP_J)}(z) & 0
\\[2mm] 0 & 0 \end{pmatrix}
\lb{4.79}
\end{equation}
with respect to the decomposition
\begin{align}
\begin{split}
L^2(\bbR\times\Omega; dx) &= \cP_J L^2(\bbR\times\Omega; dx) \oplus
\cQ_J L^2(\bbR\times\Omega; dx) \\
& = \cL_J \oplus  \cQ_J L^2(\bbR\times\Omega; dx).   \lb{4.80}
\end{split}
\end{align}
Since $ \cK_J\big|_{\ran(\cP_J)}$ is unitarily equivalent to
$\hatt K_J(z)$ via the Fourier transform, \eqref{4.79} implies that $\cK_J(z)$ 
and $\hatt \cK_J(z)$ are trace class operators at the same time, and it also implies equality of the following traces:
\begin{align}
\begin{split}
\tr_{L^2(\bbR\times\Omega; dx)}(\cK_J(z))
& = \tr_{\cL_J}\Big(\cK_J\big|_{\ran(\cP_J)}(z)\Big)   \\
& =\tr_{L^2(\bbR; dx_1;\bbC^{N_j})}\big(\hatt \cK_J(z)\big),   
\label{4.81}
\end{split}
\end{align}
proving the second equality in \eqref{trformKJone}. 
\end{proof}
 
\begin{remark}\label{remtrKJ} 
We emphasize that the sequence $\{\Theta_J(z)\}_{J\ge0}$ diverges
as $J\to\infty$ for $d\ge2$. The latter fact does not permit us to
use in the subsequent analysis the non-modified Jost function
$\cF_J(z)={\det}_{L^2(\bbR\times\Omega; dx)}
(I_{L^2(\bbR\times\Omega; dx)}- \cK_J(z))$ and pass in \eqref{4.39old}
to the limit as $J\to\infty$; instead, it forces us to work with the 
$2$-modified Jost function $\cF_{2,J}(z)$ defined in \eqref{FJ}.
\end{remark}

Denote by $\cE_J$ the Evans functions for the one-dimensional
approximate system \eqref{kevalue} obtained by Galerkin
approximation/Fourier truncation at the level $|j|\le J$,
for simplicity of discussion normalized as described in \cite{GLM07}
to agree with the corresponding (one-dimensional) $2$-modified
Fredholm determinant.
Following the approach of \cite{GLM07}, we recall that the
Evans function $\cE_J(z)$ is defined as a
$2N_J\times 2N_J$ Wronskian
\begin{equation}\label{defEvfnct}
\cE_J(z) =\det (\cY)=\det(\cY_1^+, \dots, \cY_{2N_J}^-),
\end{equation}
where the $(2N_J\times 1)$-vectors $\cY_1^+, \dots, \cY_{2N_J}^-$
are bases of solutions decaying at $x_1=+\infty$, respectively, at
$x_1=-\infty$, of the first-order system equivalent to the second-order differential
equation \eqref{matsys}, lying in appropriately prescribed directions at spatial
infinity. The solutions are chosen in \cite{GLM07}
in a way that $\cE_J$ does not depend on the choice
of coordinate system in $\bbC^{2N_J}$ and does not change
under similarity transformations of the system. Specifically, let us introduce the $(2N_J\times 2N_J)$
matrices
\begin{equation}\label{defAR}
\cA =\begin{pmatrix}0&I_{\bbC^{N_J}}\\ h&0\end{pmatrix},\quad
R(x_1) =\begin{pmatrix}0&0\\\cW_J(x_1)&0\end{pmatrix}, \quad x_1\in\bbR,
\end{equation}
where, for brevity, we denote
\begin{equation}\label{defnh}
h(z) =(V_\infty-z)I_{\bbC^{N_J}}+ \diag\big\{|j|^2\big\}_{|j|\le J}.
\end{equation} Now $d \cY/d x_1=(\cA+R(x_1))\cY$ is the first-order
system equivalent to the second-order differential
equation \eqref{matsys}. We introduce on the space
$L^2\big(\bbR;dx_1;\bbC^{2N_J}\big)$ the first order differential operators
$G^{(0)}_J=(d/d x_1)-\cA$ and $G_J=(d/d x_1)-\cA-R(x_1)$ and the
respective Birman--Schwinger-type integral operator
\begin{equation}\label{deftk}
\widetilde\cK_J(z) =-((d/dx_1)-\cA)^{-1}R(x_1).
\end{equation}

According to the main result in \cite{GLM07},
the Evans function $\cE_J(z)$ is equal  (up to the 
explicitly computed factor $e^{-\Theta_J(z)}$, that is non degenerate and 
analytic with respect to $z$)
to the $2$-modified Fredholm determinant of the
operator $I_{L^2(\bbR;dx_1;\bbC^{2N_J})}-\widetilde\cK_J(z)$ on
$L^2\big(\bbR;dx_1;\bbC^{2N_J}\big)$ corresponding
to the first-order system mentioned above. As we will see next (cf.\ also 
\cite[Theorem 4.7]{GM04}),
this $2$-modified Fredholm determinant is equal to the
$2$-modified Fredholm determinant of the operator
$I_{L^2(\bbR;dx_1;\bbC^{N_J})}-\hatt \cK(z)$ on $L^2\big(\bbR;dx_1;\bbC^{N_J}\big)$ corresponding to the second-order
operator $\cH_J$ in \eqref{matsys}. Thus,
we obtain evidently that the Evans function
$\cE_J$ coincides with the following
(non-modified(!)) Fredholm
determinant (that is, with the Jost function):
\begin{equation}\label{EJ}
\cF_J(z) ={\det}_{L^2(\bbR;dx_1;\bbC^{N_J})}
\Big(I_{L^2(\bbR;dx_1;\bbC^{N_J})} 
+ \big(\cH^{(0)}_J - z I_{L^2(\bbR;dx_1;\bbC^{N_J})}\big)^{-1} \cW _J \Big). 
\end{equation} 

Thus, we have the following main result, extending the
one-dimensional theory of \cite{GLM07}. One of its main points
can be explained as follows: Zeros of the Evans function $\cE_J$
or, equivalently, of the Jost function $\cF_J$, are the
eigenvalues of the operator $\cH_J$. The eigenvalues are also 
zeros of the 2-modified Fredholm determinant
$\cF_{2,J}$. The modified and nonmodified determinants are equal
up to the nonzero exponential factor $e^{\Theta_J}$, where 
$\Theta_J$ is the trace described in \eqref{trformKJ}. From this 
point of view the use of the nonmodified
determinant $\cF_J$ (or $\cE_J$) and the 2-modified determinant
$\cF_{2,J}$ are equivalent, as far as finding the eigenvalues of $\cH_J$ 
is concerned. However, the nonmodified determinants
have an advantage because the sequence $\{\cF_{2,J}\}$
converges to $\cF_{2}$ as $J\to\infty$, while the sequences
 $\{\cF_J\}$, $\{\cE_J\}$, and $\{\Theta_J\}$, all diverge. Thus, for
the truncated problem, the use of the
 2-modified Fredholm determinants appears to be more natural
than the use of the Evans function.

\begin{theorem}\label{equivalence}
Let $z\in\bbC$ such that $z\notin \sigma\big(\cH_J^{(0)}\big)$ and assume that 
$\hatt W_{j-m} \in \ell^1\big(L^2(\bbR;dx_1)\big) \cap C(\bbR)$
for all $m, j \in \bbZ$, $|m|, |j|\le J$. Then the Galerkin-based Evans function $\cE_J(z)$, the Jost function
$\cF_J (z)={\det}_{L^2(\bbR\times\Omega; dx)}
(I_{L^2(\bbR\times\Omega; dx)}- \cK_J(z))$, and the
approximate modified Fredholm determinants for the operators
in \eqref{defcK} and \eqref{deftk} are related as follows:
\begin{align}
\cF_{2,J}(z) &={\det}_{2,L^2(\bbR\times\Omega; dx)}
(I_{L^2(\bbR\times\Omega; dx)}- \cK_J(z))\label{f2J}\\
&={\det}_{2,L^2(\bbR; dx_1;\bbC^{2N_J})}
\big(I_{L^2(\bbR; dx_1;\bbC^{2N_J})}- \widetilde\cK_J(z)\big)  \label{tildekk}\\
& = e^{\Theta_J(z)}\cE_J(z)\label{tkevans}\\
& =e^{\Theta_J(z)}{\det}_{L^2(\bbR\times\Omega; dx)}
(I_{L^2(\bbR\times\Omega; dx)}- \cK_J(z))\label{evnmdet}\\
&=e^{\Theta_J(z)}\cF_J(z)\label{notFJ},\end{align}
where $\Theta_J(z)$ is the trace of the operator $\cK_J(z)$ given
in formula \eqref{trformKJ}.
\end{theorem}
\begin{proof} Since $\cK_J(z)$ is of trace class on $L^2(\bbR\times\Omega; dx)$, 
$\Theta_J(z)$ is the trace of $\cK_J(z)$, and $\cF_J(z)$ is just a notation for 
${\det}_{L^2(\bbR\times\Omega; dx)}(I_{L^2(\bbR\times\Omega; dx)}- \cK_J(z))$, the equality of \eqref{f2J}, \eqref{evnmdet} and \eqref{notFJ} trivially follows from 
(cf.\ \eqref{2.34})
\begin{align}
\begin{split}
&{\det}_{2,L^2(\bbR\times\Omega; dx)} 
(I_{L^2(\bbR\times\Omega; dx)}- \cK_J(z))\\
& \quad ={\det}_{L^2(\bbR\times\Omega; dx)}
(I_{L^2(\bbR\times\Omega; dx)}- \cK_J(z))
e^{\tr_{L^2(\bbR\times\Omega;dx)} (\cK_J(z))}.  \label{4.39old}
\end{split}
\end{align}
To show that the modified Fredholm determinants in
\eqref{f2J} and \eqref{tildekk} are equal, we will utilize an idea from the
proof of \cite[Proposition 8.1]{LP07}
(see also a related result in \cite[Theorem 4.7]{GM04}).
We introduce the following operator matrices
acting on $L^2\big(\bbR;dx_1;\bbC^{2N_J}\big)$:
\begin{align}\label{deftt}
T&=2^{-1/2}\begin{pmatrix}h^{1/2}&-I
\\h^{1/2}&I\end{pmatrix},\quad
T^{-1}=2^{-1/2}
\begin{pmatrix}h^{-1/2}&h^{-1/2}\\-I&I\end{pmatrix},\\
E&=\begin{pmatrix}-I&-I\\I&I\end{pmatrix},\quad Q=
\begin{pmatrix}I&0\\0&0\end{pmatrix},   \label{EQ}
\end{align}
where $I=I_{L^2(\bbR;dx_1;\bbC^{N_J})}$
and $h$ is defined in \eqref{defnh}. In addition, we use the related similarity transformation in \eqref{defAR} and  \eqref{deftk} to define the
following matrices and operators
\begin{equation}\label{simtrans}
\cA^{(1)} =T\cA T^{-1}, \quad  R^{(1)}(x_1) =TR(x_1) T^{-1}, \quad 
\widetilde\cK_J^{(1)}(z) =T\widetilde\cK_J(z) T^{-1}.
\end{equation}
A short calculation reveals:
\begin{align}
\cA^{(1)}&=\diag\big\{-h^{1/2}, h^{1/2}\big\},\label{calcA}\\
R^{(1)}(x_1)&=2^{-1}\cW_JEh^{-1/2},\label{newR}\\
\widetilde\cK_J^{(1)}(z)&=2^{-1}\begin{pmatrix}
 -((d/d x_1)+h^{1/2})^{-1}\cW_J\\
((d/dx_1)-h^{1/2})^{-1}\cW_J \end{pmatrix}
\begin{pmatrix}
h^{-1/2}& h^{-1/2}\end{pmatrix}.\label{KJone}\end{align}
Changing the order of multiplication of
the operators blocks in \eqref{KJone}, one infers 
\begin{equation}
2^{-1}\begin{pmatrix}
h^{-1/2}& h^{-1/2}\end{pmatrix}\begin{pmatrix}
 -((d/d x_1)+h^{1/2})^{-1}\cW_J\\
((d/d x_1)-h^{1/2})^{-1}\cW_J \end{pmatrix}=
\big((d^2/d x_1^2) - h\big)^{-1}\cW_J.
\end{equation}
Thus, using the standard determinant property \eqref{2.2}
and recalling \eqref{defcK},
we conclude that \eqref{f2J} and \eqref{tildekk} are equal, 
\begin{align}
&{\det}_{2,L^2(\bbR; dx_1;\bbC^{2N_J})}
\big(I_{L^2(\bbR; dx_1;\bbC^{2N_J})}- \widetilde\cK_J(z)\big)  \no \\
& \quad ={\det}_{2,L^2(\bbR; dx_1;\bbC^{2N_J})}
\big(I_{L^2(\bbR; dx_1;\bbC^{2N_J})}- \widetilde\cK^{(1)}_J(z)\big)  \no \\
& \quad ={\det}_{2,L^2(\bbR; dx_1;\bbC^{N_J})}
\big(I_{L^2(\bbR; dx_1;\bbC^{N_J})}- ((d^2/d^2 x_1^2) - h)^{-1}\cW_J\big) \no \\
& \quad ={\det}_{2,L^2(\bbR;dx_1;\bbC^{N_J})}
\Big(I+\big(\cH_J^{(0)} - z I\big)^{-1} \cW _J \Big)  \no \\
& \quad ={\det}_{2,L^2(\bbR; dx_1;\bbC^{N_J})}
\big(I_{L^2(\bbR; dx_1;\bbC^{N_J})}- \hatt \cK_J(z)\big)  \no \\
& \quad ={\det}_{2,L^2(\bbR\times\Omega; dx)}
(I_{L^2(\bbR\times\Omega; dx)}- \cK_J(z)).
\end{align}

Finally, to show that \eqref{tildekk} and
\eqref{tkevans} are equal,
we apply the similarity transformation \eqref{simtrans}
and replace the differential equation $d \cY/d x_1=(\cA+R(x_1))\cY$
by $d \cY/dx_1 =(\cA^{(1)}+R^{(1)}(x_1))\cY$.
We will now use one of the main results of \cite{GLM07}.
Since the real part of the spectrum of $h^{1/2}$ (for 
$h$ defined in \eqref{defnh}) is positive due to our convention $\Im(z)\ge 0$,
the unperturbed
equation $d \cY/dx_1 =\cA^{(1)}\cY$,
due to \eqref{calcA}, has the exponential
dichotomy on $\bbR$ with the dichotomy projection $Q$ defined
in \eqref{EQ}.  Therefore, according to \cite[Theorem 8.37]{GLM07}, under
assumption \eqref{L1} we have the formula
\begin{equation}
{\det}_{2,L^2(\bbR; dx_1;\bbC^{2N_J})}
\big(I_{L^2(\bbR; dx_1;\bbC^{2N_J})}- \widetilde\cK^{(1)}_J(z)\big)
= e^{{\widetilde\Theta}_J(z)}\cE_J(z),  
\end{equation}
where $\widetilde\Theta_J(z)$ is defined as follows:
\begin{equation}
\widetilde\Theta_J(z)=\int_0^\infty dx_1\, \tr_{\bbC^{2N_J}}
\big[ QR^{(1)}(x_1)\big]- \int_{-\infty}^0 dx_1 \, 
\tr_{\bbC^{2N_J}} \big[(I_{\bbC^{2N_J}}- Q)R^{(1)}(x_1)\big].
\end{equation}
Using \eqref{EQ} and \eqref{newR}, it follows that
$\widetilde\Theta_J(z)=\Theta_J(z)$, completing the proof. 
\end{proof}

\begin{remark}\label{vector2}
We have a similar result in the vector-valued case mentioned
in Remark \ref{vector} also in the front-type setting discussed
in Remark \ref{exp}.
\end{remark}

\subsubsection{Alternative computation}\label{alt}

For its own interest, and for reference in the following subsections,
we mention an alternative method of computing $\cF_{2,J}$
directly from the reduction of \cite{GM04}, where the Jost
function has been computed, without carrying out
the full analysis of \cite{GLM07} relating this to the Evans function.
Comparing \eqref{vsemisep}, and \cite[(1.17)]{GM04}  with $\alpha=1$,
we have the following representation:
\begin{equation}\label{Jost}
\cF_{2,J}(z)=
{\det}_{2,\bbC^{N_J}}\bigg(I_{\bbC^{N_J}}-
\int_{\bbR^2} \,dx_1\,dx'_1\, G_2(x_1)\big(I_{\bbC^{2N_J}}+\cJ(x_1,x_1')\big 
)F_2(x_1')\bigg),
\end{equation}
where
\begin{align}\label{J}
\cJ(x_1,x_1')&=C(x_1)\cU(x_1)^{-1}\cU(x_1')B(x_1'),\\
B&=\begin{pmatrix} G^\top _1\, & \,-G^\top _2\end{pmatrix}^\top ,
\quad C=\begin{pmatrix} F_1\, & \,F_2\end{pmatrix},
\label{BC}\\
A&=
\begin{pmatrix}
G_1^\top  F_1& G_1^\top  F_2\\[1mm]
-G_2^\top  F_1& -G_2^\top  F_2\\
\end{pmatrix},\label{A}
\end{align}
and $\cU$ is any nonsingular solution of the first-order system 
\begin{equation}\label{Ueq}
\f{d \cU(x_1)}{dx_1}= A(x_1)\cU(x_1).
\end{equation}

The formulation \eqref{Ueq} is not numerically useful,
since the off-diagonal elements
of $A$ are exponentially growing with rate of order
$e^{J|x_1|}$.
However, noting that $G_{1,2}^\top (x_1')$, respectively,
$F_{1,2}(x_1)$ factor
as $\diag \big\{e^{\mp (V_\infty+|j|^2- z)^{1/2} x_1'}
\cW_J \big\}$, respectively,
$\diag \big\{(V_\infty+|j|^2- z)^{-1/2}
e^{\pm (V_\infty + |j|^2-z)^{1/2}x_1}\big\}$,
we may reduce \eqref{Ueq} by the coordinate change $\cV=D\cU$,
with
\begin{equation}\label{defDchV}
D=
\diag \bigg\{\frac{e^{(V_\infty+
|j|^2-z)^{1/2}x_1}}{(V_\infty+|j|^2-z)^{1/2}}, \,
\frac{e^{-(V_\infty+
|j|^2-z)^{1/2}x_1}}{(V_\infty+|j|^2-z)^{1/2}} \bigg\} ,
\end{equation}
to a bounded-coefficient system (cf.\ \eqref{defnh}),
\begin{align}\label{altsys}
\begin{split}
\f{d \cV(x_1)}{dx_1}&=A_b(x_1)\cV(x_1),\\
 A_b(x_1)&= \begin{pmatrix}
h^{1/2}+h^{-1/2}\cW_J(x_1)&
h^{-1/2}\cW_J(x_1)\\[1mm]
-h^{-1/2}\cW_J(x_1)&
-h^{1/2}-h^{-1/2}\cW_J(x_1)
\end{pmatrix},
\end{split}
\end{align}
of a form readily solved by the same techniques used
to solve the first-order eigenvalue ODE for
basis solutions $\cY_\ell$ in \eqref{defEvfnct}.
Indeed, this can be recognized as essentially the same ODE.

\begin{remark}\label{ODEfred}
Likewise, one might start with the Evans formulation \eqref{matsys},
written as a first order system, and try to precondition by factoring out the
expected asymptotic behavior,  to obtain essentially system \eqref{Ueq}.
That is, the operations of preconditioning (viewing the Fredholm
formulation as an analogous preconditioning step of factoring
out expected spatially-asypmptotic behavior) and reduction to ODE
essentially commute, at least in this simple case.
\end{remark}

\subsubsection{Stability index computation}\label{multistab}
Following the approach of Section \ref{alt}, computation
of the multi-dimensional stability index can be carried out in
the same way, with no additional complications.
For, exactly as in \eqref{2.55}
of the one-dimensional case (but using \eqref{2.42}),
one has the formula
\begin{equation}
\cF_{2}^{\bullet}(0)
= {\det}_{2,L^2(\bbR\times\Omega;dx)}
\big(I_{L^2(\bbR\times\Omega;dx)}-\cK_0-P_0\big) \, {\det}_{P_0
L^2(\bbR\times\Omega;dx)}(P_0\cK_1P_0),
\lb{multformula}
\end{equation}
with $\cK_0 = \cK(0)= -\ol{{H^{(0)}_{\Omega, \rm p}}^{-1}uv}$ and
$\cK_1 = \cK^{\bullet}(z)\big|_{z=0}= - \ol{{H^{(0)}_{\Omega, \rm p}}^{-2}uv}$. 
(One can remove the closure symbols in the last two expressions since all operators involved are bounded.) 

The second inner-product-type factor is straightforward to evaluate,
requiring only approximation of the eigenfunctions
of the operator $I-\cK_0$ and its adjoint corresponding
to the zero eigenvalue,
which in many cases are known from the outset. We recall that in the
present case, the  eigenfunction $\Phi$
of the operator $I-\cK_0$ is $dU/dx_1$;
here $U$ is the standing wave, see Section 4.1.1. The
eigenfunction $\widetilde\Phi$
of the operator $(I-\cK_0)^\ast$ may be deduced by the fact that the original
differential operator $L$ is self-adjoint; specifically,
$\Big(I-(H^{(0)}_{\Omega, \rm p}\big)^{-1}(V-V_\infty)\Big) \Phi=0$
implies, by self-adjointness of $\Big(I+\big(H^{(0)}_{\Omega, \rm p}\big)^{-1}\Big)$ 
and $(V-V_\infty)$,
that $\widetilde\Phi=H^{(0)}_{\Omega, \rm p}\Phi$ is indeed the required
 eigenfunction:
$\Big(I+\big(H^{(0)}_{\Omega, \rm p}\big)^{-1}(V-V_\infty)\Big)^* \widetilde\Phi=0$.

The first factor in \eqref{multformula} on the other hand
is the characteristic
determinant of a rank-one perturbation at $z=0$, so can be approximated as in Section \ref{alt}
using Galerkin approximation/semi-separable reduction
by a finite dimensional determinant. Precisely, combining the steps of Sections
\ref{fred} and \ref{alt},
one reduces the computation at the $J$-th Galerkin level to the evaluation of a
$2(N_J+1)\times 2(N_J+1)$
determinant, obtained by solving a $2(N_J+1)\times 2(N_J+1)$ ODE system 
\begin{equation}
\f{d \wti \cU(x_1)}{dx_1}=\wti A(x_1)\wti\cU(x_1),
\quad \wti \cU(x_1)\in \bbC^{2(N_J+1)\times 2(N_J+1)},
\end{equation}
where, similarly to \eqref{Ueq},
\begin{align}
\begin{split}
\wti A & =
\begin{pmatrix}
\wti G_1^\top  \wti F_1& \wti G_1^\top  \wti F_2\\[1mm]
-\wti G_2^\top  \wti F_1& -\wti G_2^\top  \wti F_2\\
\end{pmatrix},   \label{tA} \\
\wti F_k& =\begin{pmatrix} F_k\, & \, \Phi\end{pmatrix},
\quad
\wti G_k =\begin{pmatrix} G_k\, & \, H^{(0)}_{\Omega, \rm p} \Phi\end{pmatrix},  
\end{split} 
\end{align}
and $\Phi$ is the eigenfunction of the operator
$\cK_0$ corresponding to the eigenvalue $1$, that is,
$\cK_0\Phi=\Phi$.
Making the change of coordinates $\wti \cU= \wti D \wti \cV$,
$ \wti D=\diag \{D, 1\}$, see \eqref{defDchV}, we obtain a system
$d \wti \cV(x_1)/dx_1 = \wti A_b(x_1) \wti \cV(x_1)$
with bounded coefficient matrix $\wti A_b$ that can be numerically
solved by standard techniques used to compute the Evans function.

By comparison, if one follows the existing Galerkin methods, working
with an approximate truncated system at level $J$, one must face
the difficulty that zero eigenvalues for the exact system perturb
to small but in general nonzero eigenvalues of the approximate system,
making difficult a straightforward numerical computation without
further analytical preparations.
On the other hand, the usual analytic preparations (see
\cite{AGJ90},\cite{PW92}) involve solving variational equations
about the zero-energy eigenfunction and also projecting out the zero
eigenmodes from other modes to obtain a well-conditioned basis.
These do not appear to be real obstructions to the computation, but
are at least complications.  Perhaps for this reason, to our knowledge
no such computation has so far been carried out, or even proposed
in full detail.\footnote{However, see the interesting analysis
\cite{OS07} in the somewhat different spatially periodic case.}

The formulation of the above multi-dimensional stability index
algorithm we thus view as a useful practical contribution of 
the present work, and its numerical realization as an important
direction for further investigation.

\subsubsection{Numerical conditioning}\label{conditioning}
Last, we examine the question of numerical conditioning.
By Theorem \ref{equivalence}, one way to compute the approximate
Fredholm determinants $\cF_{2,J}$ is to carry out a standard
Evans function computation as in \cite{LPSS00}.
However, from a numerical perspective, this
might be missing the point.
For, note that the principal, constant-coefficient diagonal,
part of the coefficient matrix of \eqref{matsys} has entries
$|j|^2$ leading to spatial growth rates $\pm |j|$ of
order up to $J$. Computing the Evans function 
thus involves solution of an ODE that becomes
infinitely stiff as $J\to \infty$.

We suggest as a possible alternative, discretizing the kernel
$\cK_J$ in variables $x_1$, $x_1'$ and directly evaluating the
determinant of the resulting $MN_J\times MN_J$ matrix,
where $M$ is the number of mesh points in the $x_1$ ($x_1'$)
discretization required to give a desired error bound.
Noting that $M$ is essentially dimension-independent for
simple first-order quadrature (since matrix norm $|\cW_J|$
is bounded, by Parseval's identity, while the first derivative
of $\diag \Big\{ \frac{e^{ -(V_\infty + |j|^2-z)^{1/2}|x_1|}}
{(V_\infty + |j|^2-z)^{1/2}} \Big\}$ is of order one),
we see that there should be a break-even point at which
the cost of order
$(MN_J)^3$ of evaluating the discretized determinant
should be better than the cost $\sim \tilde N(J) N_J^3$ of
evaluating the Evans function, or,
equivalently,
\begin{equation}
\tilde N(J)\ge M^3,
\end{equation}
where $\tilde N(J)$ denotes the number of mesh
points required to evaluate the $2N_J\times 2N_J$ Evans ODE
to the same tolerance.  For a first-order A-stable scheme,
$\tilde N(J)\sim M J^2$ (note: this is determined by truncation error,
which is estimated as proportional to second derivative of the solution
as the square of
the norm of the largest eigenvalues $\pm J$), yielding break-even at
$J\sim M$,
where typical values of $M$ are of order
$\sim 100$ \cite{HZ06}.

Though hardly conclusive, this rough calculation suggests at least
that direct Fredholm computation is worthy of further study; we
recall that $J\sim 100$ is the order studied in \cite{LPSS00}.
Alternatively, one might compute the Evans function not by
shooting, but by continuation-type algorithms
as suggested by Sandstede \cite{Sa99}, viewing
the eigenvalue equation as a two-point boundary-value problem,
avoiding stiffness by another route;
however, so far as we know, such a scheme has not
yet been implemented.
See \cite{HZ06} for further discussion of this approach.

\subsection{Functions with radial limits}\label{rad}
Finally, consider standing-wave solutions $U$ of \eqref{multeq}
on the whole space $\bbR^d$, possessing radial limits in the following sense:
Introduce spherical coordinates $x=(r,\omega)$, $r>0$, $\omega\in S^{d-1}$, and let
$U(R,\omega)$, be the restriction of $U$ to the  sphere of radius $R$. Then considered as a
function of the angle $\omega\in S^{d-1}$, $U(R,\dott)$, has an
$L^1(S^{d-1};d\omega_{d-1})$-limit as $R\to +\infty$.

\begin{remark}\label{rdecoupled}
Assuming the hypotheses of Subsection \ref{cyl}, there exist radially
symmetric solutions $U(|\dott|)$, where $U$ is the
solution of the corresponding one-dimensional problem with the second
derivative replaced by the spherical Laplacian.
Linearizing about $U$ and expanding in spherical harmonics,
one obtains a decoupled family of one-dimensional eigenvalue problems,
similarly as in Remark \ref{decoupled},
each of which possess a well-defined Evans function and stability
index. In this case, there is no zero-eigenvalue at the zeroth harmonic
(constant function), but there is a zero-eigenvalue of order $d$ at the
level of the first harmonic, with associated eigenfunctions
$d U(|\cdot|)/dx_j$, $j=1, \dots, d$, corresponding to
translation-invariance of the underlying equations.
These may be treated similarly as in Remark \ref{decoupled}.
(However, we note that this involves an Evans function on the half-line
$[0,+\infty)$, which involves some modifications and will be analyzed elsewhere.)
\end{remark}

\begin{remark}\label{rcoupled}
In the general case, Galerkin approximation in spherical harmonics
yields a finite-dimensional system for which an Evans function and
stability can again be defined, similarly as in Remark \ref{coupled}.
\end{remark}

\subsubsection{Fredholm determinant version: radial case}\label{rfred}
In the simplest situation that $U$ has a single limit
as $|x|\to + \infty$, the operator $H^{(0)}_{\Omega, \rm p}$ is again constant-coefficient,
and the procedure of Section \ref{fred} leads again to expansion
of $H^{(0)}_{\Omega, \rm p}$ in a countable sum of operators with semi-separable integral kernels corresponding to the restrictions to different spherical harmonics.
In the general case, we may proceed instead by Galerkin approximation
as described in Remark \ref{exp}.

\begin{remark}\label{integral}
The common feature of the problems discussed is the presence of
a single unbounded spatial dimension (axial for cylindrical
case, radial for the radial case), along which the semi-separable
reduction is performed. In principle, one could treat still more general
problems by truncation/disretization of a continuous Fourier integral.
In this setting, the reference to a concrete object in the form
of a Fredholm determinant might become still more useful for numerical
validation/conditioning,
as compared to ad hoc constructions like those in
Remarks \ref{decoupled} and \ref{coupled}.
However, it is not clear that there would be a computation advantage
to doing so.
\end{remark}

\subsection{General operators}\label{general}
We recall from \cite{GLM07}, that it was necessary for general first-order
operators to relate the Evans function and a $2$-modified determinant
already in the one-dimensional case, since the Birman--Schwinger
kernel is for first-order operators only Hilbert--Schmidt
(indeed, this is one of the key insights of \cite{GLM07}).
Likewise, for more general operators involving a first-order component,
in particular those arising
in the study of stability of viscous shock solutions of
hyperbolic--parabolic conservation laws appearing in continuum mechanics
\cite{Zu03}, it is necessary in dimensions $d=2$ and $3$ to relate
the Evans function to a higher modified Fredholm determinant, since the
Birman--Schwinger kernel is no longer Hilbert--Schmidt. 
Flow in a cylindrical duct has been studied for viscous shock
and detonation waves in \cite{TZ06}, \cite{TZ07}.

\bigskip
\noindent {\bf Acknowledgments.}
We would like to thank Vita Borovik and Alin Pogan for helpful discussions.
We are indebted to the organizers of the workshop ``Stability
Criteria for Multi-Dimensional Waves and Patterns'', at the American
Institute of
Mathematics (AIM) in Palo Alto (California/USA), May 16-20, 2005, for
providing a most exciting environment which led to this collaboration.
Yuri Latushkin was partially supported by the
 Research Board and the Research Council
 of the University of
 Missouri and by the EU Marie Curie
 ''Transfer of Knowledge'' program.



\begin{thebibliography}{99}
\bi{AGJ90} J.\ Alexander, R.\ Gardner, and C.\ Jones,
{\em A topological invariant arising in the stability analysis of
travelling waves}, J. reine angew. Math. {\bf 410}, 167--212 (1990).
%
\bi{Br91} C.\ Brislawn, {\it Traceable integral kernels on countably generated 
measure spaces}, Pac.\ J.\ Math. {\bf 150}, 229--240 (1991). 
%
\bi{DN06} J.\ Deng and S.\  Nii,
{\it Infinite-dimensional Evans function theory for
elliptic boundary value problems},
J. Diff. Eq. {\bf 225}, 57--89 (2006).
%
\bi{DN08} J.\ Deng and S.\  Nii,
{\it An infinite-dimensional Evans function theory for
elliptic eigenvalue problems in a channel},
J. Diff. Eq. {\bf 244}, 753--765 (2008).
%
\bi{DS88} N.~Dunford and J.~T.~Schwartz, {\it Linear Operators Part II:
Spectral Theory}, Interscience, New York, 1988.
%
\bi{Ev72} J.\ W.\ Evans, {\it Nerve axon equations. I. Linear
approximations}, Indiana Univ. Math. J. {\bf 21}, 877--885 (1972).
%
\bi{Ev72a} J.\ W.\ Evans, {\it Nerve axon equations. II. Stability at
rest}, Indiana Univ. Math. J. {\bf 22}, 75--90 (1972).  Errata: Indiana
Univ. Math. J. {\bf 25}, 301 (1976).
%
\bi{Ev72b} J.\ W.\ Evans, {\it Nerve axon equations. III: Stability of the
nerve impulse}, Indiana Univ. Math. J. {\bf 22}, 577--593 (1972).
Errata:  Indiana Univ. Math. J. {\bf 25}, 301 (1976).
%
\bi{Ev75} J.\ W.\ Evans, {\it Nerve axon equations. IV. The stable and
unstable impulse}, Indiana Univ. Math. J. {\bf 24}, 1169--1190 (1975).
%
\bi{GH03} F.\ Gesztesy and H.\ Holden, {\it Soliton Equations and
Their Algebro-Geometric Solutions. Vol. I: $(1+1)$-Dimensional
Continuous Models},  Cambridge Studies in Advanced Mathematics,
Vol.\ 79, Cambridge Univ. Press, Cambridge, 2003.
%
\bi{GLM07} F.\ Gesztesy, Y.\ Latushkin, and K.\ A.\ Makarov, {\it
Evans Functions, Jost Functions, and Fredholm Determinants},
Arch. Rat. Mech. Anal., {\bf 186},  361--421 (2007).
%
\bibitem{GLMZ05} F.\,Gesztesy, Y.\,Latushkin, M.\,Mitrea and M. Zinchenko,
{\it Nonselfadjoint operators, infinite determinants, and some
applications}, Russ. J. Math. Phys. {\bf 12}, 443--471 (2005).
%
\bi{GM04} F.\ Gesztesy and K.\ A.\ Makarov, {\it (Modified) Fredholm
Determinants for Operators with Matrix-Valued Semi-Separable Integral
Kernels Revisited}, Integral Equations and Operator Theory {\bf 47}, 457--497
(2003). (See also Erratum {\bf 48}, 425--426 (2004) and the corrected
electronic only version in {\bf 48}, 561--602  (2004).)
%
\bi{GGK96} I.\ Gohberg, S.\ Goldberg, and N.\ Krupnik, {\it Traces
and determinants of linear operators}, Integr. Eqns. Oper. Theory {\bf
26}, 136--187 (1996).
%
\bi{GGK97} I.\ Gohberg, S.\ Goldberg, and N.\ Krupnik, {\it
Hilbert--Carleman and regularized determinants for linear operators},
Integr. Equ. Oper. Theory {\bf 27}, 10--47 (1997).
%
\bi{GGK00} I.\ Gohberg, S.\ Goldberg, and N.\ Krupnik, {\it Traces and
Determinants for Linear Operators}, Operator Theory: Advances and
Applications, Vol.\ 116, Birkh\"auser, Basel, 2000.
%
\bi{GK69} I.\ Gohberg and M.\ G.\ Krein, {\it Introduction to the Theory of
Linear Nonselfadjoint Operators}, Translations of Mathematical Monographs,
Vol.\ 18, Amer. Math. Soc., Providence, RI, 1969.
%
\bi{HZ06} J. \ Humpherys and K.\ Zumbrun, {\it An efficient shooting
algorithm for {E}vans function calculations in large systems},
Phys. D {\bf 220}, 116--126 (2006).
%
\bi{He81} D.\ Henry, {\it Geometric Theory of Semilinear Parabolic Equations},
Lecture Notes in Math., Vol.\ {\bf 840}, Springer, Berlin, 1981.
%
\bi{JP51} R.\ Jost and A.\ Pais, {\it On the scattering of a particle by a
static potential}, Phys. Rev. {\bf 82}, 840--851 (1951).
%
\bi{Ka80} T.\ Kato, {\it Perturbation Theory for Linear Operators},
corr. printing of the 2nd ed., Springer, Berlin, 1980.
%
\bi{LP07} Y.\ Latushkin and A.\ Pogan,
The Dichotomy Theorem for evolution bi-families, J. Diff. Eq., to appear.
%
\bi{LPSS00}
G.\ J.\ Lord, D.\ Peterhof, B.\ Sandstede, and A.\ Scheel.
{\it Numerical computation of solitary waves in infinite cylindrical 
domains}, SIAM J. Numer. Anal. {\bf 37}, 1420--1454 (2000).
%
\bi{LLMP05} L.\ Lorenzi, A.\ Lunardi, G.\ Metafune, D.\ Pallara, {\it Analytic
Semigroups and Reaction--Diffusion Equations}, Internet Seminar, 2004--2005;
available at  \\
http://www.math.unipr.it/~lunardi/LectureNotes/I-Sem2005.pdf.
%
\bi{Lu95} A.\ Lunardi, {\it Analytic semigroups and optimal
regularity in parabolic
problems}, Birkh\"auser, Basel, 1995.
%
\bi{Ne80} R.\ G.\ Newton, {\it Inverse scattering. I. One dimension},
J. Math. Phys. {\bf 21}, 493--505 (1980).
%
\bi{N07} J.\ Niesen, {\it Evans function calculations for a
two-dimensional system}, presented talk,
SIAM Conference on Applications of Dynamical Systems, Snowbird, UT, USA,
May 2007.
%
\bi{OS07} M.\ Oh and B.\ Sandstede,
{\it Evans function for periodic waves in infinite cylindrical domain},
preprint, 2007.
%
\bi{Pa92} C.\ V.\ Pao, {\it Nonlinear Parabolic and Elliptic Equations},
Plenum, New York, 1992.
%
\bi{PW92} R.\ L.\ Pego and M.\ I.\ Weinstein, {\it Eigenvalues, and
instabilities of solitary waves}, Phil. Trans. Roy. Soc. London {\bf
A 340}, 47--94 (1992).
%
\bi{PZ04} R.\ Plaza and K.\ Zumbrun, {\it An Evans function approach
to spectral stability of small-amplitude shock profiles}, Discrete
Cont. Dyn. Syst. B {\bf  10}, 885--924 (2004).
%
\bibitem{RS78} M.\ Reed and B.\ Simon,
{\em Methods of Modern Mathematical Physics. IV:
Analysis of Operators,} Academic Press, New York, 1978.
%
\bi{Ro84} F.\ Rothe, {\it Global Solutions of Reaction-Diffusion
Systems}, Lecture Notes in Math. {\bf 1072}, Springer, Berlin, 1984.
%
\bi{Sa99} B.\ Sandstede, private communication.
%
\bi{Sa02} B.\ Sandstede, {\it Stability of travelling waves}, in {\it
Handbook of
Dynamical Systems}, Vol.\ 2, B.\ Fiedler (ed.), Elsevier, Amsterdam,
2002, pp.\ 983--1055.
%
\bi{Si77} B.\ Simon, {\it Notes on infinite determinants of Hilbert space
operators}, Adv. Math. {\bf 24}, 244--273 (1977).
%
\bi{Si00} B.\ Simon, {\it Resonances in one dimension and Fredholm
determinants}, J. Funct. Anal. {\bf 178}, 396--420 (2000).
%
\bi{Si05} B.\ Simon, {\it Trace Ideals and Their Applications}, 2nd ed.,
Mathematical Surveys and Monographs, Vol.\ 120, Amer. Math. Soc.,
Providence, RI, 2005.
%
\bi{Sm94} J.\ Smoller, {\it Shock waves and reaction-diffusion
equations}, 2nd ed., Springer, New York, 1994.
%
\bi{TZ06} B.\ Texier and K.\ Zumbrun,
{\it Galloping instability of viscous shock waves}, Physica D, to appear, 
preprint 2006, available at http://arxiv.org/abs/math.AP/0609331.
%
\bi{TZ07} B.\ Texier and K.\ Zumbrun,
{\it Hopf bifurcation of viscous shock waves in
compressible gas- and magnetohydrodynamics},
Arch. Rational Mech. Anal., to appear, DOI 10.1007/s00205-008-0112-x, 
preprint 2006, available at http://arxiv.org/abs/math.AP/0612044.
%
\bi{Tr95} H.\ Triebel, {\it Interpolation Theory, Function Spaces, Differential 
Operators}, J.\ A.\ Barth, Heidelberg, 1995. 
%
\bi{VD99} R.\ Vein and P.\ Dale, {\it Determinants and Their
Applications in Mathematical Physics}, Springer, New York, 1999.
%
\bi{Wl87} J.\ Wloka, {\it Partial Differential Equations}, Cambridge University 
Press, Cambridge, 1987.
%
\bi{Wo52} F.\ Wolf, {\it Analytic perturbation of operators in Banach
spaces}, Math. Ann. {\bf 124}, 317--333 (1952).
%
\bi{Ya92} D.\ R.\ Yafaev, {\em Mathematical Scattering Theory},
Transl. Math. Monographs,  Vol.\ 105, Amer. Math. Soc.,
Providence, RI, 1992.
%
\bi{Zu03} K.\ Zumbrun, {\it Stability of large-amplitude shock
waves of compressible {N}avier--{S}tokes equations},
In {\em Handbook of mathematical fluid dynamics. Vol.  III}, pages
311--533. North-Holland, Amsterdam, 2004.
With an appendix by Helge Kristian Jenssen and Gregory Lyng.
%
\bi{ZH98} K.\ Zumbrun and P.\ Howard, {\it Pointwise semigroup methods and
stability of viscous shock waves}, Indiana Univ. Math. J. {\bf 47},
937--992 (1998). Errata: Indiana Univ. Math. J. {\bf 51}, 1017--1021 (2002).
\end{thebibliography}
\end{document}